**Example 2.**

|   | M |   |   |   |   |   |   |   |   |   |   |   |   |   |   |   |   |   |   |   |   |
|---|---|---|---|---|---|---|---|---|---|---|---|---|---|---|---|---|---|---|---|---|---|
|   | 1 | 2 | 3 | 4 | 5 | 6 | 7 | 8 | 9 | 10 | 11 | 12 | 13 | 14 | 15 | 16 | 17 | 18 | 19 | 20 |   |
| 1 | ∞ | 88 | 72 | 97 | 14 | 38 | 9 | 59 | 39 | 46 | 52 | 50 | 29 | 17 | 48 | 65 | 2 | 72 | 86 | 65 | 1 |
| 2 | 80 | ∞ | 28 | 72 | 67 | 18 | 99 | 20 | 58 | 66 | 24 | 76 | 32 | 45 | 11 | 62 | 54 | 62 | 25 | 45 | 2 |
| 3 | 55 | 40 | ∞ | 82 | 34 | 26 | 73 | 76 | 97 | 17 | 13 | 33 | 23 | 94 | 76 | 87 | 56 | 32 | 74 | 81 | 3 |
| 4 | 29 | 76 | 77 | ∞ | 92 | 63 | 94 | 88 | 87 | 18 | 38 | 59 | 94 | 62 | 33 | 18 | 9 | 67 | 93 | 31 | 4 |
| 5 | 17 | 32 | 60 | 80 | ∞ | 49 | 21 | 64 | 77 | 54 | 41 | 18 | 91 | 4 | 35 | 29 | 10 | 19 | 99 | 35 | 5 |
| 6 | 21 | 19 | 31 | 75 | 19 | ∞ | 98 | 94 | 72 | 40 | 30 | 43 | 48 | 18 | 94 | 82 | 69 | 70 | 22 | 71 | 6 |
| 7 | 1 | 99 | 82 | 68 | 4 | 35 | ∞ | 74 | 28 | 44 | 81 | 59 | 24 | 33 | 9 | 10 | 91 | 95 | 58 | 89 | 7 |
| 8 | 1 | 49 | 71 | 49 | 52 | 30 | 52 | ∞ | 59 | 92 | 25 | 48 | 56 | 65 | 34 | 59 | 24 | 78 | 67 | 70 | 8 |
| 9 | 92 | 99 | 49 | 6 | 63 | 79 | 66 | 19 | ∞ | 63 | 90 | 24 | 92 | 21 | 4 | 43 | 77 | 68 | 84 | 66 | 9 |
| 10 | 46 | 65 | 2 | 4 | 80 | 54 | 92 | 92 | 72 | ∞ | 34 | 10 | 86 | 63 | 63 | 40 | 73 | 99 | 85 | 5 | 10 |
| 11 | 11 | 14 | 81 | 18 | 91 | 37 | 46 | 51 | 32 | 23 | ∞ | 90 | 30 | 85 | 18 | 66 | 54 | 85 | 31 | 19 | 11 |
| 12 | 43 | 1 | 52 | 64 | 6 | 39 | 79 | 89 | 39 | 44 | 77 | ∞ | 78 | 42 | 47 | 62 | 68 | 65 | 25 | 7 | 12 |
| 13 | 19 | 82 | 93 | 76 | 20 | 80 | 15 | 81 | 15 | 87 | 45 | 67 | ∞ | 54 | 9 | 92 | 16 | 67 | 7 | 4 | 13 |
| 14 | 84 | 73 | 89 | 90 | 56 | 96 | 31 | 52 | 28 | 26 | 23 | 54 | 19 | ∞ | 91 | 37 | 6 | 95 | 59 | 26 | 14 |
| 15 | 84 | 53 | 65 | 42 | 65 | 54 | 62 | 81 | 90 | 80 | 98 | 52 | 59 | 44 | ∞ | 18 | 79 | 39 | 50 | 91 | 15 |
| 16 | 13 | 5 | 77 | 60 | 81 | 5 | 88 | 17 | 58 | 48 | 62 | 12 | 59 | 20 | 48 | ∞ | 69 | 61 | 20 | 57 | 16 |
| 17 | 46 | 48 | 25 | 59 | 8 | 83 | 83 | 24 | 28 | 1 | 19 | 75 | 17 | 28 | 82 | 75 | ∞ | 71 | 31 | 6 | 17 |
| 18 | 31 | 50 | 84 | 98 | 26 | 80 | 67 | 51 | 83 | 80 | 82 | 90 | 42 | 9 | 3 | 26 | 68 | ∞ | 51 | 41 | 18 |
| 19 | 53 | 89 | 6 | 44 | 58 | 48 | 26 | 17 | 64 | 88 | 63 | 63 | 87 | 61 | 42 | 57 | 32 | 59 | ∞ | 68 | 19 |
| 20 | 43 | 6 | 73 | 51 | 49 | 52 | 14 | 56 | 35 | 18 | 24 | 80 | 65 | 47 | 6 | 55 | 91 | 85 | 84 | ∞ | 20 |
|   | 1 | 2 | 3 | 4 | 5 | 6 | 7 | 8 | 9 | 10 | 11 | 12 | 13 | 14 | 15 | 16 | 17 | 18 | 19 | 20 |   |

**MIN(M)**

|    | 1  | 2  | 3  | 4  | 5  | 6  | 7  | 8  | 9  | 10 | 11 | 12 | 13 | 14 | 15 | 16 | 17 | 18 | 19 | 20 |
|----|----|----|----|----|----|----|----|----|----|----|----|----|----|----|----|----|----|----|----|----|
| 1  | 17 | 7  | 5  | 14 | 13 | 6  | 9  | 10 | 15 | 12 | 11 | 8  | 16 | 20 | 3  | 18 | 19 | 2  | 4  | 1  |
| 2  | 15 | 6  | 8  | 11 | 19 | 3  | 13 | 14 | 20 | 17 | 9  | 16 | 18 | 10 | 5  | 4  | 12 | 1  | 7  | 2  |
| 3  | 11 | 10 | 13 | 6  | 18 | 12 | 5  | 2  | 1  | 17 | 7  | 19 | 8  | 15 | 20 | 4  | 16 | 14 | 9  | 3  |
| 4  | 17 | 10 | 16 | 1  | 20 | 15 | 11 | 12 | 14 | 6  | 18 | 2  | 3  | 9  | 8  | 5  | 7  | 13 | 19 | 4  |
| 5  | 14 | 17 | 1  | 12 | 18 | 7  | 16 | 2  | 15 | 20 | 11 | 6  | 10 | 3  | 8  | 9  | 4  | 13 | 19 | 5  |
| 6  | 14 | 2  | 5  | 1  | 19 | 11 | 3  | 10 | 12 | 13 | 17 | 18 | 20 | 9  | 4  | 16 | 8  | 15 | 7  | 6  |
| 7  | 1  | 5  | 9  | 16 | 13 | 9  | 6  | 14 | 10 | 19 | 12 | 4  | 8  | 11 | 3  | 20 | 17 | 18 | 2  | 7  |
| 8  | 1  | 17 | 11 | 6  | 15 | 12 | 2  | 4  | 5  | 7  | 13 | 9  | 16 | 14 | 19 | 20 | 3  | 18 | 10 | 8  |
| 9  | 15 | 4  | 8  | 14 | 12 | 16 | 3  | 5  | 10 | 7  | 20 | 18 | 17 | 6  | 19 | 11 | 1  | 13 | 2  | 9  |
| 10 | 3  | 4  | 20 | 12 | 11 | 16 | 1  | 6  | 14 | 15 | 2  | 9  | 17 | 5  | 19 | 13 | 7  | 8  | 18 | 10 |
| 11 | 1  | 2  | 4  | 15 | 20 | 10 | 13 | 19 | 9  | 6  | 7  | 8  | 17 | 16 | 3  | 14 | 18 | 12 | 5  | 11 |
| 12 | 2  | 5  | 20 | 19 | 6  | 9  | 14 | 1  | 10 | 15 | 3  | 16 | 4  | 17 | 18 | 11 | 13 | 7  | 8  | 12 |
| 13 | 20 | 19 | 15 | 7  | 9  | 17 | 1  | 5  | 11 | 14 | 12 | 18 | 4  | 6  | 8  | 2  | 10 | 16 | 3  | 13 |
| 14 | 17 | 13 | 11 | 10 | 20 | 9  | 7  | 16 | 8  | 12 | 5  | 19 | 2  | 1  | 3  | 4  | 15 | 18 | 6  | 14 |
| 15 | 16 | 18 | 4  | 14 | 12 | 2  | 6  | 13 | 19 | 7  | 3  | 5  | 17 | 10 | 8  | 1  | 9  | 20 | 11 | 15 |
| 16 | 2  | 6  | 12 | 1  | 8  | 14 | 19 | 10 | 15 | 20 | 9  | 13 | 4  | 18 | 11 | 17 | 3  | 5  | 7  | 16 |
| 17 | 10 | 5  | 20 | 13 | 11 | 8  | 3  | 9  | 14 | 19 | 1  | 2  | 4  | 18 | 12 | 16 | 15 | 6  | 7  | 17 |
| 18 | 15 | 14 | 5  | 16 | 1  | 20 | 13 | 2  | 8  | 19 | 7  | 17 | 6  | 10 | 11 | 9  | 3  | 12 | 4  | 18 |
| 19 | 3  | 8  | 7  | 17 | 15 | 4  | 6  | 1  | 16 | 5  | 18 | 14 | 11 | 12 | 9  | 20 | 13 | 10 | 2  | 19 |
| 20 | 2  | 15 | 7  | 10 | 11 | 9  | 1  | 14 | 5  | 4  | 6  | 16 | 8  | 13 | 3  | 12 | 19 | 18 | 17 | 20 |
|    | 1  | 2  | 3  | 4  | 5  | 6  | 7  | 8  | 9  | 10 | 11 | 12 | 13 | 14 | 15 | 16 | 17 | 18 | 19 | 20 |

D = (1 2 3 4 5 6 7 8 9 10 11 12 13 14 15 16 17 18 19 20).

From M, the values of the arcs of D are

      88    28    82    92    49    98    74    59    63    34
(1    2    3    4    5    6    7    8    9    10

      90    78    54    91    18    69    71    51    68    43
   11    12    13    14    15    16    17    18    19    20).

We now use MIN(M) to obtain the smallest value in each row of M.
MIN(M)(1, 1) = 17: d(1, 17) = 2. MIN(M)(2, 1) = 15: d(2, 15) = 11.
MIN(M)(3, 1) = 11: d(3, 11) = 13. MIN(M)(4, 1) = 17: d(4, 17) = 9.
MIN(M)(5, 1) = 14: d(5, 14) = 18. MIN(M)(6, 1) = 14: d(6, 14) = 18.
MIN(M)(7, 1) = 1: d(7, 1) = 1. MIN(M)(8, 1) = 1: d(8, 1) = 1.
MIN(M)(9, 1) = 15: d(9, 15) = 4. MIN(M)(10, 1) = 3: d(10, 3) = 2.
MIN(M)(11, 1) = 1: d(11, 1) = 11. MIN(M)(12, 1) = 2: d(12, 2) = 1.
MIN(M)(13, 1) = 20: d(13, 20) = 4. MIN(M)(14, 1) = 17: d(14, 17) = 6.
MIN(M)(15, 1) = 16: d(15, 16) = 18. MIN(M)(16, 1) = 2: d(16, 2) = 5.
MIN(M)(17, 1) = 10: d(17, 10) = 1. MIN(M)(18, 1) = 15: d(18, 15) = 3.
MIN(M)(19, 1) = 3: d(19, 3) = 6. MIN(M)(20, 1) = 2: d(20, 2) = 6.
For j = 1, 2, ... , 19, we now obtain d(j, MIN(M)(j,1)) - d(j, j+1);
also, d(20, 2) - d(20, 1). We then replace the values given the arcs of D
by the corresponding differences we've just constructed. We now have

$$D = \begin{pmatrix} 1 & 2 & 3 & 4 & 5 & 6 & 7 & 8 & 9 & 10 \\ -86 & -17 & -69 & -83 & -45 & -80 & -73 & -58 & -59 & -32 \end{pmatrix}$$

$$\begin{pmatrix} 11 & 12 & 13 & 14 & 15 & 16 & 17 & 18 & 19 & 20 \\ -79 & -77 & -50 & -85 & 0 & -64 & -70 & -48 & -62 & -37 \end{pmatrix}$$

We next construct a permutation consisting of at most two cycles which when applied to D will reduce the total sum of its arcs in M. To start the construction, we choose an arc where the absolute value of its difference is greatest. We therefore start our first cycle with 1. Before going on, we introduce an economical way to write permutations, their inverses and their differences as we proceed from derangement to derangement.

We first write two ordered lists of the numbers from 1 through n.
We then first enter the number corresponding to $D_i(a)$ underneath a in the first list. We then write the number preceding a underneath a in the second list. Simply speaking, if have

$$\begin{matrix} a \\ b \end{matrix}$$

occurring in a column in the first list, we have

$$\begin{matrix} b \\ a \end{matrix}$$

occurring in a column in the second list. If we don't know the predecessor, p, of the first point a in a cycle of a permutation, it will always occur in the second row of the permutation when we reach the arc $(p, D_i(p))$ where $D_i(p) = a$. Finally, we place the value of DIFF(a) above a in the second list. When we follow this practice for D, we must write all points of both D and $D^{-1}$ because the permutation D is itself a derangement. However, when constructing further permutations, we need only change those columns in which the point a is contained in a permutation.

We now rewrite D and $D^{-1}$ in the form mentioned earlier. For simplicity, we call it the row form of D.

1 2 3 4 5 6 7 8  9 10 11 12 13 14 15 16 17 18 19 20

2 3 4 5 6 7 8 9 10 11 12 13 14 15 16 17 18 19 20  1

$1^{-86}$ $2^{-17}$ $3^{-69}$ $4^{-83}$ $5^{-45}$ $6^{-80}$ $7^{-73}$ $8^{-58}$ $9^{-59}$ $10^{-32}$ $11^{-79}$ $12^{-77}$ $13^{-50}$ $14^{-85}$
20    1     2     3     4     5     6     7     8      9      10     11     12     13

$15^{0}$ $16^{-64}$ $17^{-70}$ $18^{-48}$ $19^{-62}$ $20^{-37}$
14  15    16    17    18    19

As we construct our lists, we keep track of the size of the values of DIFF obtained. We start off with DIFF(1). If an when we find a smaller value of DIFF, say at point a, we replace DIFF(1) by DIFF(a), replacing the preceding value of DIFF. Thus, when we are finished we have obtained both the point at which the smallest value of DIFF occurs as well as the value of DIFF there.

Using this procedure, we see that the smallest value of DIFF occurs when a = 1. DIFF(1) = -86.

The rule that we adopted in PHASE 1 was that we try the first [log n] + 1 entries MIN(M)(i, j) in row i (j = 1, 2, ... , [log n] + 1). log(20) = 2.996. Thus, in this case, we will pick the first three possibilities {MIN(M)(i, j) | j = 1,2,3}. There are two more things of which we must be aware. The first is that if an arc (a, b) of $D_i$ is chosen when constructing a permutation ,we will obtain (a, a) for the corresponding arc of the permutation. If (a, b) comes from MIN(M)(a, 1), then we must choose MIN(M)(a, 2). Otherwise, we will choose MIN(M)(a, 1). Secondly, in testing for cycles, the last arc of a permutation cycle is often chosen without knowing the arc of the value matrix from which it is obtained. Occasionally, it may occur that the arc chosen actually comes from a loop. i.e., an arc of the form (a, a). Thus, the current derangement would lead to a permutation which is not a derangement. To avoid this occurrence, we must check the terminal arc of any possibility that we have not gotten directly from an arc obtained from MIN(M). To do this, given the permutation arc, say (a, b), we obtain from the first two rows of the row form, (a, $D_i(b)$). If $D_i(b) \neq a$, then we have no difficulty. Otherwise, we go on to test the next possible terminal arc of a cycle.

TRIAL 1.

MIN(M)(1,1) = 17. d(1, 17) = 2. In order to construct a corresponding arc of the cycle, we must obtain $D^{-1}(17)$ = 16. Thus, the first arc of the first cycle is (1, 16). We now continue with the vertex 16. MIN(M)(16, 1) = 2. Thus, the second arc of M chosen is (16, 2) where d(16, 2) = 5. $D^{-1}(2)$ = 1. Thus, the first cycle we've obtained is (1 16). Adding up the values corresponding to 1 and 16, respectively, if we applied the 2-cycle (1 16) to D, we would reduce the total sum of arcs of D in M by 86 + 64 = 150.

TRIAL 2.

MIN(M)(1,2) = 7. d(1, 7) = 9. Thus, instead of - 86, we now have - 79. $D^{-1}(7)$ = 6. Thus, the first arc of the first cycle of a permutation is
(1 6).

We will now proceed with the construction of the first cycle. We must always keep in mind that the path we construct must always have a negative sum of values. Otherwise, the permutation won't reduce the total value of the derangement to which it will be applied.

        (1  7)    (1 6)  - 79
        (6  14)   (6 13)  - 80
        (13  20)  (13 19) - 50
        (19  3)   (19  2) - 62
        (2  15)   (2 14) - 17
        (14  17)  (14  16) - 85
        (16  2)   (16 1) - 64.

The cycle (1 6 13 19 2 14 16) reduces the sum of the values of the arcs of D in M by 437.

TRIAL 3.

MIN(M)(1, 3) = 5. d(1, 5) = 14. Thus, the value of the arc of D out of 1 is reduced by 74.

$$\begin{array}{ll}(1\ 5) & (1\ 4) - 74 \\ (4\ 17) & (4\ 16) - 83 \\ (16\ 2) & (16\ 1) - 64.\end{array}$$

The cycle (1 4 16) reduces the sum of the values of the arcs of D by 221. It follows that our best choice is the second one. Let $s_1$ = (1 6 13 19 2 14 16). Next, we give the row form of the new permutation, $D_1$, obtained by applying $s_1$ to D. Asterisks have been placed above each point of D which has been moved.

```
* *         *                    * *      *          *
1  2  3  4  5  6  7  8   9 10 11 12 13 14 15 16 17 18 19 20

7 15  4  5  6 14  8  9 10 11 12 13 20 17 16  2 18 19  3  1
```

$$\begin{array}{llllllllllllllll}1 & 2 & 3^{-69} & 4^{-83} & 5^{-45} & 6 & 7^{-73} & 8^{-58} & 9^{-59} & 10^{-32} & 11^{-79} & 12^{-77} & 13 & 14 & 15^0 & 16 \\ 20 & 16 & 19 & 3 & 4 & 5 & 1 & 7 & 8 & 9 & 10 & 11 & 12 & 6 & 2 & 15\end{array}$$

$$\begin{array}{llll}17^{-70} & 18^{-48} & 19 & 20^{-37} \\ 14 & 17 & 18 & 13\end{array}$$

We now must obtain the DIFF values for points under asterisks.

MIN(M)(1, 1) = 17 $\Rightarrow$ DIFF(1) = d(1, 17) - d(1, 17) = 0.
MIN(M)(2, 1) = 15 $\Rightarrow$ DIFF(2) = d(2, 15) - d(2, 15) = 0.
MIN(M)(6, 1) = 14 $\Rightarrow$ DIFF(6) = d(6, 14) - d(6, 14) = 0.
MIN(M)(13, 1) = 20 $\Rightarrow$ DIFF(13) = d(13, 20) - d(13, 20) = 0.
MIN(M)(14, 1) = 17 $\Rightarrow$ DIFF(14) = d(14, 17) - d(14, 17) = 0.
MIN(M)(16, 1) = 2 $\Rightarrow$ DIFF(16) = d(16, 2) - d(16, 2) = 0.
MIN(M)(19, 1) = 3 $\Rightarrow$ DIFF(19) = d(19, 3) - d(19, 3) = 0.

Thus, we can complete the row form above.

```
     1  2  3  4  5   6  7  8   9 10 11 12 13 14 15 16 17 18 19 20
D_1
     7 15  4  5  6  14  8  9  10 11 12 13 20 17 16  2 18 19  3  1
```

$$\begin{array}{llllllllllllll}1^0 & 2^0 & 3^{-69} & 4^{-83} & 5^{-45} & 6^0 & 7^{-73} & 8^{-58} & 9^{-59} & 10^{-32} & 11^{-79} & 12^{-77} & 13^0 & 14^0\end{array}$$

$D_1^{-1}$

$$\begin{array}{llllllllllllll}20 & 16 & 19 & 3 & 4 & 5 & 1 & 7 & 8 & 9 & 10 & 11 & 12 & 6\end{array}$$

$$\begin{array}{llllll}15^0 & 16^0 & 17^{-70} & 18^{-48} & 19^0 & 20^{-37}\end{array}$$

$$\begin{array}{llllll}2 & 15 & 14 & 17 & 18 & 13\end{array}$$

The smallest value of DIFF occurs when a = 4.

**TRIAL 1.**
**MIN(M)(4, 1) = 17.**
    (4, 17) → (4, 14)
    (14, 17) → (14, 14)
(14, 17) is an arc of $D_1$. We next try obtain MIN(M)(14, 2) = 13.
    (14, 13) → (14, 12)
    (12, 2) → (12, 16 )
    (16, 2) → (16, 16)
(16, 2) is an arc of $D_1$. MIN(M)(16, 2) = 6.
         (16, 6) → (16, 5)
         (5, 14) → (5, 6)
         (6, 14) → (6, 6)
(6, 14) is an arc of $D_1$.
         (6, 2) → (6, 16).
We thus have a path P = [4, 14, 12, 16, 5, 6, 16].
Let $s_{11}$ = (4 14 12 16 5 6).
We check to see the arc of $D_1$ associated with (6 4). Going to the first two rows of ROW FORM, we obtain (6 5). In those cases in which we didn't encounter an arc of $D_1$ in constructing an arc of $s_{11}$, we can use the values associated with the points of the cycle in ROW FORM.
Thus, $s_{11} = (4^{-83}\ 14\ 12^{-77}\ 16\ 5^{-45})$.
We now find the values of 14, 16 and 6. For each vertex a, we have the arc in M, obtained from MIN(M)(a, 2). We then subtract the arc in $D_1$ with the same initial vertex.
d(14, 13) - d(14, 17) = 19 - 6 = 13 .
d(16, 6) - d(16, 2) = 5 - 5 = 0.
d(6, 5) - d(6, 14) = 19 - 18 = 1
Thus, $s_{11} = (4^{-83}\ 14^{13}\ 12^{-77}\ 16^{0}\ 5^{-45}\ 6^{1})$.
$s_{11}$ has the value -191. Thus, when applied to $D_1$, it reduces the value of $D_1$ by 191.
We now consider the permutations obtained by including all partial paths of $s_{11}$ which begin with 4. Let
$s_{111}$ = (4 14), $s_{112}$ = (4 14 12), $s_{113}$ = (4 14 12 16),
$s_{114}$ = (4 14 12 16 5).
Since 4 is mapped into 5 by $D_1$, in order to obtain the value of each of these permutations, we have to obtain the respective values of (14, 5), (12, 5), (16, 5). Since (5,5) is a loop, we need not consider $s_{114}$.
d(14, 5) = 56, d(12, 5) = 6, d(16, 5) = 81.
We next subtract the value of each arc of form (a, $D_1$(a))
(a = 14, 12, 5) from the respective values obtained above to obtain the value of each arc in its respective permutation.
d(14, 5) - d(14, 17) = 56 - 6 = 50.
d(12, 5) - d(14, 17) = 6 - 6 = 0.
d(16, 5) - d(16, 2) = 81 - 5 = 76.

Thus, $s_{111} = (4^{-83}\ 14^{50})$, $s_{112} = (4^{-83}\ 14^{13}\ 12^{0})$,
$s_{113} = (4^{-83}\ 14^{13}\ 12^{-77}\ 16^{76})$. Thus, the respective values of these permutations are -33, -70 and -71.

We now consider the permutation obtained by using the path P to construct a permutation of two cycles.

$s_{12} = (4\ 14\ 12)(16\ 5\ 6)$.

We already have the values of all of the points of $s_{12}$ except for the arc (6 16). We now obtain the value associated with 6 in $s_{12}$.

$d(6, 16) \rightarrow d(6, 15) = 94$

$d(6, 15) - d(6, 14) = 94 - 18 = 76$.

$s_{12} = (4^{-83}\ 14^{13}\ 12^{0})(16^{0}\ 5^{-45}\ 6^{76})$.

The first cycle of $s_{12}$ is $s_{112}$. The value of the second cycle is non-negative so we needn't include it.

TRIAL 2. MIN(M)(4, 2) = 10.

$$(4, 10) \rightarrow (4, 9)$$
$$(9, 15) \rightarrow (9, 2)$$
$$(2, 15) \rightarrow (2, 2)$$

Thus, (2, 15) is an arc of $D_1$. MIN(M)(2,2) = 6

$$(2, 6) \rightarrow (2, 5)$$
$$(5, 14) \rightarrow (5, 6)$$
$$(6, 14) \rightarrow (6, 6).$$

(6, 14) is an arc of $D_1$. MIN(M)(6, 2) = 2.

$$(6, 2) \rightarrow (6, 16)$$
$$(16, 2) \rightarrow (16, 16).$$

(16, 2) is an arc of $D_1$.

$$(16, 6) \rightarrow (16, 5).$$

Let P = [4, 9, 2, 5, 6 16, 5].

$s_{21} = (4\ 9\ 2\ 5\ 6\ 16)$.

We must obtain the values associated with the initial vertices 2, 6, 16. The other arcs were all obtained using terminal vertices in the first column of MIN(M). Thus, we can read them off directly from ROW FORM. We must compute $d(16, D_1(4)) = d(16, 5)$ rather than using $d(16, 6)$.

$d(2, 6) - d(2, 15) = 18 - 11 = 7$.

$d(6, 2) - d(6, 14) = 19 - 18 = 1$.

$d(16, 5) - d(16, 2) = 81 - 5 = 76$.

Thus, $s_{21} = (4^{-83}\ 9^{-59}\ 2^{7}\ 5^{-45}\ 6^{1}\ 16^{76})$. The value of $s_{21}$ is -103.

$s_{211} = (4\ 9)$, $s_{212} = (4\ 9\ 2)$, $s_{213} = (4\ 9\ 2\ 5)$, $s_{214} = (4\ 9\ 2\ 5\ 6)$.

We must thus obtain d(9, 5), d(2, 5), d(6, 5). $s_{213}$ need not be considered since the arc (5, 4) of the permutation maps into the loop (5, 5).

$d(9, 5) = 63$, $d(2, 5) = 67$, $d(6, 5) = 19$.

The respective values in $s_{211}$, $s_{212}$, $s_{214}$ are

$d(9, 5) - d(9, 10) = 63 - 63 = 0$, $d(2, 5) - d(2, 15) = 67 - 11 = 56$,
$d(6, 5) - d(6, 14) = 19 - 18 = 1$.

$s_{211} = (4^{-83} \, 9^0), (4^{-83} \, 9^{-59} \, 2^{56})$, $s_{214} = (4^{-83} \, 9^{-59} \, 2^7 \, 5^{-45} \, 6^1)$.

The respective values of these permutations are -83, -86, -179.

We now construct $s_{22}$ from P.

$s_{22} = (4\ 9\ 2)(5\ 6\ 16)$. We have already obtained the values for all vertices except 16 previously.

$d(16, D_1(5)) = d(16, 6)$. $d(16, 6) - d(16, 2) = 5 - 5 = 0$.

Thus, we have

$s_{22} = (4^{-93} \, 9^{-59} \, 2^{56})(5^{-45} \, 6^1 \, 16^0)$.

The value of $s_{22}$ is -140.

TRIAL 3. MIN(M)(4, 3) = 16.

$$(4, 16) \rightarrow (4, 15)$$
$$(15, 16) \rightarrow (15, 15)$$

(15, 16) is an arc of $D_1$.

$$(15, 18) \rightarrow (15, 17)$$
$$(17, 10) \rightarrow (17, 9)$$
$$(9, 15) \rightarrow (9, 2)$$
$$(2, 15) \rightarrow (2, 2)$$

(2, 15) is an arc of $D_1$.

$$(2, 6) \rightarrow (2, 5)$$
$$(5, 14) \rightarrow (5, 6)$$
$$(6, 14) \rightarrow (6, 6)$$

(6, 14) is an arc of $D_1$.

$$(6, 2) \rightarrow (6, 16)$$
$$(16, 2) \rightarrow (16, 16)$$

(16, 2) is an arc of $D_1$.

$$(16, 6) \rightarrow (16, 5)$$

P = [4, 15, 17, 9, 2, 5, 6, 16, 5].

15, 2, 6 and 16 were the initial vertices of arcs in $D_1$. Thus, we had to use arcs obtained from the second column of MIN(M). We can immediately read off the values of the remaining vertices from ROW FORM..

$s_{31} = (4^{-83} \, 15 \, 17^{-70} \, 9^{-59} \, 2 \, 5^{-45} \, 6 \, 16)$.

$d(15, 18) - d(15, 16) = 39 - 18 = 21$.
$d(2, 6) - d(2, 15) = 18 - 11 = 7$.
$d(6, 2) - d(6, 14) = 19 - 18 = 1$.
$d(16, 5) - d(16, 2) = 81 - 5 = 76$.

We thus obtain $s_{31} = (4^{-83} \, 15^{21} \, 17^{-70} \, 9^{-59} \, 2^7 \, 5^{-45} \, 6^1 \, 16^{76})$. The value of $s_{31}$ is -152.

$s_{311} = (4\ 15)$. $s_{312} = (4\ 15\ 17)$. $s_{313} = (4\ 15\ 17\ 9)$. $s_{314} = (4\ 15\ 17\ 9\ 2)$.

Since (5, 4) maps into (5, 5), we need not consider (4 15 17 9 2 5).

$s_{315} = (4\ 15\ 17\ 9\ 2\ 5\ 6)$.

$d(15, 5) - d(15, 16) = 65 - 18 = 47$.
$d(17, 5) - d(17, 18) = 8 - 71 = -63$.
$d(9, 5) - d(9, 10) = 63 - 63 = 0$.
$d(2, 5) - d(2, 15) = 67 - 11 = 56$.

d(6, 5) - d(6, 14) = 19 - 18 = 1.
d(16, 5) - d(16, 2) = 81 - 5 = 76.
Thus, $s_{311} = (4^{-83} 15^{47})$, $s_{312} = (4^{-73}15^{21}17^{-73})$, $s_{313} = (4^{-73}15^{21}17^{-70}9^0)$,
$s_{314} = (4^{-83}15^{21}17^{-70}9^{-59}2^{56})$, $s_{315} = (4^{-83}15^{21}17^{-70}9^{-59}2^7 5^{-45}6^1)$.
The respective values of these permutations are: -36, -125, -122, -135, -228.
From P, $s_{32} = (4\ 15\ 17\ 9\ 2)(5\ 6\ 16)$. The first cycle is $s_{314}$. The second was obtained in trial 2 as the second factor of $s_{22}$. It follows that the value of $s_{32}$ is (-135) + (-44) = -179.
After checking all cycles obtained, $s_{315}$ has the smallest value: -228.
Thus, $s_2 = (4\ 15\ 17\ 9\ 2\ 5\ 6)$. In order to obtain $D_2$ in row form, we need only change the points 4, 15, 17, 9, 2, 5, 6. In constructing permutations, we have obtained the changes in the arcs of $D_1$. However, if we wish, we can use the row form of $D_1$ to obtain new columns in $D_2$.
Our new arcs are: (4 16), (15 18), (17 10), (9 15), (2 6), (5 14), (6 5).

```
        1  2  3  4  5  6  7  8  9  10 11 12 13 14 15 16 17 18 19 20
D_2
        7  6  4 16 14  5  8  9 15  11 12 13 20 17 18  2 10 19  3  1

        0 -7 -69 -9  0 -1 -73 -58  0 -32 -79 -77  0  0 -21  0  0 -48  0 -37
        1  2  3  4  5  6  7  8  9  10 11 12 13 14 15 16 17 18 19 20

       20 16 19  3  6  2  1  7  8  17 10 11 12  5  9  4 14 15 18 13
```

d(2, 15) - d(2, 6) = 11 - 18 = -7.
d(4, 17) - d(4, 16) = 9 - 18 = -9.
d(5, 14) - d(5, 14) = 0.
d(6, 14) - d(6, 5) = 18 - 19 = -1
d(9, 15) - d(9, 15) = 0.
d(15, 16) - d(15, 18) = 18 - 39 = -21.
d(17, 10) - d(17, 10) = 0.
We choose 11 as an initial vertex for our trials.
TRIAL 1. MIN(M)(11, 1) = 1.
$\quad$ (11, 1) → (11, 20)
$\quad$ (20, 2) → (20, 16)
$\quad$ (16, 2) → (16, 16)
(16, 2) is an arc of $D_2$.
$\quad$ (16, 6) → (16, 2)
$\quad$ (2, 15) → (2, 9)
$\quad$ (9, 15) → (9, 9).
(9, 15) is an arc of $D_2$.
$\quad$ (9, 4) → (9, 3).
$\quad$ (3, 11) → (3, 10)
$\quad$ (10, 3) → (10, 19)
$\quad$ (19, 3) → (19, 19)

(19, 3) is an arc of $D_2$.

$$(19, 8) \to (19, 7)$$
$$(7, 1) \to (7, 20)$$

P = [11, 20, 16, 2 9 3 10 19, 7, 20].

$s_{11} = (11^{-79}\ 20^{-37}\ 16\ 2\ 9\ 3^{-69}\ 10^{-32}\ 19\ 7)$

d(16, 6) - d(16, 2) = 5 - 5 = 0.
d(2, 6) - d(2, 15) = 18 - 11 = 7.
d(9, 4) - d(9, 15) = 6 - 4 = 2.
d(19, 8) - d(19, 3) = 17 - 6 = 11.
d(7, 12) - d(7, 1) = 59 - 1 = 58

$s_{11} = (11^{-79}20^{-37}16^0 2^7 9^2 3^{-69}10^{-32}19^{11}7^{58})$.

The value of $s_{11}$ is -139.

$s_{111} = (11\ 20)$, $s_{112} = (11\ 20\ 16)$, $s_{113} = (11\ 20\ 16\ 2)$, $s_{114} = (11\ 20\ 16\ 2\ 9)$,
$s_{115} = (11\ 20\ 16\ 2\ 9\ 3)$, $s_{116} = (11\ 20\ 16\ 2\ 9\ 3\ 10)$,
$s_{117} = (11\ 20\ 16\ 2\ 9\ 3\ 10\ 19)$.

d(20, 12) - d(20, 1) = 80 - 43 = 37.
d(16, 12) - d(16, 2) = 12 - 5 = 7.
d(2, 12) - d(2, 6) = 76 - 18 = 58.
d(9, 12) - d(9, 15) = 24 - 4 = 20.
d(3, 12) - d(3, 4) = 33 - 82 = -49.
d(10, 12) - d(10, 11) = 10 - 34 = -24.
d(19, 12) - d(19, 3) = 63 - 6 = 57.

$s_{111} = (11^{-79}20^{37})$, $s_{112} = (11^{-79}20^{-37}16^7)$, $s_{113} = (11^{-79}20^{-37}16^0 2^{58})$,
$s_{114} = (11^{-79}20^{-37}16^0 2^7 9^{20})$, $s_{115} = (11^{-79}20^{-37}16^0 2^7 9^2 3^{-49})$.
$s_{116} = (11^{-79}\ 20^{-37}16^0 2^7 9^2 3^{-69}10^{-24})$, $s_{117} = (11^{-79}20^{-37}16^0 2^7 9^2 3^{-69}10^{-32}19^{57})$.

The respective values of the above permutations are: -42, -109, -58, -89, -156, -200, -151.

$s_{12} = (20^{-37}16^0 2^7 9^2 3^{-69}10^{-32}19^{11}7^0)$.

The value of $s_{12}$ is -118.

TRIAL 2.  MIN(M)(11, 2) = 2.

$$(11, 2) \to (11, 16)$$

From this point on, the path is the same as that of P in trial 1 up until we reach 7. Thus, let us continue on from 7 allowing P' to be

P' = [11, 16, 2, 9, 3, 10, 19, 7]

$$(7, 1) \to (7, 20)$$
$$(20, 2) \to (20, 16).$$

Thus, P = [11, 16, 2, 9, 3, 10, 19, 7, 20, 16]

$s_{21} = (11\ 16\ 2\ 9\ 3\ 10\ 19\ 7\ 20)$

d(11, 2) - d(11, 12) = 14 - 90 = -76.
d(20, 12) - d(20, 1) = 80 - 43 = 37.

$s_{21} = (11^{-76}16^0 2^7 9^2 3^{-69}10^{-32}19^{11}7^0 20^{37})$.

The value of $s_{21}$ is -120.

$s_{211} = (11^{-76}16)$, $s_{212} = (11^{-76}16^0 2)$, $s_{213} = (11^{-76}16^0 2^7 9)$,

$s_{214} = (11^{-76}16^{0}2^{7}9^{2}3)$, $s_{215} = (11^{-76}16^{0}2^{7}9^{2}3^{-69}10)$, $s_{216} = (11^{-76}16^{0}2^{7}9^{2}3^{-69}10^{-32}19)$,

$s_{217} = (11^{-76}16^{0}2^{7}9^{2}3^{-69}10^{-32}19^{11}7)$

$d(16, 12) - d(16, 2) = 12 - 5 = 7$.
$d(2, 12) - d(2, 6) = 76 - 18 = 58$.
$d(9, 12) - d(9, 15) = 24 - 4 = 20$.
$d(3, 12) - d(3, 4) = 33 - 82 = -49$.
$d(10, 12) - d(10, 11) = 10 - 34 = -24$.
$d(19, 12) - d(19, 3) = 80 - 6 = 74$.
$d(7, 12) - d(7, 8) = 59 - 74 = -15$.

$s_{211} = (11^{-76}16^{7})$, $s_{212} = (11^{-76}16^{0}2^{58})$, $s_{213} = (11^{-76}16^{0}2^{7}9^{20})$,

$s_{214} = (11^{-76}16^{0}2^{7}9^{2}3^{-49})$, $s_{215} = (11^{-76}16^{0}2^{7}9^{2}3^{-69}10^{-24})$,

$s_{216} = (11^{-76}16^{0}2^{7}9^{2}3^{-69}10^{-32}19^{74})$, $s_{217} = (11^{-76}16^{0}2^{7}9^{2}3^{-69}10^{-32}19^{11}7^{-15})$.

The respective values of the above permutations are: -69, -18, -49, -116, -172, -94, -172.

$s_{22} = (16^{0}2^{7}9^{2}3^{-69}10^{-32}19^{11}7^{0}20)$

$d(20, 4) - d(20, 1) = 51 - 43 = 8$.

$s_{22} = (16^{0}2^{7}9^{2}3^{-69}10^{-32}19^{11}7^{0}20^{8})$.

The value of $s_{22}$ is -73.

TRIAL 3.  MIN(M)(11, 3) = 4.

$(11, 4) \rightarrow (11, 3)$
$(3, 11) \rightarrow (3, 10)$

Using $s_{22}$, we can assume that P', a subpath of P, is
P' = [11, 3, 10, 19, 7, 20].
Continuing from 20, we obtain

$(20, 2) \rightarrow (20, 16)$
$(16, 2) \rightarrow (16, 16)$

(16, 2) is an arc of $D_{2}$.

$(16, 6) \rightarrow (16, 2)$
$(2, 15) \rightarrow (2, 9)$
$(9, 15) \rightarrow (9, 9)$

(9, 15) is an arc of $D_{2}$.

$(9, 4) \rightarrow (9, 3)$.

P = [11, 3, 10 19, 7, 20, 16, 2, 9, 3].

$s_{31} = (11\ 3\ 10\ 19\ 7\ 20\ 16\ 2\ 9)$.

$d(11, 4) - d(11, 12) = 18 - 90 = -72$.
$d(20, 2) - d(20, 1) = 6 - 43 = -37$.
$d(16, 6) - d(16, 2) = 5 - 5 = 0$.
$d(2, 15) - d(2, 6) = 11 - 18 = -7$.
$d(9, 12) - d(9, 15) = 24 - 4 = 20$.

$s_{31} = (11^{-72}3^{-69}10^{-32}19^{11}7^{0}20^{-37}16^{0}2^{-7}9^{20})$.

The value of $s_{31}$ is -186.

$d(3, 12) - d(3, 4) = 33 - 82 = -59$.
$d(10, 12) - d(10, 11) = 10 - 34 = -24$.

d(19, 12) - d(19, 3) = 63 - 6 = 57.
d(7, 12) - d(7, 8) = 59 - 74 = 15.
d(20, 12) - d(20, 1) = 80 - 43 = 37.
d(16, 12) - d(16, 2) = 12 - 5 = 7.
d(2, 12) - d(2, 6) = 76 - 18 = 58.

$s_{311} = (11^{-72}3^{-59})$, $s_{312} = (11^{-72}3^{-69}10^{-24})$, $s_{313} = (11^{-72}3^{-69}10^{-32}19^{57})$,

$s_{314} = (11^{-72}3^{-59}10^{-32}19^{11}7^{15})$, $s_{315} = (11^{-72}3^{-69}10^{-32}19^{11}7^{0}20^{37})$,

$s_{316} = (11^{-72}3^{-69}10^{-32}19^{11}7^{0}20^{-37}16^{7})$, $s_{317} = (11^{-72}3^{-69}10^{-32}19^{11}20^{-37}16^{0}2^{58})$.

The respective values of the above permutations are: -131, -173, -116, -147, -141, -192, -141.

d(9, 4) - d(9, 15) = 6 - 4 = 2.

$s_{32} = (3^{-69}10^{-32}19^{11}7^{0}20^{-37}16^{0}2^{-7}9^{2})$.

The value of $s_{32}$ is -132. It follows that our best possibility for $s_{3}$ is

$s_{116} = (11^{-79}20^{-37}16^{0}2^{7}9^{2}3^{-69}10^{-24})$.

(11, 20) → (11, 1), (20, 16) → (20, 2), (16, 2) → (16, 6),
(2, 9) → (2, 15), (9, 3) → (9, 4), (3, 10) → (3, 11), (10, 11) → (10, 12).

Our new arcs in $D_{3}$ are: (11, 1), (20, 2), (16, 6), (2, 15), (9, 4), (3,11) and (10, 12).

```
   *  *            *  *  *                    *              *
1  2  3  4  5  6  7  8  9  10 11 12 13 14 15 16 17 18 19 20

7 15 11 16 14  5  8  9  4  12  1 13 20 17 18  6 10 19  3  2
```

d(2, 15) - d(2, 15) = 0.
d(3, 11) - d(3, 11) = 0.
d(9, 15) - d(9, 4) = 4 - 6 = -2.
d(10, 3) - d(10, 12) = 2 - 10 = -8.
d(11, 1) - d(11, 1) = 0.
d(16, 2) - d(16, 6) = 5 - 5 = 0.
d(20, 2) - d(20, 2) = 0.

```
 0  0  0 -9  0  0 -73 -58 -2 -8  0 -77  0  0 -21  0  0 -48  0  0
 1  2  3  4  5  6   7   8  9 10 11  12 13 14  15 16 17  18 19 20
11 20 19  9  6 16   1   7  8 17  3  10 12  5   2  4 14  15 18 13
```

The smallest value of DIFF occurs out of 12.

TRIAL 1. MIN(M)(12, 1) = 2.

(12, 2) → (12, 20)
(20, 2) → (20, 20).

(20, 2) is an arc of $D_{3}$.

(20, 15) → (20, 2)
(2, 15) → (2, 2).

(2, 15) is an arc of $D_{3}$.

(2, 6) → (2, 16)
(16, 2) → (16, 20)

P = [12, 20, 2, 16, 20].

$s_{11} = (12\ 20\ 2\ 16)$

d(20, 2) - d(20, 15) = 0.
d(2, 6) - d(2, 15) = 18 - 11 = 7.
d(16, 2) - d(16, 13) = 59 - 5 = 54.
$s_{11} = (12^{-77} 20^{0} 2^{7} 16^{54})$.
The value of $s_{11}$ is -16.
$s_{12}$ = (20 2 16).
d(16, 2) - d(16, 2) = 0.
$s_{12} = (20^{0} 2^{7} 16^{0})$.
The value of $s_{12}$ is non-negative.
TRIAL 2. MIN(M)(12, 2) = 5.
d(12, 5) - d(12, 13) = 6 - 78 = -72.
$$(12, 5) \to (12, 6)$$
$$(6, 14) \to (6, 5)$$
$$(5, 14) \to (5, 5)$$
(5, 14) is an arc of $D_3$.
$$(5, 17) \to (5, 14)$$
$$(14, 17) \to (14, 14)$$
(14, 17) is an arc of $D_3$.
$$(14, 13) \to (14, 12)$$
P = [12, 6, 5, 14, 12]
$s_{21}$ = (12 6 5 14)
d(6, 14) - d(6, 5) = 18 -19 = -1.
d(5, 17) - d(5, 14) = 10 - 4 = 6.
d(14, 13) - d(14, 17) = 19 - 6 = 13.
$s_{21} = (12^{-72} 6^{-1} 5^{6} 14^{13})$.
The value of $s_{21}$ is -54.
TRIAL 3. MIN(M)(12, 3) = 20.
d(12, 20) - d(12, 13) = 7 - 78 = -71.
$$(12, 20) \to (12, 13)$$
$$(13, 20) \to (13, 13)$$
(13, 20) is an arc of $D_3$.
$$(13, 19) \to (13, 18)$$
$$(18, 15) \to (18, 2)$$
$$(2, 15) \to (2, 2)$$
(2, 15) is an arc of $D_3$.
$$(2, 6) \to (2, 16)$$
$$(16, 2) \to (16, 20)$$
$$(20, 2) \to (20, 20)$$
(20, 2) is an arc of $D_3$.
$$(20, 15) \to (20, 2)$$
P = [12, 13, 18, 2, 16, 20, 2].
$s_{31}$ = (12 13 18 2 16 20)
(12, 13) → (12, 20), (13, 18) → (13, 19), (18, 2) → (18, 15),
(2, 16) → (2, 6), (16, 20) → (16, 2), (20, 12) → (20, 13).

d(12, 20) - d(12, 13) = 7 - 78 = -71.
d(13, 19) - d(13, 20) = 7 - 4 = 3.
d(18, 15) - d(18, 19) = 3 - 51 = -48.
d(2, 6) - d(2, 15) = 18 - 11 = 7.
d(16, 2) - d(16, 6) = 5 - 5 = 0.
d(20, 13) - d(20, 2) = 65 - 6 = 59.
(13, 12) → (13, 13). Thus, $s_{311}$ need not be considered.
$s_{312} = (12^{-71}13^{3}18^{-9})$, $s_{313} = (12^{-71}13^{3}18^{-48}2^{21})$, $s_{314} = (12^{-71}13^{3}18^{-48}2^{7}16^{54})$.
The respective values of the above permutations are: -77, -95, -55.
d(20, 15) – d(20, 2) = 6 – 6 = 0.
$s_{32} = (12^{-71}13^{3}18^{-9})(2^{7}16^{0}20^{0})$.
The value of the second cycle of $s_{32}$ is non-negative, while the first cycle is precisely $s_{312}$.
The smallest-valued permutation is $s_{313}$. Thus, $s_{3}$ = (12 13 18 2).
(12, 13) → (12, 20), (13, 18) → (13, 19), (18, 2) → (18, 15),
(2, 12) → (2, 13). Thus, the new arcs of $D_4$ are: (12, 20), (13, 19),
(18, 15) and (2, 13).

|   |   |   |   |   |   |   |   |   |   |   | * |   |   |   | * | * |   |   |   | * |
|---|---|---|---|---|---|---|---|---|---|---|---|---|---|---|---|---|---|---|---|---|
|   | 1 | 2 | 3 | 4 | 5 | 6 | 7 | 8 | 9 | 10 | 11 | 12 | 13 | 14 | 15 | 16 | 17 | 18 | 19 | 20 |
| $D_4$ | 7 | 13 | 11 | 16 | 14 | 5 | 8 | 9 | 4 | 12 | 1 | 20 | 19 | 17 | 18 | 6 | 10 | 15 | 3 | 2 |

In order to create $D_4^{-1}$, we need only invert the four new arcs and use all of the remaining arcs from $D_3^{-1}$. Before going on, we obtain DIFF(a) for a = 2, 12, 13 and 18.
d(2, 15) – d(2, 13) = 11 – 32 = -21.
d(12, 2) – d(12, 20) = 1 – 7 = -6.
d(13, 20) – d(13, 19) = 4 – 7 = -3
d(18, 15) – d(18, 15) = 0.

| 0 | –21 | 0 | -9 | 0 | 0 | –73 | -58 | -2 | -8 | 0 | -6 | -3 | 0 | -21 | 0 | 0 | 0 | 0 | 0 |
|---|---|---|---|---|---|---|---|---|---|---|---|---|---|---|---|---|---|---|---|
| 1 | 2 | 3 | 4 | 5 | 6 | 7 | 8 | 9 | 10 | 11 | 12 | 13 | 14 | 15 | 16 | 17 | 18 | 19 | 20 |
| 11 | 20 | 19 | 9 | 6 | 16 | 1 | 7 | 8 | 17 | 3 | 10 | 2 | 5 | 18 | 4 | 14 | 15 | 13 | 12 |

The smallest DIFF value occurs when a = 7.
TRIAL 1.
$$(7, 1) \rightarrow (7, 11)$$
$$(11, 1) \rightarrow (11, 11)$$
(11, 1) is an arc of $D_4$.
$$(11, 2) \rightarrow (11, 20)$$
$$(20, 2) \rightarrow (20, 20)$$
(20, 2) is an arc of $D_4$.
$$(20, 15) \rightarrow (20, 18)$$
$$(18, 15) \rightarrow (18, 18)$$
(18, 15) is an arc of $D_4$.
$$(18, 14) \rightarrow (18, 5)$$

$$(5, 14) \to (5, 5)$$

(5, 14) is an arc of $D_4$.

$$(5, 17) \to (5, 14)$$
$$(14, 17) \to (14, 14)$$

(14, 17) is an arc of $D_4$.

$$(14, 13) \to (14, 2)$$
$$(2, 15) \to (2, 18)$$

P = [7, 11, 20, 18, 5, 14, 2, 18].

$s_{11}$ = (7 11 20 18 5 14 2)

d(7, 1) − d(7, 8) = 1 − 74 = -73
d(11, 2) − d(11, 1) = 14 − 11 = 3
d(20, 15) − d(20, 2) = 6 − 6 = 0.
d(18, 14) − d(18, 15) = 9 − 3 = 6.
d(5, 17) − d(5, 14) = 10 − 4 = 6.
d(14, 13) − d(14, 17) = 19 − 6 = 13.
d(2, 8) − d(2, 13) = 20 − 32 = -12.

$s_{11} = (7^{-73} 11^3 20^0 18^6 5^6 14^{13} 2^{-12})$.

The value of $s_{11}$ is −57.

d(11, 8) − d(11, 1) = 51 − 11 = 40.
d(20, 8) − d(20, 2) = 56 − 6 = 50.
d(18, 8) − d(18, 15) = 51 − 3 = 48.
d(5, 8) − d(5, 14) = 64 − 4 = 60.
d(14, 8) − d(14, 17) = 52 − 6 = 46.

$s_{114} = (7^{-73} 11^3 20^0 18^6 5^{60})$, $s_{115} = (7^{-73} 11^3 20^0 18^6 5^6 14^{46})$.

The respective values of the above permutations are: -33, -20, -21, -14, -12.

$s_{12}$ = (7 11 20)(18 5 14) = $(7^{-73} 11^3 20^{50})(18^6 5^6 14^{85})$.

d(14, 15) − d(14, 17) = 91 − 6 = 85.

The first cycle is $s_{112}$; the second cycle is non-negative.

TRIAL 2. MIN(M)(7, 2) = 5.

d(7, 5) − d(7,8) = 4 − 74 = -70

$$(7, 5) \to (7, 6)$$
$$(6, 14) \to (6, 5)$$

In Trial 1, 5 was followed by 14, 2 and 18. We continue with 18.

$$(18, 15) \to (18, 18)$$

(18, 15) is an arc of $D_4$.

$$(18, 14) \to (18, 5)$$

P = [7, 6, 5, 14, 2, 18, 5].

$s_{21}$ = (7 6 5 14 2 18).

d(6, 14) − d(6, 5) = 18 − 19 = -1.
d(5, 17) − d(5, 14) = 6.
d(14, 13) − d(14, 17) = 13.
d(2, 15) − d(2, 13) = 11 − 32 = -21.
d(18, 8) − d(18, 15) = 51 − 3 = 48.

$s_{21} = (7^{-70} 6^{-1} 5^6 14^{13} 2^{-21} 18^{48})$.

The value of $s_{21}$ is -25.
$d(6, 8) - d(6, 5) = 94 - 19 = 75$.
$d(5, 8) - d(5, 14) = 64 - 4 = 60$.
$d(14, 8) - d(14, 17) = 52 - 6 = 46$.
$d(2, 8) - d(2, 13) = 20 - 32 = -12$.
$s_{211} = (7^{-70} 6^{75})$, $s_{212} = (7^{-70} 6^0 5^{60})$, $s_{213} = (7^{-70} 6^0 5^6 14^{46})$,
$s_{214} = (7^{-70} 6^0 5^6 14^{13} 2^{-12})$.
The respective values of the above permutations are: 5, -10, -18, -63.
$d(2, 15) - d(2, 13) = 11 - 32 = -21$.
$d(18, 14) - d(18, 15) = 9 - 3 = 6$.
$s_{22} = (7^{-70} 6^{76})(5^6 14^{13} 2^{-21} 18^6)$.
Each cycle has a non-negative value.
TRIAL 3. MIN(M)(7, 3) = 9.
$d(7, 9) - d(7, 8) = 28 - 74 = -46$.
$$(7, 9) \to (7, 4)$$
$$(4, 17) \to (4, 14)$$
In Trial 1, we obtained 14, 2, 18, 5.
We continue with 5.
$$(5, 14) \to (5, 5)$$
(5, 14) is an arc of $D_4$.
$$(5, 17) \to (5, 14)$$
P = [7, 4, 14, 2, 18, 5, 14].
$s_{31} = (7\ 4\ 14\ 2\ 18\ 5)$.
$d(4, 17) - d(4, 16) = 9 - 18 = -9$.
$d(14, 13) - d(14, 17) = 19 - 6 = 13$.
$d(2, 15) - d(2, 13) = 11 - 32 = -21$.
$d(18, 14) - d(18, 15) = 9 - 3 = 6$.
$d(5, 8) - d(5, 14) = 64 - 4 = 60$.
$s_{31} = (7^{-46} 4^{-9} 14^{13} 2^{-21} 18^6 5^{60})$.
The value of $s_{31}$ is 3.
$d(4, 8) - d(4, 16) = 88 - 18 = 70$.
$d(14, 8) - d(14, 17) = 52 - 6 = 46$.
$d(2, 8) - d(2, 13) = 20 - 32 = -12$.
$d(18, 8) - d(18, 15) = 51 - 3 = 48$.
$s_{311} = (7^{-46} 4^{70})$, $s_{312} = (7^{-46} 4^{-9} 14^{46})$, $s_{313} = (7^{-46} 4^{-9} 14^{13} 2^{-12})$,
$s_{314} = (7^{-46} 4^{-9} 14^{13} 2^{-21} 18^{48})$.
The respective values of the above permutations are: 24, 1, -54, -15.
$d(5, 17) - d(5, 14) = 10 - 4 = 6$.
$s_{32} = (7^{-46} 4^{70})(14^{13} 2^{-21} 18^6 5^6)$
Both cycles are non-negative.
The smallest value of all permutations is $s_{214} = (7^{-70} 6^0 5^6 14^{13} 2^{-12})$.
We obtain the following new arcs for $D_5$: (7, 5), (6, 14), (5, 17), (14, 13), (2, 8).

```
    *           *  *  *                              *
 1  2  3  4  5  6  7  8  9  10  11  12  13  14  15  16  17  18  19  20

 7  8 11 16 17 14  5  9  4  12   1  20  19  13  18   6  10  15   3   2
```
$d(2, 15) - d(2, 8) = 11 - 20 = -9.$
$d(5, 14) - d(5, 17) = 4 - 10 = -6.$
$d(6, 14) - d(6,14) = 0.$
$d(7, 1) - d(7, 5) = 1 - 4 = -3.$
$d(14, 17) - d(14, 13) = 6 - 19 = -13.$
```
 0 -9  0 -9 -6  0 -3 -58 -2 -8  0 -6 -3 -13 -21  0  0  0  0  0
 1  2  3  4  5  6  7  8  9  10 11 12 13  14  15 16 17 18 19 20

11 20 19  9  7 16  1  2  8  17  3  10 14   6  18  4  5 15 13 12
```
We start with a = 8.

**TRIAL 1.**

$$d(8, 1) - d(8, 9) = 1 - 59 = -58.$$
$$(8, 1) \rightarrow (8, 11)$$
$$(11, 1) \rightarrow (11, 11)$$

(11, 1) is an arc of $D_5$.
$$(11, 2) \rightarrow (11, 20)$$
$$(20, 2) \rightarrow (20, 20)$$

(20, 2) is an arc of $D_5$.
$$(20, 15) \rightarrow (20, 18)$$
$$(18, 15) \rightarrow (18, 18)$$

(18, 15) is an arc of $D_5$.
$$(18, 14) \rightarrow (18, 6)$$
$$(6, 14) \rightarrow (6, 6)$$

(6, 14) is an arc of $D_5$
$$(6, 2) \rightarrow (6, 20)$$

P = [8, 11, 20, 18, 6, 20].

$s_{11}$ = (8 11 20 18 6)

$d(8, 1) - d(8, 9) = -58.$
$d(11, 2) - d(11, 1) = 14 - 11 = 3.$
$d(20, 15) - d(20, 2) = 6 - 6 = 0.$
$d(18, 14) - d(18, 15) = 9 - 3 = 6.$
$d(6, 9) - d(6, 14) = 72 - 18 = 54.$

$s_{11} = (8^{-58} 11^3 20^0 18^6 6^{54}).$

$d(11, 9) - d(11, 1) = 32 - 11 = 21.$
$d(20, 9) - d(20, 2) = 35 - 6 = 29.$
$d(18, 9) - d(18, 15) = 83 - 3 = 80.$

$s_{111} = (8^{-58} 11^{21}),\ s_{112} = (8^{-58} 11^3 20^{29}),\ s_{113} = (8^{-58} 11^3 20^0 18^{80}).$

The respective values of the above permutations are: -37, -27, 25.

$d(6, 2) - d(6, 14) = 19 - 18 = 1.$

$s_{12} = (8^{-58} 11^{21})(20^0 18^6 6^1).$

The first cycle is $s_{111}$. The second cycle has a non-negative value.

**TRIAL 2. MIN(M)(8, 2) = 17.**

$$d(8, 17) - d(8, 9) = 24 - 58 = -34.$$
$$(8, 17) \rightarrow (8, 5)$$
$$(5, 14) \rightarrow (5, 13)$$
$$(13, 20) \rightarrow (13, 2)$$
$$(2, 15) \rightarrow (2, 18)$$
$$(18, 15) \rightarrow (18, 18)$$

$(18, 15)$ is an arc of $D_5$.

$$(18, 14) \rightarrow (18, 6)$$
$$(6, 14) \rightarrow (6, 6)$$

$(6, 14)$ is an arc of $D_5$.

$$(6, 2) \rightarrow (6, 20)$$
$$(20, 2) \rightarrow (20, 20)$$

$(20, 2)$ is an arc of $D_5$.

$$(20, 15) \rightarrow (20, 18)$$

$P = [8, 5, 13, 2, 18, 6, 20, 18]$.
$s_{21} = (8\ 5\ 13\ 2\ 18\ 6\ 20)$.
$d(8, 17) - d(8, 9) = -34.$
$d(5, 14) - d(5, 17) = 4 - 10 = -6.$
$d(13, 20) - d(13, 19) = 4 - 7 = -3.$
$d(2, 15) - d(2, 8) = 11 - 20 = -9.$
$d(18, 14) - d(18, 15) = 9 - 3 = 6.$
$d(6, 2) - d(6, 14) = 19 - 18 = 1.$
$d(20, 15) - d(20, 2) = 6 - 6 = 0.$
$s_{21} = (8^{-34} 5^{-6} 13^{-3} 2^{-9} 18^6 6^1 20^{29}).$
The value of $s_{21}$ is $-16$.
$d(5, 9) - d(5, 17) = 77 - 10 = 67.$
$d(13, 9) - d(13, 19) = 15 - 7 = 8.$
$d(2, 9) - d(2, 8) = 58 - 20 = 38.$
$d(18, 9) - d(18, 15) = 83 - 3 = 80.$
$d(6, 9) - d(6, 14) = 72 - 18 = 54.$
$s_{211} = (8^{-34} 5^{67}),\ s_{212} = (8^{-34} 5^{-6} 13^8),\ s_{213} = (8^{-34} 5^{-6} 13^{-3} 2^{38}),$
$s_{214} = (8^{-34} 5^{-6} 13^{-3} 2^{-9} 18^{80}),\ s_{215} = (8^{-34} 5^{-6} 13^{-3} 2^{-9} 18^6 6^{54}).$
**The respective values of the above permutations are: 33, -32, -5, 28, 8.**

$s_{22} = (8\ 5\ 13\ 2)(18\ 6\ 20).$
**The first cycle is $s_{213}$.**

$d(20, 15) - d(20, 2) = 6 - 6 = 0.$ Thus, the second cycle is $(18^6 6^1 20^0)$.
**Thus, it has a non-negative value.**
**TRIAL 3. MIN(M)(8, 3) = 11.**
$$d(8, 11) - d(8, 9) = 25 - 59 = -34.$$
$$(8, 11) \rightarrow (8, 3)$$
$$(3, 11) \rightarrow (3, 3)$$

(3, 11) is an arc of $D_5$.

$$(3, 10) \rightarrow (3, 17)$$
$$(17, 10) \rightarrow (17, 17)$$

(17, 10) is an arc of $D_5$.

$$(17, 5) \rightarrow (17, 7)$$
$$(7, 1) \rightarrow (7, 11)$$
$$(11, 1) \rightarrow (11, 11)$$

(11, 1) is an arc of $D_5$.

$$(11, 2) \rightarrow (11, 20)$$
$$(20, 2) \rightarrow (20, 20)$$

(20, 2) is an arc of $D_5$.

$$(20, 15) \rightarrow (20, 18)$$
$$(18, 15) \rightarrow (18, 18)$$

(18, 15) is an arc of $D_5$.

$$(18, 14) \rightarrow (18, 6)$$
$$(6, 14) \rightarrow (6, 6)$$

(6, 14) is an arc of $D_5$.

$$(6, 2) \rightarrow (6, 20)$$
$$(20, 2) \rightarrow (20, 20)$$

(20, 2) is an arc of $D_5$.

$$(20, 15) \rightarrow (20, 18)$$

P = [8, 3, 17, 7, 11, 20, 18, 6, 20].

d(3, 10) − d(3, 11) = 17 − 13 = 4.
d(17, 5) − d(17, 10) = 8 − 1 = 7.
d(7, 1) − d(7, 5) = 1 − 4 = -3.
d(11, 2) − d(11, 1) = 14 − 11 = 3.
d(20, 15) − d(20, 2) = 6 − 6 = 0.
d(18, 14) − d(18, 15) = 9 − 3 = 6.
d(6, 9) − d(6, 14) = 72 − 18 = 54.

$s_{31} = (8^{-34} 3^4 17^7 7^{-3} 11^3 20^0 18^6 6^{54})$.

d(3, 9) − d(3,11) = 97 − 13 = 84.
d(17, 9) − d(17, 10) = 28 − 1 = 27.
d(7, 9) − d(7, 5) = 28 − 4 = 24.
d(11, 9) − d(11, 2) = 32 − 14 = 18.
d(20, 9) − d(20, 2) = 35 − 6 = 29.
d(18, 9) − d(18, 15) = 83 − 3 = 80.

$s_{311} = (8^{-34} 3^{84})$, $s_{312} = (8^{-34} 3^4 17^{27})$, $s_{313} = (8^{-34} 3^4 17^7 7^{24})$,

$s_{314} = (8^{-34} 3^4 17^7 7^{-3} 11^{18})$, $s_{315} = (8^{-34} 3^4 17^7 7^{-3} 11^3 20^{29})$,

$s_{316} = (8^{-34} 3^4 17^7 7^{-3} 11^3 20^0 18^{80})$.

$s_{32} = (8^{-34} 3^4 17^7 7^{-3} 11^{18})(20^0 18^6 6^1)$.

The first cycle is $s_{314}$. The second is non-negative.

The smallest negative cycle is $s_{111} = (8^{-58} 11^{21})$. Thus, $s_4 = (8\ 11)$.

| 1 | 2 | 3 | 4 | 5 | 6 | 7 | 8 | 9 | 10 | 11 | 12 | 13 | 14 | 15 | 16 | 17 | 18 | 19 | 20 |
|---|---|---|---|---|---|---|---|---|----|----|----|----|----|----|----|----|----|----|----|
| 7 | 8 | 11 | 16 | 17 | 14 | 5 | 1 | 4 | 12 | 9 | 20 | 19 | 13 | 18 | 6 | 10 | 15 | 3 | 2 |

| -7 | -9 | 0 | -9 | -6 | 0 | -3 | 0 | -2 | -8 | -21 | -6 | -3 | -13 | -21 | 0 | 0 | 0 | 0 | 0 |
|----|----|---|----|----|---|----|---|----|----|-----|----|----|-----|-----|---|---|---|---|---|
| 1 | 2 | 3 | 4 | 5 | 6 | 7 | 8 | 9 | 10 | 11 | 12 | 13 | 14 | 15 | 16 | 17 | 18 | 19 | 20 |

| 8 | 20 | 19 | 9 | 7 | 16 | 1 | 2 | 11 | 17 | 3 | 10 | 14 | 6 | 18 | 4 | 5 | 15 | 13 | 12 |
|---|----|----|---|---|----|---|---|----|----|---|----|----|---|----|---|---|----|----|----|

TRIAL 1. MIN(M)(15, 1) = 16.

$d(15, 16) - d(15, 18) = 18 - 39 = -21$.
$(15, 16) \rightarrow (15, 4)$
$(4, 17) \rightarrow (4, 5)$
$(5, 14) \rightarrow (5, 6)$
$(6, 14) \rightarrow (6, 6)$

$(6, 14)$ is an arc of $D_6$.

$(6, 2) \rightarrow (6, 20)$
$(20, 2) \rightarrow (20, 20)$

$(20, 2)$ is an arc of $D_6$.

$(20, 15) \rightarrow (20, 18)$
$(18, 15) \rightarrow (18, 18)$

$(18, 15)$ is an arc of $D_6$.

$(18, 14) \rightarrow (18, 6)$

P = [15, 4, 5, 6, 20, 18, 6].
$d(15, 16) - d(15, 18) = 18 - 39 = -21$.
$d(4, 17) - d(4, 16) = 9 - 18 = -9$
$d(5, 14) - d(5, 17) = 4 - 9 = -5$.
$d(6, 2) - d(6, 14) = 19 - 18 = 1$.
$d(20, 15) - d(20, 2) = 6 - 6 = 0$.
$s_{11} = (15\ 4\ 5\ 6\ 20\ 18)$.

We can't use $s_{11}$ since $D_6(18) = 18$. Thus, $D_6 s_{11}$ wouldn't be a derangement.

$d(4, 18) - d(4, 16) = 67 - 18 = 49$.
$d(5, 18) - d(5, 17) = 19 - 10 = 9$.
$d(6, 18) - d(6, 14) = 70 - 18 = 52$.
$d(20, 18) - d(20, 2) = 85 - 6 = 79$.

$s_{111} = (15^{-21} 4^{49})$, $s_{112} = (15^{-21} 4^{-9} 5^9)$, $s_{113} = (15^{-21} 4^{-9} 5^{-5} 6^{52})$,
$s_{114} = (15^{-21} 4^{-9} 5^{-5} 6^1 20^{79})$.

The respective values of the above permutations are: 28, -21, 17, 45.
$s_{22} = (15^{-21} 4^{-9} 5^9)(6^1 20^0 18^6)$.

The first cycle is $s_{112}$. The second cycle is non-negative.

TRIAL 2. MIN(M)(15, 2) = 18.

$(15, 18) - (15, 18) = 0$.

Thus, we can proceed no further in trials 2 and 3.

**The smallest negative permutation obtained is $s_{112} = (15^{-21}4^{-9}5^9)$.**
**Thus, $s_6 = (15^{-21}4^{-9}5^9)$. The new arcs of $D_7$ are: (15, 16), (4, 17), (5, 18).**

```
         *  *                                      *
  1  2  3  4  5  6  7  8  9 10 11 12 13 14 15 16 17 18 19 20

  7 20 11 17 18 14  5  1  4 12  9 20 19 13 16  6 10 15  3  2
```
$d(15, 16) - d(15, 16) = 0$.
$d(4, 17) - d(4, 17) = 0$.
$d(5, 14) - d(5, 18) = 4 - 19 = -15$.

```
 -7 -9  0  0 -15  0 -3  0 -2 -8 -21 -6 -3 -13  0  0  0  0  0  0
  1  2  3  4  5  6  7  8  9 10  11 12 13 14 15 16 17 18 19 20

  8 20 19  9  7 16  1  2 11 17  3 10 14  6 18 15  4  5 13 12
```
We let a = 11.
TRIAL 1. MIN(M)(11, 1) = 1.
$\qquad\qquad\qquad (11, 1) \rightarrow (11, 8)$
$\qquad\qquad\qquad (8, 1) \rightarrow (8, 8)$
(8, 1) is an arc of $D_7$.
$\qquad\qquad\qquad (8, 17) \rightarrow (8, 10)$
$\qquad\qquad\qquad (10, 3) \rightarrow (10, 19)$
$\qquad\qquad\qquad (19, 3) \rightarrow (19, 19)$
(19, 3) is an arc of $D_7$.
$\qquad\qquad\qquad (19, 8) \rightarrow (19, 2)$
$\qquad\qquad\qquad (2, 15) \rightarrow (2, 18)$
$\qquad\qquad\qquad (18, 15) \rightarrow (18, 18)$
(18, 15) is an arc of $D_7$.
$\qquad\qquad\qquad (18, 14) \rightarrow (18, 6)$
$\qquad\qquad\qquad (6, 14) \rightarrow (6, 6)$
(6, 14) is an arc of $D_7$.
$\qquad\qquad\qquad (6, 2) \rightarrow (6, 20)$
$\qquad\qquad\qquad (20, 2) \rightarrow (20, 20)$
(20, 2) is an arc of $D_7$.
$\qquad\qquad\qquad (20, 15) \rightarrow (20, 18)$
P = [11, 8, 10, 19, 2, 18. 6, 20, 18]
$d(8, 17) - d(8, 1) = 24 - 1 = 23$.
The path [11, 8, 10] has a positive value $-21 + 23 = 2$.
We need proceed no further in Trial 1. The only possible permutation is
(11 8).
$d(11, 1) - d(11, 9) = -21$.
$d(8, 9) - d(8, 1) = 59 - 1 = 58$.
Thus, the value of (11 8) is non-negative.
TRIAL 2. MIN(M )(11, 2) = 2.
$\qquad\qquad\qquad d(11, 2) - d(11, 9) = 14 - 32 = -18$.
$\qquad\qquad\qquad (11, 2) \rightarrow (11, 20)$
$\qquad\qquad\qquad (20, 2) \rightarrow (20, 20)$

(20, 2) is an arc of $D_7$.

$$(20, 15) \rightarrow (20, 18)$$

$d(20, 15) - d(20, 2) = 6 - 6 = 0$.

$$(18, 15) \rightarrow (18, 18)$$

(18, 15) is an arc of $D_7$.

$$(18, 14) \rightarrow (18, 6)$$

$d(18, 14) - d(18, 15) = 9 - 3 = 6$.

$$(6, 14) \rightarrow (6, 6)$$

(6, 14) is an arc of $D_7$.

$$(6, 2) \rightarrow (6, 20)$$

$d(6, 2) - d(6, 14) = 19 - 18 = 1$.

$P = [11, 20, 18, 6, 20]$

$d(6, 9) - d(6, 14) = 72 - 18 = 54$.

$s_{21} = (11^{-18} 20^0 18^6 6^{54})$.

The value of $s_{21}$ is non-negative.

$d(20, 9) - d(20, 2) = 35 - 6 = 29$.

$d(18, 9) - d(18, 15) = 83 - 3 = 80$.

$s_{211} = (11^{-21} 20^{29})$, $s_{212} = (11^{-21} 20^0 18^{80})$.

The value of both of the above permutations are non-negative.

$s_{22} = (20^0 18^6 6^1)$.

The value of $s_{22}$ is non-negative.

TRIAL 3. MIN(M)(11, 3) = 4.

$d(11, 4) - d(11, 9) = 18 - 32 = -14$.

$$(11, 4) \rightarrow (11, 9)$$
$$(9, 15) \rightarrow (9, 18)$$

$d(9, 15) - d(9, 4) = 4 - 6 = -2$.

$$(18, 15) \rightarrow (18, 18)$$

(18, 15) is an arc of $D_7$.

$$(18, 14) \rightarrow (18, 6)$$

$d(18, 14) - d(18, 15) = 9 - 3 = 6$.

$$(6, 14) \rightarrow (6, 6)$$

(6, 14) is an arc of $D_7$.

$$(6, 2) \rightarrow (6, 20)$$

$d(6, 2) - d(6, 14) = 19 - 18 = 1$.

$$(20, 2) \rightarrow (20, 20)$$

(20, 2) is an arc of $D_7$.

$$(20, 15) \rightarrow (20, 18)$$

$d(20, 15) - d(20, 2) = 6 - 6 = 0$.

$P = [11, 9, 18, 6, 20, 18]$.

$d(20, 9) - d(20, 2) = 35 - 6 = 29$.

$s_{31} = (11^{-14} 9^{-2} 18^6 6^1 20^{29})$.

The value of $s_{31}$ is non-negative.

$d(9, 9)$ is not permissible.

$d(18, 9) - d(18, 15) = 83 - 3 = 80$.

d(6, 9) − d(6, 14) = 72 − 18 = 54.
$s_{312}$ = ($11^{-14}9^{-2}18^{80}$), $s_{313}$ = ($11^{-14}9^{-2}18^{6}6^{54}$).
Each of the above permutations has a non-negative value.
$s_{32}$ = ($18^{6}6^{1}20^{0}$).
The value of $s_{32}$ is non-negative.
We next start with 15.
TRIAL 1.  MIN(M)(15, 1) = 16.
                (15, 16) → (15, 4)   -21
                (4, 17) → (4, 9)   -9
                (9, 15) → (9, 18)   -2
                (18, 15) → (18, 18).
(18, 15) is an arc of $D_7$.
                d(18, 14) - d(18, 15) = 6
                (18, 14) → (18, 9)
P = [15, 4, 9, 18, 9].
(18, 15) is derived from the loop (18, 18).
d(9, 18) - d(9, 4) = 68 - 6 = 62.
d(4, 18) - d(4, 16) = 67 - 18 = 49.
Thus, $s_{11}$ = ($15^{-21}4^{-9}9^{62}$) and ($15^{-21}4^{49}$) both have non-negative values.
$s_{21}$ = ($15^{-21}4^{49}$)($9^{-2}18^{6}$). Both cycles have non-negative values.
TRIAL 2. MIN(M)(15, 2) = 18. (15, 18) is an arc of $D_7$. Thus, we can proceed no farther with 15.
We start with 14.
TRIAL 1. MIN(M)(14, 1) = 17.
                (14, 17) → (14, 5)   -13
                (5, 14) → (5, 6)   -6
                (6, 14) → (6, 6).
(6, 14) is an arc of $D_7$.
                d(6, 2) - d(6, 14) = 19 - 18 = 1
                (6, 2) → (6, 20)   1
                (20, 2) → (20, 20).
(20, 2) is an arc of $D_7$.
                d(20, 15) - d(20, 2) = 0
                (20, 15) → (20, 18)   0
                (18, 15) → (18, 18).
(18, 15) is an arc of $D_7$.
                d(18, 14) - d(18, 15) = 9 - 3 = 6.
                (18, 14) → (18, 6)   6
P = [14, 5, 6, 20, 18, 6].
(18, 14) is derived from (18, 13) in M.
d(18, 13) - d(18, 15) = 42 - 3 = 39.
d(20, 13) - d(20, 2) = 65 - 6 = 59.
d(6, 13) - d(6, 14) = 48 - 18 = 30.
d(5, 13) - d(5, 17) = 91 - 10 = 81.

$s_{11} = (14^{-13}5^{-6}6^120^018^{39})$, $s_{111}= (14^{-13}5^{-6}6^120^{59})$, $s_{112} = (14^{-13}5^{-6}6^{20})$,

$s_{113} = (14^{-13}5^{81})$. The values of all of these permutations are non-negative.

$s_{21} = (14^{-13}5^{81})(6^120^018^6)$. Both cycles have non-negative values.

**TRIAL 2.** MIN(M)(14, 2) = 13. (14, 13) is an arc of $D_7$.

We have tested [log n] + 1 = 3 vertices and have obtained no negatively-valued permutations. Thus, we proceed to Phase 2,

$$D_7^{-1}M$$

|    | 7 | 8 | 11 | 17 | 18 | 14 | 5 | 1 | 4 | 12 | 9 | 20 | 19 | 13 | 16 | 6 | 10 | 15 | 3 | 2 |    |
|----|---|---|----|----|----|----|---|---|---|----|---|----|----|----|----|---|----|----|---|---|----|
|    | 1 | 2 | 3  | 4  | 5  | 6  | 7 | 8 | 9 | 10 | 11| 12 | 13 | 14 | 15 | 16| 17 | 18 | 19| 20|    |
| 1  | 9 | 59| 52 | 2  | 72 | 17 | 14| ∞ | 97| 50 | 39| 65 | 86 | 29 | 65 | 38| 46 | 48 | 72| 88| 1  |
| 2  | 99| 20| 24 | 54 | 62 | 45 | 67| 80| 72| 76 | 58| 45 | 25 | 32 | 62 | 18| 66 | 11 | 28| ∞ | 2  |
| 3  | 73| 76| 13 | 56 | 32 | 94 | 34| 55| 82| 33 | 97| 81 | 74 | 23 | 87 | 26| 17 | 76 | ∞ | 40| 3  |
| 4  | 94| 88| 38 | 9  | 67 | 62 | 92| 29| ∞ | 59 | 87| 31 | 93 | 94 | 18 | 63| 18 | 33 | 77| 76| 4  |
| 5  | 21| 64| 41 | 10 | 19 | 4  | ∞ | 17| 80| 18 | 77| 35 | 99 | 91 | 29 | 49| 54 | 35 | 60| 32| 5  |
| 6  | 98| 94| 30 | 69 | 70 | 18 | 19| 21| 75| 43 | 72| 71 | 22 | 48 | 82 | ∞ | 40 | 94 | 31| 19| 6  |
| 7  | ∞ | 74| 81 | 91 | 95 | 33 | 4 | 1 | 68| 59 | 28| 89 | 58 | 24 | 10 | 35| 44 | 9  | 82| 99| 7  |
| 8  | 52| ∞ | 25 | 24 | 78 | 65 | 52| 1 | 49| 48 | 59| 70 | 67 | 56 | 59 | 30| 92 | 34 | 71| 49| 8  |
| 9  | 66| 19| 90 | 77 | 68 | 21 | 63| 92| 6 | 24 | ∞ | 66 | 84 | 92 | 43 | 79| 63 | 4  | 49| 99| 9  |
| 10 | 92| 92| 34 | 73 | 99 | 63 | 80| 46| 4 | 10 | 72| 5  | 85 | 86 | 40 | 54| ∞  | 63 | 2 | 65| 10 |
| 11 | 46| 51| ∞  | 54 | 85 | 85 | 91| 11| 18| 90 | 32| 19 | 31 | 30 | 66 | 37| 23 | 18 | 81| 14| 11 |
| 12 | 79| 89| 77 | 68 | 65 | 42 | 6 | 43| 64| ∞  | 39| 7  | 25 | 78 | 62 | 39| 44 | 47 | 52| 1 | 12 |
| 13 | 15| 81| 45 | 16 | 67 | 54 | 20| 19| 76| 67 | 15| 4  | 7  | ∞  | 92 | 80| 87 | 9  | 93| 82| 13 |
| 14 | 31| 52| 23 | 6  | 95 | ∞  | 56| 84| 90| 54 | 28| 26 | 59 | 19 | 37 | 96| 26 | 91 | 89| 73| 14 |
| 15 | 62| 81| 98 | 79 | 39 | 44 | 65| 84| 42| 52 | 90| 91 | 50 | 59 | 18 | 54| 80 | ∞  | 65| 53| 15 |
| 16 | 88| 17| 62 | 69 | 61 | 20 | 81| 13| 60| 12 | 58| 57 | 20 | 59 | ∞  | 5 | 48 | 48 | 77| 5 | 16 |
| 17 | 83| 24| 19 | ∞  | 71 | 28 | 8 | 46| 59| 75 | 28| 6  | 31 | 17 | 75 | 83| 1  | 82 | 25| 48| 17 |
| 18 | 67| 51| 82 | 68 | ∞  | 9  | 26| 31| 98| 90 | 83| 41 | 51 | 42 | 26 | 80| 80 | 3  | 84| 50| 18 |
| 19 | 26| 17| 63 | 32 | 59 | 61 | 58| 53| 44| 63 | 64| 68 | ∞  | 87 | 57 | 48| 88 | 42 | 6 | 89| 19 |
| 20 | 14| 56| 24 | 91 | 85 | 47 | 49| 43| 51| 80 | 35| ∞  | 84 | 65 | 55 | 52| 18 | 6  | 73| 6 | 20 |
|    | 1 | 2 | 3  | 4  | 5  | 6  | 7 | 8 | 9 | 10 | 11| 12 | 13 | 14 | 15 | 16| 17 | 18 | 19| 20|    |

$$D_7^{-1}M^-$$

|   | 7 | 8 | 11 | 17 | 18 | 14 | 5 | 1 | 4 | 12 | 9 | 20 | 19 | 13 | 16 | 6 | 10 | 15 | 3 | 2 |   |
|---|---|---|----|----|----|----|---|---|---|----|---|----|----|----|----|---|----|----|---|---|---|
|   | 1 | 2 | 3  | 4  | 5  | 6  | 7 | 8 | 9 | 10 | 11| 12 | 13 | 14 | 15 | 16| 17 | 18 | 19| 20|   |
| 1 | 0 | 50| 43 | -7 | 63 | 8  | 5 | ∞ | 88| 42 | 30| 56 | 77 | 20 | 56 | 29| 37 | 39 | 63| 79| 1 |
| 2 | 79| 0 | 4  | 34 | 42 | 25 | 47| 60| 52| 56 | 38| 25 | 5  | 12 | 42 | -2| 46 | -9 | 8 | ∞ | 2 |
| 3 | 60| 63| 0  | 43 | 19 | 81 | 21| 42| 69| 20 | 64| 68 | 61 | 10 | 74 | 13| 4  | 63 | ∞ | 27| 3 |
| 4 | 85| 79| 29 | 0  | 58 | 53 | 83| 20| ∞ | 50 | 78| 22 | 84 | 85 | 9  | 54| 9  | 24 | 68| 67| 4 |
| 5 | 2 | 45| 22 | -9 | 0  | -15| ∞ | -2| 61| -1 | 58| 16 | 80 | 72 | 10 | 30| 35 | 16 | 41| 13| 5 |
| 6 | 80| 76| 12 | 51 | 52 | 0  | 1 | 3 | 57| 25 | 54| 53 | 4  | 30 | 64 | ∞ | 22 | 76 | 13| 1 | 6 |
| 7 | ∞ | 70| 77 | 87 | 91 | 29 | 0 | -3| 64| 55 | 24| 85 | 54 | 20 | 6  | 31| 40 | 5  | 78| 95| 7 |
| 8 | 53| ∞ | 24 | 23 | 77 | 64 | 51| 0 | 48| 47 | 58| 69 | 66 | 55 | 58 | 29| 91 | 33 | 70| 48| 8 |
| 9 | 60| 13| 84 | 71 | 62 | 15 | 57| 86| 0 | 18 | ∞ | 60 | 78 | 86 | 37 | 73| 57 | -2 | 43| 93| 9 |
| 10| 82| 82| 24 | 63 | 89 | 53 | 70| 36| -6| 0  | 62| -5 | 75 | 76 | 30 | 44| ∞  | 53 | -8| 55| 10|
| 11| 14| 19| ∞  | 22 | 53 | 53 | 59| -21| -14| 58| 0 | -13| -1 | -2 | 34 | 5 | -9 | -14| 49| -18| 11|
| 12| 72| 82| 70 | 61 | 58 | 35 | -1| 37| 57| ∞  | 32| 0  | 18 | 71 | 55 | 32| 37 | 40 | 43| -6| 12|
| 13| 8 | 74| 38 | 9  | 60 | 47 | 13| 12| 69| 60 | 8 | -3 | 0  | ∞  | 85 | 73| 80 | 2  | 86| 75| 13|
| 14| 12| 33| 4  | -13| 76 | ∞  | 37| 65| 71| 35 | 9 | 7  | 40 | 0  | 18 | 77| 7  | 72 | 70| 54| 14|
| 15| 44| 63| 80 | 61 | 21 | 26 | 47| 66| 24| 34 | 72| 73 | 32 | 41 | 0  | 36| 62 | ∞  | 47| 35| 15|
| 16| 83| 12| 57 | 64 | 56 | 15 | 76| 8 | 55| 7  | 53| 52 | 15 | 54 | ∞  | 0 | 43 | 43 | 72| 0 | 16|
| 17| 82| 23| 18 | ∞  | 70 | 27 | 7 | 45| 58| 74 | 27| 5  | 30 | 16 | 74 | 82| 0  | 81 | 24| 47| 17|
| 18| 64| 63| 79 | 65 | ∞  | 6  | 23| 28| 95| 87 | 80| 38 | 48 | 39 | 23 | 77| 77 | 0  | 81| 47| 18|
| 19| 20| 48| 57 | 26 | 59 | 55 | 52| 47| 38| 57 | 58| 62 | ∞  | 81 | 51 | 42| 82 | 36 | 0 | 83| 19|
| 20| 8 | 11| 18 | 85 | 79 | 41 | 43| 37| 45| 74 | 29| ∞  | 78 | 59 | 49 | 46| 12 | 0  | 67| 0 | 20|
|   | 1 | 2 | 3  | 4  | 5  | 6  | 7 | 8 | 9 | 10 | 11| 12 | 13 | 14 | 15 | 16| 17 | 18 | 19| 20|   |

In what follows call an entry of $D_8^{-1}M^-(20i)$, i = 2, 3, 4, ... an *initial negative entry* if it belongs to $D_8^{-1}M^-$, an *inactive negative entry* if its value represents that of a negative path which cannot be extended to a longer negative path. An *active negative entry* is one which represents a path which can be extended in the current $D_8^{-1}M^-(20i)$ to a path in $D_8^{-1}M^-(20(i+1))$. An inactive negative entry is expressed in italics. An active negative entry which has extended a path from (a, b) in $D_8^{-1}M^-(20i)$ to (a, c) in $D_8^{-1}M^-(20(i+1))$ has the entry (a, c) underlined in $D_8^{-1}M^-(20(i+1))$. If an active entry, say (a, b), has been extended to a new path, (a, c), the original value at (a, b) is written in italics unless its value has been superceded by a smaller one during the iteration. The reason for doing this is that each active entry (a, b) represents the smallest-valued simple negative path going from a to b *up until that point in the algorithm*. Thus, two different types of entries are written in italics. An *iteration* of the algorithm occurs when we go from $D_8^{-1}M^-(20i)$ to $D_8^{-1}M^-(20(i+1))$.

Ex. Let (7, 12) be an entry with a value -15. Assume that (12, 19) has the value 5 and that the value for (7, 19) in $D_8^{-1}M^-$ is -8. Since (7, 12)(12, 19) yields a path starting at 7 and ending at 19, (7, 12) is an active entry that may be extended to the path [7 ... 12 19] that has a value -10. On the other hand, we now consider the value -15 of (7, 12). This entry itself now becomes an inactive entry with a value of -15. It can become no smaller unless we can obtain some path in the future starting at 7 and ending at 12 that has a value less than -15.

An initial negative entry is presented in ordinary type without underlining. We write entries in italics to immediately tell us if the possible addition of an arc to form a path [a, b] has a value which is no greater than a previously obtained value for the path whose value is given at (a, b).

We now present negative paths in $D_7^{-1}M^-$.

j = 4.
(14, 4)(4, 15): -4;  (14, 4)(4, 17): -4,
j = 6.
(5, 6)(6, 3): -3; (5, 6)(6, 7): -14;  (5, 6)(6, 8): -12;  (5, 6)(6, 13): -11; (5, 6)(6, 19): -2; (5, 6)(6, 20): -14 .
j = 7.
(5, 7)(7, 8): -17;  (5, 7)(7, 15): -8;  (5, 7)(7, 18): -9;  (12, 7)(7, 8): -4.
j = 9.
(10, 9)(9, 18): -8;  (11, 9)(9, 2): -1;  (11, 9)(9, 18): -16.
j = 10.
(5, 10)(10, 9): -7;  (5, 10)(10, 12): -6;  (5, 10)(10, 18): -9.
j = 12.
(5, 12)(12, 7): -7;  (5, 12)(12, 8): -10;  (5, 12)(12, 20): -12;
(10, 12)(12, 8): -9;  (10, 12)(12, 20): -11;  (11, 12)(12, 7): -14;
(11, 12)(12, 8): -17;  (11, 12)(12, 20): -19; (13, 12)(12, 7): -4;
(13, 12)(12, 8): -7;  (13, 12)(12, 20): -9.
j = 13.
(5, 13)(13, 1): -3;  (5, 13)(13, 4): -2;  (5, 13)(13, 7): -15;
(5, 13)(13, 8): -7;  (5, 13)(13, 11): -3;  (5, 13)(13, 18): -9;
(5, 13)(13, 20): -20; (11, 13)(13, 7): -5;  (11, 13)(13, 12): -4;
(11, 13)(13, 20): -10.

**j = 14.**
**(11, 14)(14, 4): -15;  (11, 14)(14, 15): -6;  (11, 14)(14, 17): -6.**
**j = 16.**
**(2, 16)(16, 20): -2.**
**j = 17.**
**(11, 17)(17, 7): -2;  (11, 17)(17, 12): -4.**
**j = 18.**
**(5, 18)(18, 6): -3;  (11, 18)(18, 6): -10.**
**j = 20.**
**(5, 20)(20, 1): -12;  (5, 20)(20, 2): -9;  (5, 20)(20, 3): -2;**
**(5, 20) (20, 17): -8;  (5, 20)(20, 18): -20;  (10, 20)(20, 1): -3;**
**(10, 20)(20, 18): -11;  (10, 20)(20, 1): -3;  (10, 20)(20, 18): -11;**
**(11, 20)(20, 1): -11;  (11, 20)(20, 2): -8;  (11, 20)(20, 3): -1;**
**(11, 20)(20, 17): -7;  (11, 20)(20, 18): -13;  (13, 20)(20, 1): -1;**
**(13, 20)(20, 18): -9.**
**None of the negative paths given has the property that if an arbitrary path is (a, b)(b, c), then (a, b)(b, c)(c, a) is a negative cycle. In order to illustrate this better, we give the path matrix, $P_{20}$, which contains the smallest negatively-valued paths obtainable using the F-W algorithm. We note that as we have shown in Theorem 3 and its proof that a simple path can only become non-simple (form a cycle from a subpath) by adding an arc to it if the value of the cycle is negative. But if j =r when the cycle is formed, then the negative cycle formed would have been obtained when j ≤ r.**

**P$_{20}$**

|    | 1  | 2  | 3  | 4  | 5  | 6  | 7  | 8  | 9  | 10 | 11 | 12 | 13 | 14 | 15 | 16 | 17 | 18 | 19 | 20 |    |
|----|----|----|----|----|----|----|----|----|----|----|----|----|----|----|----|----|----|----|----|----|----|
| 1  |    |    |    |    |    |    |    |    |    |    |    |    |    |    |    |    |    |    |    |    | 1  |
| 2  |    |    |    |    | 18 |    |    |    |    |    |    |    |    |    |    |    |    |    | 16 |    | 2  |
| 3  |    |    |    |    |    |    |    |    |    |    |    |    |    |    |    |    |    |    |    |    | 3  |
| 4  |    |    |    |    |    |    |    |    |    |    |    |    |    |    |    |    |    |    |    |    | 4  |
| 5  | 20 | 20 | 6  | 13 |    | 13 |    | 10 |    | 13 | 13 | 6  |    | 7  |    |    | 20 | 10 | 12 |    | 5  |
| 6  |    |    |    |    |    |    |    |    |    |    |    |    |    |    |    |    |    |    |    |    | 6  |
| 7  |    |    |    |    |    |    |    |    |    |    |    |    |    |    |    |    |    |    |    |    | 7  |
| 8  |    |    |    |    |    |    |    |    |    |    |    |    |    |    |    |    |    |    |    |    | 8  |
| 9  |    |    |    |    |    |    |    |    |    |    |    |    |    |    |    |    |    |    |    |    | 9  |
| 10 | 20 |    |    |    |    |    | 12 | 12 |    |    |    |    |    |    |    |    | 20 |    | 12 |    | 10 |
| 11 | 20 | 9  | 20 | 14 |    | 18 | 12 |    |    |    |    |    |    | 14 |    |    | 20 |    | 12 |    | 11 |
| 12 |    |    |    |    |    |    | 7  |    |    |    |    |    |    |    |    |    |    |    |    |    | 12 |
| 13 | 20 |    |    |    |    |    |    | 12 |    |    |    |    |    |    |    |    | 20 |    | 12 |    | 13 |
| 14 |    |    |    |    |    |    |    |    |    |    |    |    |    | 4  |    | 4  |    |    |    |    | 14 |
| 15 |    |    |    |    |    |    |    |    |    |    |    |    |    |    |    |    |    |    |    |    | 15 |
| 16 |    |    |    |    |    |    |    |    |    |    |    |    |    |    |    |    |    |    |    |    | 16 |
| 17 |    |    |    |    |    |    |    |    |    |    |    |    |    |    |    |    |    |    |    |    | 17 |
| 18 |    |    |    |    |    |    |    |    |    |    |    |    |    |    |    |    |    |    |    |    | 18 |
| 19 |    |    |    |    |    |    |    |    |    |    |    |    |    |    |    |    |    |    |    |    | 19 |
| 20 |    |    |    |    |    |    |    |    |    |    |    |    |    |    |    |    |    |    |    |    | 20 |
|    | 1  | 2  | 3  | 4  | 5  | 6  | 7  | 8  | 9  | 10 | 11 | 12 | 13 | 14 | 15 | 16 | 17 | 18 | 19 | 20 |    |

We next give $D_7^{-1}M^-(20)$. It will contain all of the negative paths given in $P_{20}$. In order to guarantee that there are no negative cycles in $D_7^{-1}M^-$, we need to use all of the paths given in $P_{20}$ to see if we can add an arc to each path considered. Note that each path is actually represented by an entry which gives the value of the path together with the initial and terminal vertices of the path. If we go beyond $D_7^{-1}M^-(20)$, we can only add one arc of $D_7^{-1}M^-$ at a time to extend a path. We *do not* add more than one arc of $D_7^{-1}M^-$ to a current path of $D_7^{-1}M^-(20i)$ (i = 2,3, ...) h which we are trying to extend.

$$D_7^{-1}M^-(20)$$

| | 7 | 8 | 11 | 17 | 18 | 14 | 5 | 1 | 4 | 12 | 9 | 20 | 19 | 13 | 16 | 6 | 10 | 15 | 3 | 2 | |
|---|---|---|---|---|---|---|---|---|---|---|---|---|---|---|---|---|---|---|---|---|---|
| | 1 | 2 | 3 | 4 | 5 | 6 | 7 | 8 | 9 | 10 | 11 | 12 | 13 | 14 | 15 | 16 | 17 | 18 | 19 | 20 | |
| 1 | 0 | 50 | 43 | -7 | 63 | 8 | 5 | ∞ | 88 | 42 | 30 | 56 | 77 | 20 | 56 | 29 | 37 | 39 | 63 | 79 | 1 |
| 2 | 79 | 0 | 4 | 34 | 42 | -3 | 47 | 60 | 52 | 56 | 38 | 25 | 5 | 12 | 42 | -2 | 46 | -9 | 8 | -2 | 2 |
| 3 | 60 | 63 | 0 | 43 | 19 | 81 | 21 | 42 | 69 | 20 | 64 | 68 | 61 | 10 | 74 | 13 | 4 | 63 | ∞ | 27 | 3 |
| 4 | 85 | 79 | 29 | 0 | 58 | 53 | 83 | 20 | ∞ | 50 | 78 | 22 | 84 | 85 | 9 | 54 | 9 | 24 | 68 | 67 | 4 |
| 5 | -12 | -9 | -3 | -9 | 0 | -15 | -15 | -18 | -7 | -1 | -3 | -6 | -11 | 72 | -8 | 30 | -8 | -20 | -2 | -20 | 5 |
| 6 | 80 | 76 | 12 | 51 | 52 | 0 | 1 | 3 | 57 | 25 | 54 | 53 | 4 | 30 | 64 | ∞ | 22 | 76 | 13 | 1 | 6 |
| 7 | ∞ | 70 | 77 | 87 | 91 | 29 | 0 | -3 | 64 | 55 | 24 | 85 | 54 | 20 | 6 | 31 | 40 | 5 | 78 | 95 | 7 |
| 8 | 53 | ∞ | 24 | 23 | 77 | 64 | 51 | 0 | 48 | 47 | 58 | 69 | 66 | 55 | 58 | 29 | 91 | 33 | 70 | 48 | 8 |
| 9 | 60 | 13 | 84 | 71 | 62 | 15 | 57 | 86 | 0 | 18 | ∞ | 60 | 78 | 86 | 37 | 73 | 57 | -2 | 43 | 93 | 9 |
| 10 | -3 | 82 | 24 | 63 | 89 | 53 | -6 | -9 | -6 | 0 | 62 | -5 | 75 | 76 | 30 | 44 | ∞ | -11 | -8 | -11 | 10 |
| 11 | -11 | -8 | -1 | -15 | 53 | -10 | -14 | -21 | -14 | 58 | 0 | -13 | -1 | -2 | -6 | 5 | -9 | -19 | 49 | -19 | 11 |
| 12 | 72 | 82 | 70 | 61 | 58 | 35 | -1 | -4 | 57 | ∞ | 32 | 0 | 18 | 71 | 55 | 32 | 37 | 40 | 43 | -6 | 12 |
| 13 | -1 | 74 | 38 | 9 | 60 | 47 | -15 | -7 | 69 | 60 | 8 | -3 | 0 | ∞ | 85 | 73 | 80 | -9 | 86 | -9 | 13 |
| 14 | 12 | 33 | 4 | -13 | 76 | ∞ | 37 | 65 | 71 | 35 | 9 | 7 | 40 | 0 | 18 | 77 | 7 | 72 | 70 | 54 | 14 |
| 15 | 44 | 63 | 80 | 61 | 21 | 26 | 47 | 66 | 24 | 34 | 72 | 73 | 32 | 41 | 0 | 36 | 62 | ∞ | 47 | 35 | 15 |
| 16 | 83 | 12 | 57 | 64 | 56 | 15 | 76 | 8 | 55 | 7 | 53 | 52 | 15 | 54 | ∞ | 0 | 43 | 43 | 72 | 0 | 16 |
| 17 | 82 | 23 | 18 | ∞ | 70 | 27 | 7 | 45 | 58 | 74 | 27 | 5 | 30 | 16 | 74 | 82 | 0 | 81 | 24 | 47 | 17 |
| 18 | 64 | 63 | 79 | 65 | ∞ | 6 | 23 | 28 | 95 | 87 | 80 | 38 | 48 | 39 | 23 | 77 | 77 | 0 | 81 | 47 | 18 |
| 19 | 20 | 48 | 57 | 26 | 59 | 55 | 52 | 47 | 38 | 57 | 58 | 62 | ∞ | 81 | 51 | 42 | 82 | 36 | 0 | 83 | 19 |
| 20 | 8 | 11 | 18 | 85 | 79 | 41 | 43 | 37 | 45 | 74 | 29 | ∞ | 78 | 59 | 49 | 46 | 12 | 0 | 67 | 0 | 20 |
| | 1 | 2 | 3 | 4 | 5 | 6 | 7 | 8 | 9 | 10 | 11 | 12 | 13 | 14 | 15 | 16 | 17 | 18 | 19 | 20 | |

We now extend the paths of $D_7^{-1}M^-(20)$ by using the F-W negatively-valued path algorithm to add arcs of $D_7^{-1}M^-$, one at a time, to each path. Using MIN(M), we will have gone through all columns a second time, i.e., the new value matrix, $D_7^{-1}M^-(40)$ is an iteration of $D_7^{-1}M^-(20)$. We include those underlined values, d(a, b), of $D_7^{-1}M^-(20)$ that have not been superceded by smaller values during the iteration as italicized elements of $D_7^{-1}M^-(40)$.

j = 1.
(5, 1)(1, 4): -19;  (10, 1)(1, 4): -10;  (11, 1)(1, 4): -18;  (13, 1)(1, 4): -8.
j = 2.
(5, 2)(2, 3): -5;  (5, 2)(2, 16): -11;  (11, 2)(2, 3): -4;  (11, 2)(2, 16): -10.
j = 6.
(2, 6)(6, 7): -2;  (11, 6)(6, 13): -9.
j = 7.
(5, 7)(7, 15): -9;  (11, 7)(7, 15): -8;  (13, 7)(7, 8): -18;  (13, 7)(7, 15): -9.
j = 11.
(5, 11)(11, 8): -24;  (5, 11)(11, 9): -17;  (5, 11)(11, 12): -16;   (5, 11)(11, 14): -5;
(5, 11)(11, 17): -12;  (5, 11)(11, 20): -21.
j = 13.
(5, 13)(13, 12): -14.
j = 18.
(10, 18)(18, 6): -5;  (11, 18)(18, 6): -9;  (13, 18)(18, 6): -3.

Using the paths given in $P_{20}$, we now add arcs to obtain $P_{40}$. After constructing each path, we check to see if adding the arc connecting the terminal vertex of a path to its initial vertex will yield a negative cycle.

| | 1 | 2 | 3 | 4 | 5 | 6 | 7 | 8 | 9 | 10 | 11 | 12 | 13 | 14 | 15 | 16 | 17 | 18 | 19 | 20 | |
|---|---|---|---|---|---|---|---|---|---|---|---|---|---|---|---|---|---|---|---|---|---|
| 1 | | | | | | | | | | | | | | | | | | | | | 1 |
| 2 | | | | | | 18 | 6 | | | | | | | | | | | | | 16 | 2 |
| 3 | | | | | | | | | | | | | | | | | | | | | 3 |
| 4 | | | | | | | | | | | | | | | | | | | | | 4 |
| 5 | 20 | 20 | 2 | 1 | | | 13 | 11 | 10 | | 13 | 13 | 6 | 11 | 7 | 2 | 11 | 20 | 10 | 11 | 5 |
| 6 | | | | | | | | | | | | | | | | | | | | | 6 |
| 7 | | | | | | | | | | | | | | | | | | | | | 7 |
| 8 | | | | | | | | | | | | | | | | | | | | | 8 |
| 9 | | | | | | | | | | | | | | | | | | | | | 9 |
| 10 | | | | 1 | | 18 | 12 | 12 | | | | | | | | | 20 | | 12 | | 10 |
| 11 | 20 | 9 | 2 | 1 | | 18 | 12 | | | | | | 6 | | 7 | 2 | | 20 | | 12 | 11 |
| 12 | | | | | | | 7 | | | | | | | | | | | | | | 12 |
| 13 | 20 | | | 1 | | 18 | 12 | 7 | | | | | | | | 7 | | 20 | | 12 | 13 |
| 14 | | | | | | | | | | | | | | | 4 | | 4 | | | | 14 |
| 15 | | | | | | | | | | | | | | | | | | | | | 15 |
| 16 | | | | | | | | | | | | | | | | | | | | | 16 |
| 17 | | | | | | | | | | | | | | | | | | | | | 17 |
| 18 | | | | | | | | | | | | | | | | | | | | | 18 |
| 19 | | | | | | | | | | | | | | | | | | | | | 19 |
| 20 | | | | | | | | | | | | | | | | | | | | | 20 |
| | 1 | 2 | 3 | 4 | 5 | 6 | 7 | 8 | 9 | 10 | 11 | 12 | 13 | 14 | 15 | 16 | 17 | 18 | 19 | 20 | |

$$D_7^{-1}M^-(40)$$

|   | 7 | 8 | 11 | 17 | 18 | 14 | 5 | 1 | 4 | 12 | 9 | 20 | 19 | 13 | 16 | 6 | 10 | 15 | 3 | 2 |   |
|---|---|---|---|---|---|---|---|---|---|---|---|---|---|---|---|---|---|---|---|---|---|
|   | 1 | 2 | 3 | 4 | 5 | 6 | 7 | 8 | 9 | 10 | 11 | 12 | 13 | 14 | 15 | 16 | 17 | 18 | 19 | 20 |   |
| 1 | 0 | 50 | 43 | *-7* | 63 | 8 | 5 | ∞ | 88 | 42 | 30 | 56 | 77 | 20 | 56 | 29 | 37 | 39 | 63 | 79 | 1 |
| 2 | 79 | 0 | 4 | 34 | 42 | *-3* | *-2* | 60 | 52 | 56 | 38 | 25 | 5 | 12 | 42 | *-2* | 46 | *-9* | 8 | *-2* | 2 |
| 3 | 60 | 63 | 0 | 43 | 19 | 81 | 21 | 42 | 69 | 20 | 64 | 68 | 61 | 10 | 74 | 13 | 4 | 63 | ∞ | 27 | 3 |
| 4 | 85 | 79 | 29 | 0 | 58 | 53 | 83 | 20 | ∞ | 50 | 78 | 22 | 84 | 85 | 9 | 54 | 9 | 24 | 68 | 67 | 4 |
| 5 | *-12* | *-9* | *-5* | *-19* | 0 | -15 | -15 | *-24* | *-17* | -1 | *-3* | *-16* | *-16* | *-5* | *-9* | *-11* | *-12* | *-20* | -2 | *-21* | 5 |
| 6 | 80 | 76 | 12 | 51 | 52 | 0 | 1 | 3 | 57 | 25 | 54 | 53 | 4 | 30 | 64 | ∞ | 22 | 76 | 13 | 1 | 6 |
| 7 | ∞ | 70 | 77 | 87 | 91 | 29 | 0 | *-3* | 64 | 55 | 24 | 85 | 54 | 20 | 6 | 31 | 40 | 5 | 78 | 95 | 7 |
| 8 | 53 | ∞ | 24 | 23 | 77 | 64 | 51 | 0 | 48 | 47 | 58 | 69 | 66 | 55 | 58 | 29 | 91 | 33 | 70 | 48 | 8 |
| 9 | 60 | 13 | 84 | 71 | 62 | 15 | 57 | 86 | 0 | 18 | ∞ | 60 | 78 | 86 | 37 | 73 | 57 | *-2* | 43 | 93 | 9 |
| 10 | *-3* | 82 | 24 | *-10* | 89 | *-5* | *-6* | *-9* | *-6* | 0 | 62 | *-5* | 75 | 76 | 30 | 44 | ∞ | *-11* | *-8* | *-11* | 10 |
| 11 | *-11* | *-8* | *-4* | *-18* | 53 | *-13* | *-14* | -21 | -14 | 58 | 0 | *-13* | *-12* | -2 | *-8* | *-10* | -9 | *-19* | 49 | *-19* | 11 |
| 12 | 72 | 82 | 70 | 61 | 58 | 35 | *-1* | *-4* | 57 | ∞ | 32 | 0 | 18 | 71 | 55 | 32 | 37 | 40 | 43 | *-6* | 12 |
| 13 | *-1* | 74 | 38 | *-8* | 60 | *-3* | *-15* | *-18* | 69 | 60 | 8 | *-3* | 0 | ∞ | *-9* | 73 | 80 | *-10* | 86 | *-9* | 13 |
| 14 | 12 | 33 | 4 | *-13* | 76 | ∞ | 37 | 65 | 71 | 35 | 9 | 7 | 40 | 0 | 18 | 77 | 7 | 72 | 70 | 54 | 14 |
| 15 | 44 | 63 | 80 | 61 | 21 | 26 | 47 | 66 | 24 | 34 | 72 | 73 | 32 | 41 | 0 | 36 | 62 | ∞ | 47 | 35 | 15 |
| 16 | 83 | 12 | 57 | 64 | 56 | 15 | 76 | 8 | 55 | 7 | 53 | 52 | 15 | 54 | ∞ | 0 | 43 | 43 | 72 | 0 | 16 |
| 17 | 82 | 23 | 18 | ∞ | 70 | 27 | 7 | 45 | 58 | 74 | 27 | 5 | 30 | 16 | 74 | 82 | 0 | 81 | 24 | 47 | 17 |
| 18 | 64 | 63 | 79 | 65 | ∞ | 6 | 23 | 28 | 95 | 87 | 80 | 38 | 48 | 39 | 23 | 77 | 77 | 0 | 81 | 47 | 18 |
| 19 | 20 | 48 | 57 | 26 | 59 | 55 | 52 | 47 | 38 | 57 | 58 | 62 | ∞ | 81 | 51 | 42 | 82 | 36 | 0 | 83 | 19 |
| 20 | 8 | 11 | 18 | 85 | 79 | 41 | 43 | 37 | 45 | 74 | 29 | ∞ | 78 | 59 | 49 | 46 | 12 | 0 | 67 | 0 | 20 |
|   | 1 | 2 | 3 | 4 | 5 | 6 | 7 | 8 | 9 | 10 | 11 | 12 | 13 | 14 | 15 | 16 | 17 | 18 | 19 | 20 |   |

Using $D_7^{-1}M^-(40)$ and arcs from $D_7^{-1}M^-$, we construct $D_7^{-1}M^-(60)$.

After we obtain each new path, say P = [a, ..., z], we check to see if the 2-cycle (a  z) has a negative value. We then place the new arc obtained in $P_{40}$ in order to proceed to construct $P_{60}$. If (a z) has a negative value, we extract the cycle directly from our partially constructed $P_{60}$.

j = 4.
(5, 4)(4, 15): -10;  (10, 4)(4, 15): -1;  (10, 4)(4, 17): -1;  (11, 4)(4, 15): -9.
j = 6.
(10, 6)(6, 13): -1;  (11, 6)(6, 13): -9.
But the arc (13, 11) in $D_7^{-1}M^-$ has the value 8. Thus, (11 13) is a 2-cycle with value −1. We now go to $P_{40}$ where we have placed the new arcs obtained so far.

$P_{60}$

|    | 1  | 2  | 3  | 4 | 5 | 6  | 7  | 8  | 9  | 10 | 11 | 12 | 13 | 14 | 15 | 16 | 17 | 18 | 19 | 20 |    |
|----|----|----|----|---|---|----|----|----|----|----|----|----|----|----|----|----|----|----|----|----|----|
| 1  |    |    |    |   |   |    |    |    |    |    |    |    |    |    |    |    |    |    |    |    | 1  |
| 2  |    |    |    |   |   | 18 | 6  |    |    |    |    |    |    |    |    |    |    |    |    | 16 | 2  |
| 3  |    |    |    |   |   |    |    |    |    |    |    |    |    |    |    |    |    |    |    |    | 3  |
| 4  |    |    |    |   |   |    |    |    |    |    |    |    |    |    |    |    |    |    |    |    | 4  |
| 5  | 20 | 20 | 6  | 1 |   |    | 13 | 13 | 11 |    | 13 | 13 | 6  | 11 | 4  |    | 4  | 20 | 10 | 11 | 5  |
| 6  |    |    |    |   |   |    |    |    |    |    |    |    |    |    |    |    |    |    |    |    | 6  |
| 7  |    |    |    |   |   |    |    |    |    |    |    |    |    |    |    |    |    |    |    |    | 7  |
| 8  |    |    |    |   |   |    |    |    |    |    |    |    |    |    |    |    |    |    |    |    | 8  |
| 9  |    |    |    |   |   |    |    |    |    |    |    |    |    |    |    |    |    |    |    |    | 9  |
| 10 |    |    |    | 1 |   | 18 | 12 | 12 |    |    |    |    |    |    | 4  |    | 4  | 20 |    | 12 | 10 |
| 11 | 20 | 9  | 20 | 1 |   | 18 | 6  | 6  |    |    | 13 |    | 6  |    | 7  | 2  |    | 20 |    | 12 | 11 |
| 12 |    |    |    |   |   |    | 7  |    |    |    |    |    |    |    |    |    |    |    |    |    | 12 |
| 13 | 20 |    |    | 1 |   | 18 |    | 7  |    |    |    |    |    |    |    |    |    | 20 |    | 12 | 13 |
| 14 |    |    |    |   |   |    |    |    |    |    |    |    |    |    | 4  |    | 4  |    |    |    | 14 |
| 15 |    |    |    |   |   |    |    |    |    |    |    |    |    |    |    |    |    |    |    |    | 15 |
| 16 |    |    |    |   |   |    |    |    |    |    |    |    |    |    |    |    |    |    |    |    | 16 |
| 17 |    |    |    |   |   |    |    |    |    |    |    |    |    |    |    |    |    |    |    |    | 17 |
| 18 |    |    |    |   |   |    |    |    |    |    |    |    |    |    |    |    |    |    |    |    | 18 |
| 19 |    |    |    |   |   |    |    |    |    |    |    |    |    |    |    |    |    |    |    |    | 19 |
| 20 |    |    |    |   |   |    |    |    |    |    |    |    |    |    |    |    |    |    |    |    | 20 |
|    | 1  | 2  | 3  | 4 | 5 | 6  | 7  | 8  | 9  | 10 | 11 | 12 | 13 | 14 | 15 | 16 | 17 | 18 | 19 | 20 |    |

Given the negative path [11, 6] and the arc (6, 13) of $D_7^{-1}M^-$, we illustrate how recover the negative path [11, 6] from $P_{60}$. We first see that 18 goes into 6 in row 11. Therefore, our path expands to [11, 18, 6]. Next, checking column 18 of row 11, we have 20 in that entry. Our path is now [11, 20, 18, 6]. Next, 12 goes into 20 yielding [11, 12, 20, 18, 6]. Since there is no entry at (11, 12), we have our complete path. Thus, our negative cycle is $s_7$ = (11 12 20 18 6 13). Since $D_8 = D_7 s_7$, it follows that $D_8^{-1} = s_7^{-1} D_7^{-1}$. Thus, we apply the action of $s_7^{-1}$ to those columns of $D_7^{-1}M$ corresponding to the points moved by $s_7^{-1}$.

$s_7^{-1}$ = (11 13 6 18 20 12).
9 → 11 → 13, 19 → 13 → 6, 14 → 6 → 18, 15 → 18 → 20,
2 → 20 → 12, 20 → 12 → 11.
We now construct $D_8^{-1}M$.

$$D_8^{-1}M$$

| | 7 | 8 | 11 | 17 | 18 | 19* | 5 | 1 | 4 | 12 | 20* | 2* | 9* | 13 | 16 | 6 | 10 | 14* | 3 | 15* | |
|---|---|---|---|---|---|---|---|---|---|---|---|---|---|---|---|---|---|---|---|---|---|
| | 1 | 2 | 3 | 4 | 5 | 6 | 7 | 8 | 9 | 10 | 11 | 12 | 13 | 14 | 15 | 16 | 17 | 18 | 19 | 20 | |
| 1 | 9 | 59 | 52 | 2 | 72 | 86 | 14 | ∞ | 97 | 50 | 65 | 88 | 39 | 29 | 65 | 38 | 46 | 17 | 72 | 48 | 1 |
| 2 | 99 | 20 | 24 | 54 | 62 | 25 | 67 | 80 | 72 | 76 | 45 | ∞ | 58 | 32 | 62 | 18 | 66 | 45 | 28 | 11 | 2 |
| 3 | 73 | 76 | 13 | 56 | 32 | 74 | 34 | 55 | 82 | 33 | 81 | 40 | 97 | 23 | 87 | 26 | 17 | 94 | ∞ | 76 | 3 |
| 4 | 94 | 88 | 38 | 9 | 67 | 93 | 92 | 29 | ∞ | 59 | 31 | 76 | 87 | 94 | 18 | 63 | 18 | 62 | 77 | 33 | 4 |
| 5 | 21 | 64 | 41 | 10 | 19 | 99 | ∞ | 17 | 80 | 18 | 35 | 32 | 77 | 91 | 29 | 49 | 54 | 4 | 60 | 35 | 5 |
| 6 | 98 | 94 | 30 | 69 | 70 | 22 | 19 | 21 | 75 | 43 | 71 | 19 | 72 | 48 | 82 | ∞ | 40 | 18 | 31 | 94 | 6 |
| 7 | ∞ | 74 | 81 | 91 | 95 | 58 | 4 | 1 | 68 | 59 | 89 | 99 | 28 | 24 | 10 | 35 | 44 | 33 | 82 | 9 | 7 |
| 8 | 52 | ∞ | 25 | 24 | 78 | 67 | 52 | 1 | 49 | 48 | 70 | 49 | 59 | 56 | 59 | 30 | 92 | 65 | 71 | 34 | 8 |
| 9 | 66 | 19 | 90 | 77 | 68 | 84 | 63 | 92 | 6 | 24 | 66 | 99 | ∞ | 92 | 43 | 79 | 63 | 21 | 49 | 4 | 9 |
| 10 | 92 | 92 | 34 | 73 | 99 | 85 | 80 | 46 | 4 | 10 | 5 | 65 | 72 | 86 | 40 | 54 | ∞ | 63 | 2 | 63 | 10 |
| 11 | 46 | 51 | ∞ | 54 | 85 | 31 | 91 | 11 | 18 | 90 | 19 | 14 | 32 | 30 | 66 | 37 | 23 | 85 | 81 | 18 | 11 |
| 12 | 79 | 89 | 77 | 68 | 65 | 25 | 6 | 43 | 64 | ∞ | 7 | 1 | 39 | 78 | 62 | 39 | 44 | 42 | 52 | 47 | 12 |
| 13 | 15 | 81 | 45 | 16 | 67 | 7 | 20 | 19 | 76 | 67 | 4 | 82 | 15 | ∞ | 92 | 80 | 87 | 54 | 93 | 9 | 13 |
| 14 | 31 | 52 | 23 | 6 | 95 | 59 | 56 | 84 | 90 | 54 | 26 | 73 | 28 | 19 | 37 | 96 | 26 | ∞ | 89 | 91 | 14 |
| 15 | 62 | 81 | 98 | 79 | 39 | 50 | 65 | 84 | 42 | 52 | 91 | 53 | 90 | 59 | 18 | 54 | 80 | 44 | 65 | ∞ | 15 |
| 16 | 88 | 17 | 62 | 69 | 61 | 20 | 81 | 13 | 60 | 12 | 57 | 5 | 58 | 59 | ∞ | 5 | 48 | 20 | 77 | 48 | 16 |
| 17 | 83 | 24 | 19 | ∞ | 71 | 31 | 8 | 46 | 59 | 75 | 6 | 48 | 28 | 17 | 75 | 83 | 1 | 28 | 25 | 82 | 17 |
| 18 | 67 | 51 | 82 | 68 | ∞ | 51 | 26 | 31 | 98 | 90 | 41 | 50 | 83 | 42 | 26 | 80 | 80 | 9 | 84 | 3 | 18 |
| 19 | 26 | 17 | 63 | 32 | 59 | ∞ | 58 | 53 | 44 | 63 | 68 | 89 | 64 | 87 | 57 | 48 | 88 | 61 | 6 | 42 | 19 |
| 20 | 14 | 56 | 24 | 91 | 85 | 84 | 49 | 43 | 51 | 80 | ∞ | 6 | 35 | 65 | 55 | 52 | 18 | 47 | 73 | 6 | 20 |
| | 1 | 2 | 3 | 4 | 5 | 6 | 7 | 8 | 9 | 10 | 11 | 12 | 13 | 14 | 15 | 16 | 17 | 18 | 19 | 20 | |

$$D_8^{-1}M$$

|    | 7 | 8 | 11 | 17 | 18 | 19 | 5 | 1 | 4 | 12 | 20 | 2 | 9 | 13 | 16 | 6 | 10 | 14 | 3 | 15 |    |
|    | 1 | 2 | 3 | 4 | 5 | 6 | 7 | 8 | 9 | 10 | 11 | 12 | 13 | 14 | 15 | 16 | 17 | 18 | 19 | 20 |    |
|---|---|---|---|---|---|---|---|---|---|---|---|---|---|---|---|---|---|---|---|---|---|
| 1 | 0 | 50 | 43 | -7 | 63 | 77 | 5 | ∞ | 88 | 41 | 56 | 79 | 30 | 20 | 56 | 29 | 37 | 8 | 63 | 39 | 1 |
| 2 | 79 | 0 | 4 | 34 | 42 | 5 | 47 | 60 | 52 | 56 | 25 | ∞ | 38 | 12 | 42 | -2 | 46 | 25 | 8 | -9 | 2 |
| 3 | 60 | 63 | 0 | 43 | 19 | 61 | 21 | 42 | 49 | 20 | 68 | 27 | 84 | 10 | 74 | 13 | 4 | 81 | ∞ | 63 | 3 |
| 4 | 85 | 79 | 29 | 0 | 58 | 84 | 83 | 20 | ∞ | 50 | 22 | 67 | 78 | 85 | 9 | 54 | 9 | 53 | 68 | 24 | 4 |
| 5 | 2 | 45 | 22 | -9 | 0 | 80 | ∞ | -2 | 61 | -1 | 16 | 13 | 58 | 72 | 10 | 30 | 35 | -15 | 41 | 16 | 5 |
| 6 | 76 | 72 | 8 | 47 | 48 | 0 | -3 | -1 | 53 | 21 | 49 | -3 | 50 | 26 | 60 | ∞ | 18 | -2 | 9 | 72 | 6 |
| 7 | ∞ | 70 | 77 | 87 | 91 | 54 | 0 | -3 | 64 | 55 | 85 | 95 | 24 | 20 | 6 | 31 | 40 | 29 | 78 | 5 | 7 |
| 8 | 51 | ∞ | 24 | 23 | 77 | 66 | 51 | 0 | 48 | 47 | 69 | 48 | 58 | 55 | 58 | 29 | 91 | 64 | 70 | 33 | 8 |
| 9 | 60 | 13 | 84 | 71 | 62 | 78 | 57 | 86 | 0 | 18 | 60 | 93 | ∞ | 86 | 37 | 73 | 57 | 15 | 43 | -2 | 9 |
| 10 | 82 | 82 | 24 | 63 | 89 | 75 | 70 | 36 | -6 | 0 | -5 | 55 | 62 | 76 | 30 | 44 | ∞ | 53 | -8 | 53 | 10 |
| 11 | 27 | 32 | ∞ | 35 | 66 | 12 | 72 | -8 | -1 | 71 | 0 | -5 | 13 | 11 | 47 | 18 | 4 | 66 | 62 | -1 | 11 |
| 12 | 78 | 88 | 76 | 67 | 64 | 24 | 5 | 42 | 63 | ∞ | 6 | 0 | 38 | 77 | 61 | 38 | 43 | 41 | 51 | 44 | 12 |
| 13 | 0 | 66 | 30 | 1 | 52 | -8 | 5 | 4 | 61 | 52 | -11 | 67 | 0 | ∞ | 77 | 65 | 72 | 39 | 78 | -6 | 13 |
| 14 | 12 | 33 | 4 | -15 | 76 | 40 | 37 | 65 | 71 | 35 | 7 | 54 | 9 | 0 | 18 | 77 | 7 | ∞ | 70 | 72 | 14 |
| 15 | 44 | 63 | 80 | 61 | 21 | 32 | 47 | 66 | 24 | 34 | 73 | 35 | 72 | 41 | 0 | 36 | 62 | 26 | 47 | ∞ | 15 |
| 16 | 83 | 12 | 57 | 64 | 56 | 15 | 76 | 8 | 55 | 7 | 52 | 0 | 53 | 54 | ∞ | 0 | 43 | 15 | 72 | 43 | 16 |
| 17 | 82 | 23 | 18 | ∞ | 70 | 30 | 7 | 45 | 58 | 74 | 5 | 47 | 27 | 16 | 74 | 82 | 0 | 27 | 24 | 81 | 17 |
| 18 | 58 | 42 | 73 | 59 | ∞ | 42 | 17 | 22 | 89 | 81 | 32 | 41 | 74 | 33 | 17 | 71 | 71 | 0 | 75 | -6 | 18 |
| 19 | 20 | 11 | 57 | 26 | 53 | ∞ | 52 | 47 | 38 | 57 | 62 | 83 | 58 | 81 | 51 | 42 | 82 | 55 | 0 | -6 | 19 |
| 20 | 8 | 50 | 18 | 85 | 79 | 78 | 43 | 37 | 45 | 74 | ∞ | 0 | 29 | 59 | 49 | 46 | 12 | 41 | 67 | 0 | 20 |
|    | 1 | 2 | 3 | 4 | 5 | 6 | 7 | 8 | 9 | 10 | 11 | 12 | 13 | 14 | 15 | 16 | 17 | 18 | 19 | 20 |    |

We now present the negative paths in $D_8^{-1}M^-(20)$ followed by $P_{20}$.

j = 4.
(14, 4)(4, 15): -6;  (14, 4)(4, 17): -6.
j = 6.

(13, 6)(6, 7): -11;  (13, 6)(6, 8): -9;  (13, 6)(6, 12): -11;  (13, 6)(6, 18): -10.
j = 7.
(6, 7)(7, 8): -6;  (13, 7)(7, 8): -14;  (13, 7)(7, 15): -5;  (13, 7)(7, 20): -6.
j = 9.
(10, 9)(9, 20): -8;  (11, 9)(9, 20): -3.
j = 10.
(5, 10)(10, 9): -7;  (5, 10)(10, 11): -6;  (5, 10)(10, 19): -9.
j = 11.
(5, 11)(11, 8): -14;  (5, 11)(11, 9): -7;  (5, 11)(11, 12): -11;
(5, 11)(11, 20): -7; (10, 11)(11, 8): -13;  (10, 11)(11, 9): -6;  (10, 11)(11, 12): -10;
(10, 11)(11, 20): -8; (13, 11)(11, 8): -19; (13, 11)(11, 9): -12;
(13, 11)(11, 12): -16;  (13, 11)(11, 20): -14.
j = 12.
(5, 12)(12, 7): -6;  (5, 12)(12, 11): -5;  (10, 12)(12, 7):  -5;  (10, 12)(12, 7): -5;
(10, 12)(12, 11): -4;  (13, 12)(12, 7): -11;
(13, 12)(12, 11): -10.
j = 16.
(2, 16)(16, 12): -2.
j = 17.
(14, 17)(17, 11): -1.
j = 18.
(5, 18)(18, 20): -21;  (13, 18)(18, 20): -16.
j = 19.
(5, 19)(19, 20): -15;  (10, 19)(19. 20): -14.
j = 20.
(2, 20)(20, 1): -9;  (2, 20)(20, 12): -9;  (5, 20)(20, 1): -13;  (5, 20)(20, 12): -21;
(9, 20)(20, 12): -2;  (10, 20)(20, 1): -6;  (10, 20)(20, 12): -14;
(10, 20)(20, 17): -2;  (5, 20)(20, 17): -9; (11, 20)(20, 12): -1;
(13, 20)(20, 12): -6; (18, 20)(20, 12): -6; (19, 20)(20, 12): -6.

$$D_8^{-1} M(20)$$

| | 7 | 8 | 11 | 17 | 18 | 19 | 5 | 1 | 4 | 12 | 20 | 2 | 9 | 13 | 16 | 6 | 10 | 14 | 3 | 15 | |
|---|---|---|---|---|---|---|---|---|---|---|---|---|---|---|---|---|---|---|---|---|---|
| | 1 | 2 | 3 | 4 | 5 | 6 | 7 | 8 | 9 | 10 | 11 | 12 | 13 | 14 | 15 | 16 | 17 | 18 | 19 | 20 | |
| 1 | 0 | 50 | 43 | -7 | 63 | 77 | 5 | ∞ | 88 | 41 | 56 | 79 | 30 | 20 | 56 | 29 | 37 | 8 | 63 | 39 | 1 |
| 2 | -1 | 0 | 4 | 34 | 42 | 5 | 47 | 60 | 52 | 56 | 25 | -9 | 38 | 12 | 42 | -2 | 46 | 25 | 8 | -9 | 2 |
| 3 | 60 | 63 | 0 | 43 | 19 | 61 | 21 | 42 | 49 | 20 | 68 | 27 | 84 | 10 | 74 | 13 | 4 | 81 | ∞ | 63 | 3 |
| 4 | 85 | 79 | 29 | 0 | 58 | 84 | 83 | 20 | ∞ | 50 | 22 | 67 | 78 | 85 | 9 | 54 | 9 | 53 | 68 | 24 | 4 |
| 5 | -13 | 45 | -3 | -9 | 0 | 80 | -6 | -14 | -7 | -1 | -6 | -21 | 58 | 72 | 10 | 30 | 35 | -15 | -9 | -21 | 5 |
| 6 | 76 | 72 | 8 | 47 | 48 | 0 | -3 | -6 | 53 | 21 | 49 | -3 | 50 | 26 | 60 | ∞ | 18 | -2 | 9 | 72 | 6 |
| 7 | ∞ | 70 | 77 | 87 | 91 | 54 | 0 | -3 | 64 | 55 | 85 | 95 | 24 | 20 | 6 | 31 | 40 | 29 | 78 | 5 | 7 |
| 8 | 51 | ∞ | 24 | 23 | 77 | 66 | 51 | 0 | 48 | 47 | 69 | 48 | 58 | 55 | 58 | 29 | 91 | 64 | 70 | 33 | 8 |
| 9 | 60 | 13 | 84 | 71 | 62 | 78 | 57 | 86 | 0 | 18 | 60 | -2 | ∞ | 86 | 37 | 73 | 57 | 15 | 43 | -2 | 9 |
| 10 | -6 | 82 | 24 | 63 | 89 | 75 | -5 | -13 | -6 | 0 | -5 | -14 | 62 | 76 | 30 | 44 | ∞ | 53 | -8 | -14 | 10 |
| 11 | 27 | 32 | ∞ | 35 | 66 | 12 | 72 | -8 | -1 | 71 | 0 | -5 | 13 | 11 | 47 | 18 | 4 | 66 | 62 | -3 | 11 |
| 12 | 78 | 88 | 76 | 67 | 64 | 24 | 5 | 42 | 63 | ∞ | 6 | 0 | 38 | 77 | 61 | 38 | 43 | 41 | 51 | 44 | 12 |
| 13 | 0 | 66 | 30 | 1 | 52 | -8 | -11 | -19 | -12 | 52 | -11 | -16 | 0 | -6 | -5 | 65 | 72 | -10 | 78 | -16 | 13 |
| 14 | 12 | 33 | 4 | -15 | 76 | 40 | 37 | 65 | 71 | 35 | 7 | 54 | 9 | 0 | -6 | 77 | -6 | ∞ | 70 | 72 | 14 |
| 15 | 44 | 63 | 80 | 61 | 21 | 32 | 47 | 66 | 24 | 34 | 73 | 35 | 72 | 41 | 0 | 36 | 62 | 26 | 47 | ∞ | 15 |
| 16 | 83 | 12 | 57 | 64 | 56 | 15 | 76 | 8 | 55 | 7 | 52 | 0 | 53 | 54 | ∞ | 0 | 43 | 15 | 72 | 43 | 16 |
| 17 | 82 | 23 | 18 | ∞ | 70 | 30 | 7 | 45 | 58 | 74 | 5 | 47 | 27 | 16 | 74 | 82 | 0 | 27 | 24 | 81 | 17 |
| 18 | 58 | 42 | 73 | 59 | ∞ | 42 | 17 | 22 | 89 | 81 | 32 | -6 | 74 | 33 | 17 | 71 | 71 | 0 | 75 | -6 | 18 |
| 19 | 20 | 11 | 57 | 26 | 53 | ∞ | 52 | 47 | 38 | 57 | 62 | -6 | 58 | 81 | 51 | 42 | 82 | 55 | 0 | -6 | 19 |
| 20 | 8 | 50 | 18 | 85 | 79 | 78 | 43 | 37 | 45 | 74 | ∞ | 0 | 29 | 59 | 49 | 46 | 12 | 41 | 67 | 0 | 20 |
| | 1 | 2 | 3 | 4 | 5 | 6 | 7 | 8 | 9 | 10 | 11 | 12 | 13 | 14 | 15 | 16 | 17 | 18 | 19 | 20 | |

**P$_{20}$**

|     | 1  | 2 | 3 | 4 | 5 | 6 | 7    | 8  | 9     | 10 | 11 | 12    | 13 | 14 | 15 | 16 | 17 | 18 | 19 | 20 |     |
|-----|----|---|---|---|---|---|------|----|-------|----|----|-------|----|----|----|----|----|----|----|----|-----|
| 1   |    |   |   |   |   |   |      |    |       |    |    |       |    |    |    |    |    |    |    |    | 1   |
| 2   | 20 |   |   |   |   |   |      |    |       |    | 12 | 20    |    |    |    |    |    |    |    |    | 2   |
| 3   |    |   |   |   |   |   |      |    |       |    |    |       |    |    |    |    |    |    |    |    | 3   |
| 4   |    |   |   |   |   |   |      |    |       |    |    |       |    |    |    |    |    |    |    |    | 4   |
| 5   | 20 |   |   |   |   |   | 12   | 11 | 10,11 |    | 10 | 20    |    |    |    |    | 20 |    | 10 | 18 | 5   |
| 6   |    |   |   |   |   |   |      | 7  |       |    |    |       |    |    |    |    |    |    |    |    | 6   |
| 7   |    |   |   |   |   |   |      |    |       |    |    |       |    |    |    |    |    |    |    |    | 7   |
| 8   |    |   |   |   |   |   |      |    |       |    |    |       |    |    |    |    |    |    |    |    | 8   |
| 9   |    |   |   |   |   |   |      |    |       |    |    | 20    |    |    |    |    |    |    |    |    | 9   |
| 10  | 20 |   |   |   |   |   | 12   | 11 |       |    | 12 | 20    |    |    | 4  |    | 20 |    |    | 19 | 10  |
| 11  |    |   |   |   |   |   |      |    |       |    |    | 20    |    |    |    |    |    |    | 9  |    | 11  |
| 12  |    |   |   |   |   |   |      |    |       |    |    |       |    |    |    |    |    |    |    |    | 12  |
| 13  |    |   |   |   |   |   | 6,12 | 11 | 11    |    |    | 11,20 |    |    | 7  |    |    | 6  |    | 18 | 13  |
| 14  |    |   |   |   |   |   |      |    |       |    | 17 |       |    |    | 4  |    | 4  |    |    |    | 14  |
| 15  |    |   |   |   |   |   |      |    |       |    |    |       |    |    |    |    |    |    |    |    | 15  |
| 16  |    |   |   |   |   |   |      |    |       |    |    |       |    |    |    |    |    |    |    |    | 16  |
| 17  |    |   |   |   |   |   |      |    |       |    |    |       |    |    |    |    |    |    |    |    | 17  |
| 18  |    |   |   |   |   |   | 12   |    |       |    |    | 20    |    |    |    |    |    |    |    |    | 18  |
| 19  |    |   |   |   |   |   | 12   |    |       |    |    | 20    |    |    |    |    |    |    |    |    | 19  |
| 20  |    |   |   |   |   |   |      |    |       |    |    |       |    |    |    |    |    |    |    |    | 20  |
|     | 1  | 2 | 3 | 4 | 5 | 6 | 7    | 8  | 9     | 10 | 11 | 12    | 13 | 14 | 15 | 16 | 17 | 18 | 19 | 20 |     |

Before going on, we note that unlike the usual procedure in the regular F-W algorithm, when we obtain two or more paths which extend paths given in $D_8^{-1}M^-(20i)$, if both new paths have the *same* value, we include *all* of them in $P_{20i}$. It really doesn't affect the outcome of the algorithm except that we may be able to obtain a negative cycle that has a smaller value than the one we would obtain if we arbitrarily used only the first better value that occurred. The reasoning is as follows: As long as both paths, say

$$P(1) = [a, \ldots, x, \ldots, d], \quad P(2) = [a, \ldots, d]$$

are obtained in the same iteration, P(1) may have a point, x, which is repeated, while x is not reached by P(2). Since a point of P(1) has been repeated where P(1) consists of arcs which have always yielded a shortest path as its has passed through points in N = {1,2, ... , n}, a subpath of P(1) beginning and ending at x is a negative cycle, C. From Theorem 5, we can obtain C independently of P(1) in at least one more iteration. On the other hand, P(2) doesn't contain x and (assuming that it goes through no repeated points) may be extended to obtain a negative cycle whose terminal point is t = a. We thus have the option of strictly following the regular F-W algorithm with respect to negative paths (which is the simpler approach) or using the optional procedure just presented.

We now present the negative paths that are extensions of paths in $D_8^{-1}M^-(20)$. It is followed by $D_8^{-1}M^-(40)$ and $P_{40}$. Each path given must be of smaller value then the previous path at that entry.

j = 1.
(2, 1)(1, 4): -8;  (5, 1)(1, 4): -20;  (5, 1)(1, 7): -8;  (10, 1)(1, 4): -13; (10, 1)(1, 7): -1.
j = 4.
(10, 4)(4, 15): -4;  (10, 4)(4, 17): -4.
j = 12.
(2, 12)(12, 7): -4;  (2, 12)(12, 11): -3;  (5, 12)(12, 11): -15; (10, 12)(12, 11): - 8; (18, 12)(12, 7): -1; (19, 12)(12, 7): -1.
j =17.
(14, 17)(17, 11): -1.

## $D_8^{-1} M(40)$

| | 7 | 8 | 11 | 17 | 18 | 19 | 5 | 1 | 4 | 12 | 20 | 2 | 9 | 13 | 16 | 6 | 10 | 14 | 3 | 15 | |
|---|---|---|---|---|---|---|---|---|---|---|---|---|---|---|---|---|---|---|---|---|---|
| | 1 | 2 | 3 | 4 | 5 | 6 | 7 | 8 | 9 | 10 | 11 | 12 | 13 | 14 | 15 | 16 | 17 | 18 | 19 | 20 | |
| 1 | 0 | 50 | 43 | -7 | 63 | 77 | 5 | ∞ | 88 | 41 | 56 | 79 | 30 | 20 | 56 | 29 | 37 | 8 | 63 | 39 | 1 |
| 2 | -1 | 0 | 4 | -8 | 42 | 5 | -4 | 60 | 52 | 56 | -3 | -9 | 38 | 12 | 42 | -2 | 46 | 25 | 8 | -9 | 2 |
| 3 | 60 | 63 | 0 | 43 | 19 | 61 | 21 | 42 | 49 | 20 | 68 | 27 | 84 | 10 | 74 | 13 | 4 | 81 | ∞ | 63 | 3 |
| 4 | 85 | 79 | 29 | 0 | 58 | 84 | 83 | 20 | ∞ | 50 | 22 | 22 | 78 | 85 | 9 | 54 | 9 | 53 | 68 | 24 | 4 |
| 5 | -13 | 45 | -3 | -20 | 0 | 80 | -8 | -14 | -7 | -1 | -15 | -21 | 58 | 72 | 10 | 30 | -9 | -15 | -9 | -21 | 5 |
| 6 | 76 | 72 | 8 | 47 | 48 | 0 | -3 | -6 | 53 | 21 | 49 | -3 | 50 | 26 | 60 | ∞ | 18 | -2 | 9 | 72 | 6 |
| 7 | ∞ | 70 | 77 | 87 | 91 | 54 | 0 | -3 | 64 | 55 | 85 | 95 | 24 | 20 | 6 | 31 | 40 | 29 | 78 | 5 | 7 |
| 8 | 51 | ∞ | 24 | 23 | 77 | 66 | 51 | 0 | 48 | 47 | 69 | 48 | 58 | 55 | 58 | 29 | 91 | 64 | 70 | 33 | 8 |
| 9 | 60 | 13 | 84 | 71 | 62 | 78 | 57 | 86 | 0 | 18 | 60 | -2 | ∞ | 86 | 37 | 73 | 57 | 15 | 43 | -2 | 9 |
| 10 | -6 | 82 | 24 | -13 | 89 | 75 | -9 | -13 | -6 | 0 | -8 | -14 | 62 | 76 | -4 | 44 | -4 | 53 | -8 | -14 | 10 |
| 11 | 27 | 32 | ∞ | 35 | 66 | 12 | 72 | -8 | -1 | 71 | 0 | -5 | 13 | 11 | 47 | 18 | 4 | 66 | 62 | -3 | 11 |
| 12 | 78 | 88 | 76 | 67 | 64 | 24 | 5 | 42 | 63 | ∞ | 6 | 0 | 38 | 77 | 61 | 38 | 43 | 41 | 51 | 44 | 12 |
| 13 | -8 | 66 | 30 | 1 | 52 | -8 | -11 | -19 | -12 | 52 | -11 | -16 | 0 | -6 | -5 | 65 | 72 | -10 | 78 | -16 | 13 |
| 14 | 12 | 33 | 4 | -15 | 76 | 40 | 37 | 65 | 71 | 35 | -1 | 54 | 9 | 0 | -6 | 77 | -6 | ∞ | 70 | 72 | 14 |
| 15 | 44 | 63 | 80 | 61 | 21 | 32 | 47 | 66 | 24 | 34 | 73 | 35 | 72 | 41 | 0 | 36 | 62 | 26 | 47 | ∞ | 15 |
| 16 | 83 | 12 | 57 | 64 | 56 | 15 | 76 | 8 | 55 | 7 | 52 | 0 | 53 | 54 | ∞ | 0 | 43 | 15 | 72 | 43 | 16 |
| 17 | 82 | 23 | 18 | ∞ | 70 | 30 | 7 | 45 | 58 | 74 | 5 | 47 | 27 | 16 | 74 | 82 | 0 | 27 | 24 | 81 | 17 |
| 18 | 58 | 42 | 73 | 59 | ∞ | 42 | -1 | 22 | 89 | 81 | 32 | -6 | 74 | 33 | 17 | 71 | 71 | 0 | 75 | -6 | 18 |
| 19 | 20 | 11 | 57 | 26 | 53 | ∞ | -1 | 47 | 38 | 57 | 62 | -6 | 58 | 81 | 51 | 42 | 82 | 55 | 0 | -6 | 19 |
| 20 | 8 | 50 | 18 | 85 | 79 | 78 | 43 | 37 | 45 | 74 | ∞ | 0 | 29 | 59 | 49 | 46 | 12 | 41 | 67 | 0 | 20 |
| | 1 | 2 | 3 | 4 | 5 | 6 | 7 | 8 | 9 | 10 | 11 | 12 | 13 | 14 | 15 | 16 | 17 | 18 | 19 | 20 | |

$P_{40}$

|    | 1  | 2 | 3 | 4 | 5 | 6 | 7    | 8  | 9  | 10 | 11 | 12    | 13 | 14 | 15 | 16 | 17 | 18 | 19 | 20 |    |
|----|----|---|---|---|---|---|------|----|----|----|----|-------|----|----|----|----|----|----|----|----|----|
| 1  |    |   |   |   |   |   |      |    |    |    |    |       |    |    |    |    |    |    |    |    | 1  |
| 2  | 20 |   |   | 1 |   |   | 12   |    |    |    | 12 | 20    |    |    |    |    |    |    |    |    | 2  |
| 3  |    |   |   |   |   |   |      |    |    |    |    |       |    |    |    |    |    |    |    |    | 3  |
| 4  |    |   |   |   |   |   |      |    |    |    |    |       |    |    |    |    |    |    |    |    | 4  |
| 5  | 20 |   |   | 1 |   |   | 1    | 11 | 11 |    | 12 | 20    |    |    |    |    | 20 |    | 10 | 18 | 5  |
| 6  |    |   |   |   |   |   |      | 7  |    |    |    |       |    |    |    |    |    |    |    |    | 6  |
| 7  |    |   |   |   |   |   |      |    |    |    |    |       |    |    |    |    |    |    |    |    | 7  |
| 8  |    |   |   |   |   |   |      |    |    |    |    |       |    |    |    |    |    |    |    |    | 8  |
| 9  |    |   |   |   |   |   |      |    |    |    |    | 20    |    |    |    |    |    |    |    |    | 9  |
| 10 | 20 |   |   | 1 |   |   | 12   | 11 |    |    | 12 | 20    |    |    | 4  |    | 4  |    |    | 19 | 10 |
| 11 |    |   |   |   |   |   |      |    |    |    |    | 20    |    |    |    |    |    |    |    | 9  | 11 |
| 12 |    |   |   |   |   |   |      |    |    |    |    |       |    |    |    |    |    |    |    |    | 12 |
| 13 |    |   |   |   |   |   | 6,12 | 11 | 11 |    |    | 11,20 |    |    | 7  |    |    | 6  |    | 18 | 13 |
| 14 |    |   |   |   |   |   |      |    |    |    | 17 |       |    |    | 4  |    | 4  |    |    |    | 14 |
| 15 |    |   |   |   |   |   |      |    |    |    |    |       |    |    |    |    |    |    |    |    | 15 |
| 16 |    |   |   |   |   |   |      |    |    |    |    |       |    |    |    |    |    |    |    |    | 16 |
| 17 |    |   |   |   |   |   |      |    |    |    |    |       |    |    |    |    |    |    |    |    | 17 |
| 18 |    |   |   |   |   |   | 12   |    |    |    |    | 20    |    |    |    |    |    |    |    |    | 18 |
| 19 |    |   |   |   |   |   | 12   |    |    |    |    | 20    |    |    |    |    |    |    |    |    | 19 |
| 20 |    |   |   |   |   |   |      |    |    |    |    |       |    |    |    |    |    |    |    |    | 20 |
|    | 1  | 2 | 3 | 4 | 5 | 6 | 7    | 8  | 9  | 10 | 11 | 12    | 13 | 14 | 15 | 16 | 17 | 18 | 19 | 20 |    |

We now present the extensions of paths in $D_8^{-1}M^-(40)$. It is followed by $D_8^{-1}M^-(60)$ and $P_{60}$.

i = 1.
(13, 1)(1, 4): -15.
j = 4.
(5, 4)(4, 15): -11;  (5, 4)(4, 17): -11.
j = 7.
(18, 7)(7, 8): -4.
j = 11.
(2, 11)(11, 8): -11;  (2, 11)(11, 9): -4;  (5, 11)(11, 8): -23;
(5, 11)(11, 9): -16;  (5, 11)(11, 12): -20;  ((5, 11)(11, 13): -2;
(5, 11)(11, 14): -4;  (5, 11)(11, 17): -11;  (10, 11)(11, 8): -16;
(10, 11)(11, 9): -9;  (14, 11)(11, 8): -9;  (14, 11)(11, 9): -2;
(14, 11)(11, 12): -6;  (14, 11)(11, 20): -2.

$$D_8^{-1}M(60)$$

|   | 7 | 8 | 11 | 17 | 18 | 19 | 5 | 1 | 4 | 12 | 20 | 2 | 9 | 13 | 16 | 6 | 10 | 14 | 3 | 15 |   |
|---|---|---|---|---|---|---|---|---|---|---|---|---|---|---|---|---|---|---|---|---|---|
|   | 1 | 2 | 3 | 4 | 5 | 6 | 7 | 8 | 9 | 10 | 11 | 12 | 13 | 14 | 15 | 16 | 17 | 18 | 19 | 20 |   |
| 1 | 0 | 50 | 43 | *-7* | 63 | 77 | 5 | ∞ | 88 | 41 | 56 | 79 | 30 | 20 | 56 | 29 | 37 | 8 | 63 | 39 | 1 |
| 2 | 79 | 0 | 4 | *-8* | 42 | 5 | *-4* | <u>-11</u> | <u>-4</u> | 56 | 25 | ∞ | 38 | 12 | 42 | *-2* | 46 | 25 | 8 | *-9* | 2 |
| 3 | 60 | 63 | 0 | 43 | 19 | 61 | 21 | 42 | 49 | 20 | 68 | 27 | 84 | 10 | 74 | 13 | 4 | 81 | ∞ | 63 | 3 |
| 4 | 85 | 79 | 29 | 0 | 58 | 84 | 83 | 20 | ∞ | 50 | 22 | 22 | 78 | 85 | 9 | 54 | 9 | 53 | 68 | 24 | 4 |
| 5 | *-13* | 45 | *-3* | *-20* | 0 | <u>-3</u> | *-8* | <u>-23</u> | <u>-16</u> | -1 | -15 | -21 | <u>-2</u> | <u>-4</u> | <u>-11</u> | 30 | <u>-11</u> | -15 | *-9* | *-21* | 5 |
| 6 | 76 | 72 | 8 | 47 | 48 | 0 | -3 | *-6* | 53 | 21 | 49 | *-3* | 50 | 26 | 60 | ∞ | 18 | *-2* | 9 | 72 | 6 |
| 7 | ∞ | 70 | 77 | 87 | 91 | 54 | 0 | *-3* | 64 | 55 | 85 | 95 | 24 | 20 | 6 | 31 | 40 | 29 | 78 | 5 | 7 |
| 8 | 51 | ∞ | 24 | 23 | 77 | 66 | 51 | 0 | 48 | 47 | 69 | 48 | 58 | 55 | 58 | 29 | 91 | 64 | 70 | 33 | 8 |
| 9 | 60 | 13 | 84 | 71 | 62 | 78 | 57 | 86 | 0 | 18 | 60 | *-2* | ∞ | 86 | 37 | 73 | 57 | 15 | 43 | *-2* | 9 |
| 10 | *-6* | 82 | 24 | *-13* | 89 | 75 | *-9* | *-13* | -6 | 0 | *-8* | *-14* | 62 | 76 | *-4* | 44 | *-4* | 53 | *-8* | *-14* | 10 |
| 11 | 27 | 32 | ∞ | 35 | 66 | 12 | 72 | *-8* | -1 | 71 | 0 | *-5* | 13 | 11 | 47 | 18 | 4 | 66 | 62 | *-3* | 11 |
| 12 | 78 | 88 | 76 | 67 | 64 | 24 | 5 | 42 | 63 | ∞ | 6 | 0 | 38 | 77 | 61 | 38 | 43 | 41 | 51 | 44 | 12 |
| 13 | *-8* | 66 | 30 | <u>-15</u> | 52 | *-8* | *-11* | *-19* | *-12* | 52 | -11 | *-16* | 0 | *-6* | *-5* | 65 | 72 | *-10* | 78 | *-16* | 13 |
| 14 | 12 | 33 | 4 | *-15* | 76 | 40 | 37 | 65 | 71 | 35 | *-1* | 54 | 9 | 0 | *-6* | 77 | *-6* | ∞ | 70 | <u>-2</u> | 14 |
| 15 | 44 | 63 | 80 | 61 | 21 | 32 | 47 | 66 | 24 | 34 | 73 | 35 | 72 | 41 | 0 | 36 | 62 | 26 | 47 | ∞ | 15 |
| 16 | 83 | 12 | 57 | 64 | 56 | 15 | 76 | 8 | 55 | 7 | 52 | 0 | 53 | 54 | ∞ | 0 | 43 | 15 | 72 | 43 | 16 |
| 17 | 82 | 23 | 18 | ∞ | 70 | 30 | 7 | 45 | 58 | 74 | 5 | 47 | 27 | 16 | 74 | 82 | 0 | 27 | 24 | 81 | 17 |
| 18 | 58 | 42 | 73 | 59 | ∞ | 42 | *-1* | <u>-4</u> | 89 | 81 | 32 | *-6* | 74 | 33 | 17 | 71 | 71 | 0 | 75 | *-6* | 18 |
| 19 | 20 | 11 | 57 | 26 | 53 | ∞ | *-1* | <u>-4</u> | 38 | 57 | 62 | *-6* | 58 | 81 | 51 | 42 | 82 | 55 | 0 | *-6* | 19 |
| 20 | 8 | 50 | 18 | 85 | 79 | 78 | 43 | 37 | 45 | 74 | ∞ | 0 | 29 | 59 | 49 | 46 | 12 | 41 | 67 | 0 | 20 |
|   | 1 | 2 | 3 | 4 | 5 | 6 | 7 | 8 | 9 | 10 | 11 | 12 | 13 | 14 | 15 | 16 | 17 | 18 | 19 | 20 |   |

**P₆₀**

|    | 1  | 2 | 3 | 4 | 5 | 6 | 7    | 8  | 9     | 10 | 11 | 12    | 13 | 14 | 15 | 16 | 17 | 18 | 19 | 20 |    |
|----|----|---|---|---|---|---|------|----|-------|----|----|-------|----|----|----|----|----|----|----|----|----|
| 1  |    |   |   |   |   |   |      |    |       |    |    |       |    |    |    |    |    |    |    |    | 1  |
| 2  | 20 |   |   | 1 |   |   | 12   | 11 | 11    |    |    | 12    | 20 |    |    |    |    |    |    |    | 2  |
| 3  |    |   |   |   |   |   |      |    |       |    |    |       |    |    |    |    |    |    |    |    | 3  |
| 4  |    |   |   |   |   |   |      |    |       |    |    |       |    |    |    |    |    |    |    |    | 4  |
| 5  | 20 |   |   | 1 |   |   | 1    | 11 | 10,11 |    |    | 12    | 20 |    | 11 | 11 | 4  |    | 4  |    | 10 | 18 | 5  |
| 6  |    |   |   |   |   |   |      | 7  |       |    |    |       |    |    |    |    |    |    |    |    | 6  |
| 7  |    |   |   |   |   |   |      |    |       |    |    |       |    |    |    |    |    |    |    |    | 7  |
| 8  |    |   |   |   |   |   |      |    |       |    |    |       |    |    |    |    |    |    |    |    | 8  |
| 9  |    |   |   |   |   |   |      |    |       |    |    | 20    |    |    |    |    |    |    |    |    | 9  |
| 10 | 20 |   |   | 1 |   |   | 12   | 11 | 11    |    |    | 12    | 20 |    |    | 4  |    | 4  |    |    | 19 | 10 |
| 11 |    |   |   |   |   |   |      |    |       |    |    | 20    |    |    |    |    |    |    |    | 9  | 11 |
| 12 |    |   |   |   |   |   |      |    |       |    |    |       |    |    |    |    |    |    |    |    | 12 |
| 13 |    |   |   | 1 |   |   | 6,12 | 11 | 11    |    |    | 11,20 |    |    | 7  |    |    | 6  |    | 18 | 13 |
| 14 |    |   |   |   |   |   |      | 11 |       |    | 17 |       |    |    | 4  |    | 4  |    |    | 11 | 14 |
| 15 |    |   |   |   |   |   |      |    |       |    |    |       |    |    |    |    |    |    |    |    | 15 |
| 16 |    |   |   |   |   |   |      |    |       |    |    |       |    |    |    |    |    |    |    |    | 16 |
| 17 |    |   |   |   |   |   |      |    |       |    |    |       |    |    |    |    |    |    |    |    | 17 |
| 18 |    |   |   |   |   |   | 12   | 7  |       |    |    | 20    |    |    |    |    |    |    |    |    | 18 |
| 19 |    |   |   |   |   |   | 12   |    |       |    |    | 20    |    |    |    |    |    |    |    |    | 19 |
| 20 |    |   |   |   |   |   |      |    |       |    |    |       |    |    |    |    |    |    |    |    | 20 |
|    | 1  | 2 | 3 | 4 | 5 | 6 | 7    | 8  | 9     | 10 | 11 | 12    | 13 | 14 | 15 | 16 | 17 | 18 | 19 | 20 |    |

$P_{60}$

|    | 1  | 2 | 3 | 4 | 5 | 6 | 7    | 8  | 9     | 10 | 11 | 12    | 13 | 14 | 15 | 16 | 17 | 18 | 19 | 20 |    |
|----|----|---|---|---|---|---|------|----|-------|----|----|-------|----|----|----|----|----|----|----|----|----|
| 1  |    |   |   |   |   |   |      |    |       |    |    |       |    |    |    |    |    |    |    |    | 1  |
| 2  | 20 |   |   | 1 |   |   | 12   | 11 | 11    |    | 12 | 20    |    |    |    |    |    |    |    |    | 2  |
| 3  |    |   |   |   |   |   |      |    |       |    |    |       |    |    |    |    |    |    |    |    | 3  |
| 4  |    |   |   |   |   |   |      |    |       |    |    |       |    |    |    |    |    |    |    |    | 4  |
| 5  | 20 |   |   | 1 |   |   | 1    | 11 | 10,11 |    | 12 | 20    | 11 | 11 | 4  |    | 4  |    | 10 | 18 | 5  |
| 6  |    |   |   |   |   |   |      | 7  |       |    |    |       |    |    |    |    |    |    |    |    | 6  |
| 7  |    |   |   |   |   |   |      |    |       |    |    |       |    |    |    |    |    |    |    |    | 7  |
| 8  |    |   |   |   |   |   |      |    |       |    |    |       |    |    |    |    |    |    |    |    | 8  |
| 9  |    |   |   |   |   |   |      |    |       |    |    | 20    |    |    |    |    |    |    |    |    | 9  |
| 10 | 20 |   |   | 1 |   |   | 12   | 11 | 11    |    | 12 | 20    |    |    | 4  |    | 4  |    |    | 19 | 10 |
| 11 |    |   |   |   |   |   |      |    |       |    |    | 20    |    |    |    |    |    |    |    | 9  | 11 |
| 12 |    |   |   |   |   |   |      |    |       |    |    |       |    |    |    |    |    |    |    |    | 12 |
| 13 |    |   |   | 1 |   |   | 6,12 | 11 | 11    |    |    | 11,20 |    |    | 7  |    |    | 6  |    | 18 | 13 |
| 14 |    |   |   |   |   |   |      | 11 |       | 17 |    |       |    |    | 4  |    | 4  |    |    | 11 | 14 |
| 15 |    |   |   |   |   |   |      |    |       |    |    |       |    |    |    |    |    |    |    |    | 15 |
| 16 |    |   |   |   |   |   |      |    |       |    |    |       |    |    |    |    |    |    |    |    | 16 |
| 17 |    |   |   |   |   |   |      |    |       |    |    |       |    |    |    |    |    |    |    |    | 17 |
| 18 |    |   |   |   |   |   | 12   | 7  |       |    |    | 20    |    |    |    |    |    |    |    |    | 18 |
| 19 |    |   |   |   |   |   | 12   |    |       |    |    | 20    |    |    |    |    |    |    |    |    | 19 |
| 20 |    |   |   |   |   |   |      |    |       |    |    |       |    |    |    |    |    |    |    |    | 20 |
|    | 1  | 2 | 3 | 4 | 5 | 6 | 7    | 8  | 9     | 10 | 11 | 12    | 13 | 14 | 15 | 16 | 17 | 18 | 19 | 20 |    |

The following are extensions of negative paths in $D_8^{-1}M^-(60)$.

j = 9.
(5, 9)(9, 2): -3.
j = 4.
(13, 4)(4, 15): -6;  (13, 4)(4, 17): -6.
j = 20.

**(14, 20)(20, 11): -2.**

$$D_8^{-1} M^-(80)$$

|   | 7 | 8 | 11 | 17 | 18 | 19 | 5 | 1 | 4 | 12 | 20 | 2 | 9 | 13 | 16 | 6 | 10 | 14 | 3 | 15 |   |
|---|---|---|---|---|---|---|---|---|---|---|---|---|---|---|---|---|---|---|---|---|---|
|   | 1 | 2 | 3 | 4 | 5 | 6 | 7 | 8 | 9 | 10 | 11 | 12 | 13 | 14 | 15 | 16 | 17 | 18 | 19 | 20 |   |
| 1 | 0 | 50 | 43 | *-7* | 63 | 77 | 5 | ∞ | 88 | 41 | 56 | 79 | 30 | 20 | 56 | 29 | 37 | 8 | 63 | 39 | 1 |
| 2 | 79 | 0 | 4 | *-8* | 42 | 5 | *-4* | *-11* | *-4* | 56 | *-3* | ∞ | 38 | 12 | 42 | *-2* | 46 | 25 | 8 | *-9* | 2 |
| 3 | 60 | 63 | 0 | 43 | 19 | 61 | 21 | 42 | 49 | 20 | 68 | 27 | 84 | 10 | 74 | 13 | 4 | 81 | ∞ | 63 | 3 |
| 4 | 85 | 79 | 29 | 0 | 58 | 84 | 83 | 20 | ∞ | 50 | 22 | 22 | 78 | 85 | 9 | 54 | 9 | 53 | 68 | 24 | 4 |
| 5 | *-13* | <u>*-3*</u> | *-3* | *-20* | 0 | *-3* | *-8* | *-23* | *-16* | *-1* | *-15* | *-21* | *-2* | *-4* | *-11* | 30 | *-11* | *-15* | 41 | *-21* | 5 |
| 6 | 76 | 72 | 8 | 47 | 48 | 0 | *-3* | *-6* | 53 | 21 | 49 | *-3* | 50 | 26 | 60 | ∞ | 18 | *-2* | 9 | 72 | 6 |
| 7 | ∞ | 70 | 77 | 87 | 91 | 54 | 0 | *-3* | 64 | 55 | 85 | 95 | 24 | 20 | 6 | 31 | 40 | 29 | 78 | 5 | 7 |
| 8 | 51 | ∞ | 24 | 23 | 77 | 66 | 51 | 0 | 48 | 47 | 69 | 48 | 58 | 55 | 58 | 29 | 91 | 64 | 70 | 33 | 8 |
| 9 | 60 | 13 | 84 | 71 | 62 | 78 | 57 | 86 | 0 | 18 | 60 | *-2* | ∞ | 86 | 37 | 73 | 57 | 15 | 43 | *-2* | 9 |
| 10 | *-6* | 82 | 24 | *-13* | 89 | 75 | *-9* | *-13* | *-6* | 0 | *-8* | *-14* | 62 | 76 | *-4* | 44 | *-4* | 53 | *-8* | *-14* | 10 |
| 11 | 27 | 32 | ∞ | 35 | 66 | 12 | 72 | *-8* | *-1* | 71 | 0 | *-5* | 13 | 11 | 47 | 18 | 4 | 66 | 62 | *-3* | 11 |
| 12 | 78 | 88 | 76 | 67 | 64 | 24 | 5 | 42 | 63 | ∞ | 6 | 0 | 38 | 77 | 61 | 38 | 43 | 41 | 51 | 44 | 12 |
| 13 | *-8* | 66 | 30 | *-15* | 52 | *-8* | *-11* | *-19* | *-12* | 52 | *-11* | 67 | 0 | ∞ | <u>*-6*</u> | 65 | <u>*-6*</u> | *-10* | 78 | *-16* | 13 |
| 14 | 12 | 33 | 4 | *-15* | 76 | 40 | 37 | 65 | 71 | 35 | *-1* | <u>*-2*</u> | 9 | 0 | *-6* | 77 | *-6* | ∞ | 70 | 72 | 14 |
| 15 | 44 | 63 | 80 | 61 | 21 | 32 | 47 | 66 | 24 | 34 | 73 | 35 | 72 | 41 | 0 | 36 | 62 | 26 | 47 | ∞ | 15 |
| 16 | 83 | 12 | 57 | 64 | 56 | 15 | 76 | 8 | 55 | 7 | 52 | 0 | 53 | 54 | ∞ | 0 | 43 | 15 | 72 | 43 | 16 |
| 17 | 82 | 23 | 18 | ∞ | 70 | 30 | 7 | 45 | 58 | 74 | 5 | 47 | 27 | 16 | 74 | 82 | 0 | 27 | 24 | 81 | 17 |
| 18 | 58 | 42 | 73 | 59 | ∞ | 42 | *-1* | 22 | 89 | 81 | 32 | *-6* | 74 | 33 | 17 | 71 | 71 | 0 | 75 | *-6* | 18 |
| 19 | 20 | 11 | 57 | 26 | 53 | ∞ | 52 | 47 | 38 | 57 | 62 | *-6* | 58 | 81 | 51 | 42 | 82 | 55 | 0 | *-6* | 19 |
| 20 | 8 | 50 | 18 | 85 | 79 | 78 | 43 | 37 | 45 | 74 | ∞ | 0 | 29 | 59 | 49 | 46 | 12 | 41 | 67 | 0 | 20 |
|   | 1 | 2 | 3 | 4 | 5 | 6 | 7 | 8 | 9 | 10 | 11 | 12 | 13 | 14 | 15 | 16 | 17 | 18 | 19 | 20 |   |

**P$_{80}$**

|   | 1 | 2 | 3 | 4 | 5 | 6 | 7 | 8 | 9 | 10 | 11 | 12 | 13 | 14 | 15 | 16 | 17 | 18 | 19 | 20 |   |
|---|---|---|---|---|---|---|---|---|---|----|----|----|----|----|----|----|----|----|----|----|---|
| 1 |   |   |   |   |   |   |   |   |   |    |    |    |    |    |    |    |    |    |    |    | 1 |
| 2 | 20 |  |   | 1 |   |   | 12 | 11 | 11 |   |   | 12 | 20 |   |   |   |   |   |   |   | 2 |
| 3 |   |   |   |   |   |   |   |   |   |    |    |    |    |    |    |    |    |    |    |    | 3 |
| 4 |   |   |   |   |   |   |   |   |   |    |    |    |    |    |    |    |    |    |    |    | 4 |
| 5 | 20 | 9 |   | 1 |   |   | 1 | 11 | 10,11 |   |   | 12 | 20 | 11 | 11 | 4 |   | 4 |   | 10 | 18 | 5 |
| 6 |   |   |   |   |   |   |   | 7 |   |    |    |    |    |    |    |    |    |    |    |    | 6 |
| 7 |   |   |   |   |   |   |   |   |   |    |    |    |    |    |    |    |    |    |    |    | 7 |
| 8 |   |   |   |   |   |   |   |   |   |    |    |    |    |    |    |    |    |    |    |    | 8 |
| 9 |   |   |   |   |   |   |   |   |   |    |    | 20 |    |    |    |    |    |    |    |    | 9 |
| 10 | 20 |  |   | 1 |   |   | 12 | 11 | 11 |   |   | 12 | 20 |   |   | 4 |   | 4 |   | 19 | 10 |
| 11 |   |   |   |   |   |   |   |   |   |    |    | 20 |    |    |    |    |    |    |    | 9 | 11 |
| 12 |   |   |   |   |   |   |   |   |   |    |    |    |    |    |    |    |    |    |    |    | 12 |
| 13 |   |   |   | 1 |   |   | 6,12 | 11 | 11 |   |   | 11,20 |   |   | 7 |   |   | 6 |   | 18 | 13 |
| 14 |   |   |   |   |   |   |   | 11 |   | 20 |   |   |   |   | 4 |   | 4 |   |   | 11 | 14 |
| 15 |   |   |   |   |   |   |   |   |   |    |    |    |    |    |    |    |    |    |    |    | 15 |
| 16 |   |   |   |   |   |   |   |   |   |    |    |    |    |    |    |    |    |    |    |    | 16 |
| 17 |   |   |   |   |   |   |   |   |   |    |    |    |    |    |    |    |    |    |    |    | 17 |
| 18 |   |   |   |   |   |   | 12 | 7 |   |    |    | 20 |    |    |    |    |    |    |    |    | 18 |
| 19 |   |   |   |   |   |   | 12 |   |   |    |    | 20 |    |    |    |    |    |    |    |    | 19 |
| 20 |   |   |   |   |   |   |   |   |   |    |    |    |    |    |    |    |    |    |    |    | 20 |
|   | 1 | 2 | 3 | 4 | 5 | 6 | 7 | 8 | 9 | 10 | 11 | 12 | 13 | 14 | 15 | 16 | 17 | 18 | 19 | 20 |   |

The only path that can be extended is at (5, 2). (5, 2)(2, 16): − 5. Adding the value -5 of the entry at (5, 16) to each of the entries in row 16 yields no negative sum in the corresponding column's entry in row 5 which is smaller than the corresponding entry in row 5. It follows that we can no longer obtain a derangement whose value is less than that of D$_8$. Therefore,

$\sigma_{APOPT} = D_8 = (1\ 7\ 5\ 18\ 14\ 13\ 9\ 4\ 17\ 10\ 12\ 2\ 8)(3\ 11\ 20\ 15\ 16\ 6\ 19)$.

**PHASE 3.** Here we use the corollary to Theorem 1 in order to find non-negative cycles in $D_8^{-1}M^-$ that have a value no greater than a fixed positive number. As a last resort, we will add up the first log n numbers in each row and find their mean value, m. We will then use $M_1 = |\sigma_{APOPT}| + m$ as a first upper bound for $|\sigma_{TSPOPT}|$. If we can found no set of disjoint cycles whose sum of values is no greater than m such that the action of the cycles on $\sigma_{APOPT}$ leads to an n-cycle, we define $M_2 = |\sigma_{APOPT}| + 2m$. Using the cycles we've already obtained, we obtain all cycles whose sum is no greater than 2m. We then use all cycles obtained to try to obtain a disjoint set of cycles whose sum is less than 2m whose action on $\sigma_{APOPT}$ leads to an n-cycle, etc. . We always try to obtain disjoint cycles with smallest non-negative values. We may first try 0, than 1, than 2, ... , etc. . Since we know that $D_8^{-1}M^-$ contains no negative cycles, we use the F-W bounded-value subpath algorithm (F - W bvs algorithm) restricted to those cycles whose value is less than a fixed non-negative number. From the corollary to Theorem 1, we know that we need only choose as initial value, i, an entry whose value is no greater than the sum of the values in the cycle, i. One of the advantages of our procedure for finding $\sigma_{APOPT}$ is that one or more of the derangements obtained may be an n-cycle. In Phase 1, we always start with a randomly chosen n-cycle. Thus, if we repeat the first part of PHASE 1 n(log n) times, the probability of obtaining at least one n-cycle is

$$1 - (1 - \frac{e}{n})^{n(\log n)} \approx 1 - \frac{1}{n^e} \to 1$$

for large n. Of course, there is no guarantee that the n-cycle obtained would be close in value to $|\sigma_{TSPOPT}|$. In the case given here, we are in luck. $|D_7| = 213$, while $|D_8| = 212$. Thus, either $D_7$ is an optimal solution to the TSP or else there exists a set of disjoint cycles each of value 0 in $D_8^{-1}M^-$ which when acted upon $D_8$ yields an n-cycle. In Phase 3(a), we *approximate* a solution to the TSP by using the F-W n-nvs algorithm to obtain non-negative cycles such that some subset that is disjoint forms a permutation, s, such that $D_i s$ is an n-cycle. In checking to see if $D_8$ is an optimal solution for the assignment problem, we obtain a large number of negative paths which may be used as *part* of the search for non-negative cycles. However, it is best to use the F-W n-nvs algorithm so that we include a large number of possible cases. My guess is that the determining arc for most of the cases will generally be negative. If the upper bound is small, we may use Phase 3(b) to search for *all* cycles less than a fixed bound.  The procedure used is to construct tree-like structures that contain no cycles. Each branch is treated as a separate entity. Thus, we may have the same vertex on different branches. Theorems 1 and 6 reduce the number of branches obtained.

Going on to Phase 3(a), we can use the result that $D_8 = \sigma_{APOPT}$ to help show that the tour of length 213 is probably $\sigma_{TSPOPT}$. Our reason for equivocation is the following: we are using the F-W n-nvs algorithm that allows us to get at least some of the cycles of value 0 (if any exist) such that some product of disjoint cycles, s, has the property that $\sigma_{APOPT}s$ is a tour. If that is the case, $\sigma_{APOPT}s = \sigma_{TSPOPT}$. We will use $D_8^{-1}M^-(20i)$ (i = 1, 2, ... , ) suitably revised to include paths of value 0.
We now present the non-positive paths in $D_8^{-1}M^-(20)$ followed by $P_{20}$.

j = 1.
(13, 1)(1, 5): -4.
j = 4.
(5, 4)(4, 15): 0;  (5, 4)(4, 17): 0;  (14, 4)(4, 15): -6;  (14, 4)(4, 17): -6.
j = 6.
(13, 6)(6, 7): -11;  (13, 6)(6, 8): -9;  (13, 6)(6, 12): -11;  (13, 6)(6, 18): -10.
j = 7.
(6, 7)(7, 8): -6;  (13, 7)(7, 8): -14;  (13, 7)(7, 15): -5;  (13, 7)(7, 20): -6.
j = 9.
(10, 9)(9, 20): -8;  (11, 9)(9, 20): -3.
j = 10.
(5, 10)(10, 9): -7;  (5, 10)(10, 11): -6;  (5, 10)(10, 19): -9.
j = 11.
(5, 11)(11, 8): -14;  (5, 11)(11, 9): -7;  (5, 11)(11, 12): -11;
(5, 11)(11, 20): -7; (10, 11)(11, 8): -13;  (10, 11)(11, 9): -6;
(10, 11)(11, 12): -10; (10, 11)(11, 20): -8; (13, 11)(11, 8): -19;
(13, 11)(11, 9): -12;  (13, 11)(11, 12): -16;  (13, 11)(11, 20): -14.
j = 12.
(5, 12)(12, 7): -6;  (5, 12)(12, 11): -5;  (10, 12)(12, 7):  -5;
(10, 12)(12, 7): -5;  (10, 12)(12, 11): -4;  (13, 12)(12, 7): -11;
(13, 12)(12, 11): -10.
j = 16.
(2, 16)(16, 12): -2.
j = 17.
(14, 17)(17, 11): -1.
j = 18.
(5, 18)(18, 20): -21;  (13, 18)(18, 20): -16.
j = 19.
(5, 19)(19, 20): -15;  (10, 19)(19. 20): -14.
j = 20.
(2, 20)(20, 1): -9;  (2, 20)(20, 12): -9;  (5, 20)(20, 1): -13;
(5, 20)(20, 12): -21;  (9, 20)(20, 12): -2;  (10, 20)(20, 1): -6;
(10, 20)(20, 12): -14;  (10, 20)(20, 17): -2;  (5, 20)(20, 17): -9;
(11, 20)(20, 12): -1;  (13, 20)(20, 12): -6; (18, 20)(20, 12): -6;
(19, 20)(20, 12): -6.

$$D_8^{-1}M(20)$$

|    | 7 | 8 | 11 | 17 | 18 | 19 | 5 | 1 | 4 | 12 | 20 | 2 | 9 | 13 | 16 | 6 | 10 | 14 | 3 | 15 |    |
|----|---|---|----|----|----|----|---|---|---|----|----|---|---|----|----|---|----|----|---|----|----|
|    | 1 | 2 | 3  | 4  | 5  | 6  | 7 | 8 | 9 | 10 | 11 | 12| 13| 14 | 15 | 16| 17 | 18 | 19| 20 |    |
| 1  | 0 | 50 | 43 | -7 | 63 | 77 | 5 | ∞ | 88 | 41 | 56 | 79 | 30 | 20 | 56 | 29 | 37 | 8 | 63 | 39 | 1 |
| 2  | -1 | 0 | 4 | 34 | 42 | 5 | 47 | 60 | 52 | 56 | 25 | -9 | 38 | 12 | 42 | -2 | 46 | 25 | 8 | -9 | 2 |
| 3  | 60 | 63 | 0 | 43 | 19 | 61 | 21 | 42 | 49 | 20 | 68 | 27 | 84 | 10 | 74 | 13 | 4 | 81 | ∞ | 63 | 3 |
| 4  | 85 | 79 | 29 | 0 | 58 | 84 | 83 | 20 | ∞ | 50 | 22 | 67 | 78 | 85 | 9 | 54 | 9 | 53 | 68 | 24 | 4 |
| 5  | -13 | 45 | -3 | -9 | 0 | 80 | -6 | -14 | -7 | -1 | -6 | -21 | 58 | 72 | 0 | 30 | 0 | -15 | -9 | -21 | 5 |
| 6  | 76 | 72 | 8 | 47 | 48 | 0 | -3 | -6 | 53 | 21 | 49 | -3 | 50 | 26 | 60 | ∞ | 18 | -2 | 9 | 72 | 6 |
| 7  | ∞ | 70 | 77 | 87 | 91 | 54 | 0 | -3 | 64 | 55 | 85 | 95 | 24 | 20 | 6 | 31 | 40 | 29 | 78 | 5 | 7 |
| 8  | 51 | ∞ | 24 | 23 | 77 | 66 | 51 | 0 | 48 | 47 | 69 | 48 | 58 | 55 | 58 | 29 | 91 | 64 | 70 | 33 | 8 |
| 9  | 60 | 13 | 84 | 71 | 62 | 78 | 57 | 86 | 0 | 18 | 60 | -2 | ∞ | 86 | 37 | 73 | 57 | 15 | 43 | -2 | 9 |
| 10 | -6 | 82 | 24 | 63 | 89 | 75 | -5 | -13 | -6 | 0 | -5 | -14 | 62 | 76 | 30 | 44 | ∞ | 53 | -8 | -14 | 10 |
| 11 | 27 | 32 | ∞ | 35 | 66 | 12 | 72 | -8 | -1 | 71 | 0 | -5 | 13 | 11 | 47 | 18 | 4 | 66 | 62 | -3 | 11 |
| 12 | 78 | 88 | 76 | 67 | 64 | 24 | 5 | 42 | 63 | ∞ | 6 | 0 | 38 | 77 | 61 | 38 | 43 | 41 | 51 | 44 | 12 |
| 13 | 0 | 66 | 30 | 1 | 52 | -8 | -11 | -19 | -12 | 52 | -11 | -16 | 0 | -6 | -5 | 65 | 72 | -10 | 78 | -16 | 13 |
| 14 | 12 | 33 | 4 | -15 | 76 | 40 | 37 | 65 | 71 | 35 | 7 | 54 | 9 | 0 | -6 | 77 | -6 | ∞ | 70 | 72 | 14 |
| 15 | 44 | 63 | 80 | 61 | 21 | 32 | 47 | 66 | 24 | 34 | 73 | 35 | 72 | 41 | 0 | 36 | 62 | 26 | 47 | ∞ | 15 |
| 16 | 83 | 12 | 57 | 64 | 56 | 15 | 76 | 8 | 55 | 7 | 52 | 0 | 53 | 54 | ∞ | 0 | 43 | 15 | 72 | 43 | 16 |
| 17 | 82 | 23 | 18 | ∞ | 70 | 30 | 7 | 45 | 58 | 74 | 5 | 47 | 27 | 16 | 74 | 82 | 0 | 27 | 24 | 81 | 17 |
| 18 | 58 | 42 | 73 | 59 | ∞ | 42 | 17 | 22 | 89 | 81 | 32 | -6 | 74 | 33 | 17 | 71 | 71 | 0 | 75 | -6 | 18 |
| 19 | 20 | 11 | 57 | 26 | 53 | ∞ | 52 | 47 | 38 | 57 | 62 | -6 | 58 | 81 | 51 | 42 | 82 | 55 | 0 | -6 | 19 |
| 20 | 8 | 50 | 18 | 85 | 79 | 78 | 43 | 37 | 45 | 74 | ∞ | 0 | 29 | 59 | 49 | 46 | 12 | 41 | 67 | 0 | 20 |
|    | 1 | 2 | 3 | 4 | 5 | 6 | 7 | 8 | 9 | 10 | 11 | 12 | 13 | 14 | 15 | 16 | 17 | 18 | 19 | 20 |    |

$P_{20}$

|    | 1  | 2 | 3 | 4 | 5 | 6 | 7    | 8  | 9     | 10 | 11 | 12    | 13 | 14 | 15 | 16 | 17 | 18 | 19 | 20 |    |
|----|----|---|---|---|---|---|------|----|-------|----|----|-------|----|----|----|----|----|----|----|----|----|
| 1  |    |   |   |   |   |   |      |    |       |    |    |       |    |    |    |    |    |    |    |    | 1  |
| 2  | 20 |   |   |   |   |   |      |    |       |    | 12 | 20    |    |    |    |    |    |    |    |    | 2  |
| 3  |    |   |   |   |   |   |      |    |       |    |    |       |    |    |    |    |    |    |    |    | 3  |
| 4  |    |   |   |   |   |   |      |    |       |    |    |       |    |    |    |    |    |    |    |    | 4  |
| 5  | 20 |   |   |   |   |   | 12   | 11 | 10,11 |    | 10 | 20    |    |    | 4  |    | 20 |    | 10 | 18 | 5  |
| 6  |    |   |   |   |   |   |      | 7  |       |    |    |       |    |    |    |    |    |    |    |    | 6  |
| 7  |    |   |   |   |   |   |      |    |       |    |    |       |    |    |    |    |    |    |    |    | 7  |
| 8  |    |   |   |   |   |   |      |    |       |    |    |       |    |    |    |    |    |    |    |    | 8  |
| 9  |    |   |   |   |   |   |      |    |       |    |    | 20    |    |    |    |    |    |    |    |    | 9  |
| 10 | 20 |   |   |   |   |   | 12   | 11 |       |    | 12 | 20    |    |    | 4  |    | 20 |    |    | 19 | 10 |
| 11 |    |   |   |   |   |   |      |    |       |    |    | 20    |    |    |    |    |    |    |    | 9  | 11 |
| 12 |    |   |   |   |   |   |      |    |       |    |    |       |    |    |    |    |    |    |    |    | 12 |
| 13 |    |   |   |   |   |   | 6,12 | 11 | 11    |    |    | 11,20 |    |    | 7  |    |    | 6  |    | 18 | 13 |
| 14 |    |   |   |   |   |   |      |    |       |    | 17 |       |    |    | 4  |    | 4  |    |    |    | 14 |
| 15 |    |   |   |   |   |   |      |    |       |    |    |       |    |    |    |    |    |    |    |    | 15 |
| 16 |    |   |   |   |   |   |      |    |       |    |    |       |    |    |    |    |    |    |    |    | 16 |
| 17 |    |   |   |   |   |   |      |    |       |    |    |       |    |    |    |    |    |    |    |    | 17 |
| 18 |    |   |   |   |   |   | 12   |    |       |    |    | 20    |    |    |    |    |    |    |    |    | 18 |
| 19 |    |   |   |   |   |   | 12   |    |       |    |    | 20    |    |    |    |    |    |    |    |    | 19 |
| 20 |    |   |   |   |   |   |      |    |       |    |    |       |    |    |    |    |    |    |    |    | 20 |
|    | 1  | 2 | 3 | 4 | 5 | 6 | 7    | 8  | 9     | 10 | 11 | 12    | 13 | 14 | 15 | 16 | 17 | 18 | 19 | 20 |    |

We now present the non-positive paths that are extensions of paths in $D_8^{-1}M^-(20)$ followed by $D_8^{-1}M^-(40)$ and $P_{40}$. Each path given must be of smaller value then the previous path at that entry.

j = 1.
(2, 1)(1, 4): -8;  (5, 1)(1, 4): -20;  (5, 1)(1, 7): -8;  (10, 1)(1, 4): -13;
(10, 1)(1, 7): -1.
j = 11.
(5, 11)(11, 17): -2.
j = 12.

(2, 12)(12, 7): -4;  (2, 12)(12, 11): -3;  (5, 12)(12, 7): -16;  (5, 12)(12, 11): -15; (10, 12)(12, 11): - 8;  (18, 12)(12, 7): -1; (19, 12)(12, 7): -1.
j =17.
(14, 17)(17, 11): -1.
j = 20.
(10, 20)(20, 17): -2.

$$D_8^{-1} M^- (40)$$

| | 7 | 8 | 11 | 17 | 18 | 19 | 5 | 1 | 4 | 12 | 20 | 2 | 9 | 13 | 16 | 6 | 10 | 14 | 3 | 15 | |
|---|---|---|---|---|---|---|---|---|---|---|---|---|---|---|---|---|---|---|---|---|---|
| | 1 | 2 | 3 | 4 | 5 | 6 | 7 | 8 | 9 | 10 | 11 | 12 | 13 | 14 | 15 | 16 | 17 | 18 | 19 | 20 | |
| 1 | 0 | 50 | 43 | -7 | 63 | 77 | 5 | ∞ | 88 | 41 | 56 | 79 | 30 | 20 | 56 | 29 | 37 | 8 | 63 | 39 | 1 |
| 2 | *-1* | 0 | 4 | *-8* | 42 | 5 | *-4* | 60 | 52 | 56 | *-3* | *-9* | 38 | 12 | 42 | -2 | 46 | 25 | 8 | -9 | 2 |
| 3 | 60 | 63 | 0 | 43 | 19 | 61 | 21 | 42 | 49 | 20 | 68 | 27 | 84 | 10 | 74 | 13 | 4 | 81 | ∞ | 63 | 3 |
| 4 | 85 | 79 | 29 | 0 | 58 | 84 | 83 | 20 | ∞ | 50 | 22 | 22 | 78 | 85 | 9 | 54 | 9 | 53 | 68 | 24 | 4 |
| 5 | *-13* | 45 | *-3* | *-20* | 0 | 80 | *-16* | -14 | -7 | -1 | *-15* | *-21* | 58 | 72 | 10 | 30 | *-9* | -15 | *-9* | *-21* | 5 |
| 6 | 76 | 72 | 8 | 47 | 48 | 0 | -3 | *-6* | 53 | 21 | 49 | -3 | 50 | 26 | 60 | ∞ | 18 | -2 | 9 | 72 | 6 |
| 7 | ∞ | 70 | 77 | 87 | 91 | 54 | 0 | -3 | 64 | 55 | 85 | 95 | 24 | 20 | 6 | 31 | 40 | 29 | 78 | 5 | 7 |
| 8 | 51 | ∞ | 24 | 23 | 77 | 66 | 51 | 0 | 48 | 47 | 69 | 48 | 58 | 55 | 58 | 29 | 91 | 64 | 70 | 33 | 8 |
| 9 | 60 | 13 | 84 | 71 | 62 | 78 | 57 | 86 | 0 | 18 | 60 | *-2* | ∞ | 86 | 37 | 73 | 57 | 15 | 43 | -2 | 9 |
| 10 | *-6* | 82 | 24 | *-13* | 89 | 75 | *-9* | *-13* | -6 | 0 | *-8* | *-14* | 62 | 76 | 30 | 44 | *-2* | 53 | -8 | *-14* | 10 |
| 11 | 27 | 32 | ∞ | 35 | 66 | 12 | 72 | -8 | -1 | 71 | 0 | -5 | 13 | 11 | 47 | 18 | 4 | 66 | 62 | *-3* | 11 |
| 12 | 78 | 88 | 76 | 67 | 64 | 24 | 5 | 42 | 63 | ∞ | 6 | 0 | 38 | 77 | 61 | 38 | 43 | 41 | 51 | 44 | 12 |
| 13 | *-8* | 66 | 30 | 1 | 52 | -8 | *-11* | *-19* | *-12* | 52 | *-11* | *-16* | 0 | -6 | -5 | 65 | 72 | *-10* | 78 | *-16* | 13 |
| 14 | 12 | 33 | 4 | *-15* | 76 | 40 | 37 | 65 | 71 | 35 | *-1* | 54 | 9 | 0 | *-6* | 77 | *-6* | ∞ | 70 | 72 | 14 |
| 15 | 44 | 63 | 80 | 61 | 21 | 32 | 47 | 66 | 24 | 34 | 73 | 35 | 72 | 41 | 0 | 36 | 62 | 26 | 47 | ∞ | 15 |
| 16 | 83 | 12 | 57 | 64 | 56 | 15 | 76 | 8 | 55 | 7 | 52 | 0 | 53 | 54 | ∞ | 0 | 43 | 15 | 72 | 43 | 16 |
| 17 | 82 | 23 | 18 | ∞ | 70 | 30 | 7 | 45 | 58 | 74 | 5 | 47 | 27 | 16 | 74 | 82 | 0 | 27 | 24 | 81 | 17 |
| 18 | 58 | 42 | 73 | 59 | ∞ | 42 | *-1* | 22 | 89 | 81 | 32 | *-6* | 74 | 33 | 17 | 71 | 71 | 0 | 75 | *-6* | 18 |
| 19 | 20 | 11 | 57 | 26 | 53 | ∞ | *-1* | 47 | 38 | 57 | 62 | *-6* | 58 | 81 | 51 | 42 | 82 | 55 | 0 | *-6* | 19 |
| 20 | 8 | 50 | 18 | 85 | 79 | 78 | 43 | 37 | 45 | 74 | ∞ | 0 | 29 | 59 | 49 | 46 | 12 | 41 | 67 | 0 | 20 |
| | 1 | 2 | 3 | 4 | 5 | 6 | 7 | 8 | 9 | 10 | 11 | 12 | 13 | 14 | 15 | 16 | 17 | 18 | 19 | 20 | |

$P_{40}$

|    | 1  | 2 | 3 | 4 | 5 | 6 | 7    | 8  | 9  | 10 | 11 | 12    | 13 | 14 | 15 | 16 | 17 | 18 | 19 | 20 |    |
|----|----|---|---|---|---|---|------|----|----|----|----|-------|----|----|----|----|----|----|----|----|----|
| 1  |    |   |   |   |   |   |      |    |    |    |    |       |    |    |    |    |    |    |    |    | 1  |
| 2  | 20 |   |   | 1 |   |   | 12   |    |    |    | 12 | 20    |    |    |    |    |    |    |    |    | 2  |
| 3  |    |   |   |   |   |   |      |    |    |    |    |       |    |    |    |    |    |    |    |    | 3  |
| 4  |    |   |   |   |   |   |      |    |    |    |    |       |    |    |    |    |    |    |    |    | 4  |
| 5  | 20 |   |   | 1 |   |   | 1    | 11 | 11 |    | 12 | 20    |    |    |    |    | 20 |    | 10 | 18 | 5  |
| 6  |    |   |   |   |   |   |      | 7  |    |    |    |       |    |    |    |    |    |    |    |    | 6  |
| 7  |    |   |   |   |   |   |      |    |    |    |    |       |    |    |    |    |    |    |    |    | 7  |
| 8  |    |   |   |   |   |   |      |    |    |    |    |       |    |    |    |    |    |    |    |    | 8  |
| 9  |    |   |   |   |   |   |      |    |    |    |    | 20    |    |    |    |    |    |    |    |    | 9  |
| 10 | 20 |   |   | 1 |   |   | 12   | 11 |    |    | 12 | 20    |    |    |    |    | 1  |    |    | 19 | 10 |
| 11 |    |   |   |   |   |   |      |    |    |    |    | 20    |    |    |    |    |    |    |    | 9  | 11 |
| 12 |    |   |   |   |   |   |      |    |    |    |    |       |    |    |    |    |    |    |    |    | 12 |
| 13 |    |   |   |   |   |   | 6,12 | 11 | 11 |    |    | 11,20 |    |    | 7  |    |    | 6  |    | 18 | 13 |
| 14 |    |   |   |   |   |   |      |    |    |    | 17 |       |    |    | 4  |    | 4  |    |    |    | 14 |
| 15 |    |   |   |   |   |   |      |    |    |    |    |       |    |    |    |    |    |    |    |    | 15 |
| 16 |    |   |   |   |   |   |      |    |    |    |    |       |    |    |    |    |    |    |    |    | 16 |
| 17 |    |   |   |   |   |   |      |    |    |    |    |       |    |    |    |    |    |    |    |    | 17 |
| 18 |    |   |   |   |   |   | 12   |    |    |    |    | 20    |    |    |    |    |    |    |    |    | 18 |
| 19 |    |   |   |   |   |   | 12   |    |    |    |    | 20    |    |    |    |    |    |    |    |    | 19 |
| 20 |    |   |   |   |   |   |      |    |    |    |    |       |    |    |    |    |    |    |    |    | 20 |
|    | 1  | 2 | 3 | 4 | 5 | 6 | 7    | 8  | 9  | 10 | 11 | 12    | 13 | 14 | 15 | 16 | 17 | 18 | 19 | 20 |    |

**We now present the extensions of paths in $D_8^{-1}M^-(40)$. They are followed by $D_8^{-1}M^-(60)$ and $P_{60}$.**

i = 1.
(13, 1)(1, 4): -15.
j = 4.
(5, 4)(4, 15): -11;  (5, 4)(4, 17): -11; (10, 4)(4, 15): -4;  (10, 4)(4, 17): -4.
j = 7.
(18, 7)(7, 8): -4.
j = 11.
(2, 11)(11, 8): -11;  (2, 11)(11, 9): -4;  (2, 11)(11, 12): -16;
(2, 11)(11, 14): 0;  (2, 11)(11, 17): -7;  (2, 11)(11, 20): -12;
(5, 11)(11, 6): -3;  ((5, 11)(11, 8): -23;  (5, 11)(11, 9): -16;
(5, 11)(11, 13): -2;  (5, 11)(11, 14): -4;  (5, 11)(11, 17): -11;
(10, 11)(11, 8): -16;  (10, 11)(11, 9): -9;  (14, 11)(11, 8): -9;
(14, 11)(11, 9): -2;  (14, 11)(11, 12): -6;  (14, 11)(11, 20): -2.

$$D_8^{-1}M^-(60)$$

|    | 7  | 8  | 11 | 17  | 18 | 19 | 5   | 1   | 4   | 12 | 20  | 2   | 9  | 13 | 16  | 6  | 10  | 14 | 3   | 15  |    |
|----|----|----|----|-----|----|----|-----|-----|-----|----|-----|-----|----|----|-----|----|-----|----|-----|-----|----|
|    | 1  | 2  | 3  | 4   | 5  | 6  | 7   | 8   | 9   | 10 | 11  | 12  | 13 | 14 | 15  | 16 | 17  | 18 | 19  | 20  |    |
| 1  | 0  | 50 | 43 | -7  | 63 | 77 | 5   | ∞   | 88  | 41 | 56  | 79  | 30 | 20 | 56  | 29 | 37  | 8  | 63  | 39  | 1  |
| 2  | -1 | 0  | 4  | -8  | 42 | 5  | -4  | -11 | -4  | 56 | -3  | -16 | 38 | 0  | 42  | -2 | -7  | 25 | 8   | -12 | 2  |
| 3  | 60 | 63 | 0  | 43  | 19 | 61 | 21  | 42  | 49  | 20 | 68  | 27  | 84 | 10 | 74  | 13 | 4   | 81 | ∞   | 63  | 3  |
| 4  | 85 | 79 | 29 | 0   | 58 | 84 | 83  | 20  | ∞   | 50 | 22  | 22  | 78 | 85 | 9   | 54 | 9   | 53 | 68  | 24  | 4  |
| 5  | -13| 45 | -3 | -20 | 0  | -3 | -16 | -23 | -16 | -1 | -15 | -21 | -2 | -4 | -11 | 30 | -11 | -15| -9  | -21 | 5  |
| 6  | 76 | 72 | 8  | 47  | 48 | 0  | -3  | -6  | 53  | 21 | 49  | -3  | 50 | 26 | 60  | ∞  | 18  | -2 | 9   | 72  | 6  |
| 7  | ∞  | 70 | 77 | 87  | 91 | 54 | 0   | -3  | 64  | 55 | 85  | 95  | 24 | 20 | 6   | 31 | 40  | 29 | 78  | 5   | 7  |
| 8  | 51 | ∞  | 24 | 23  | 77 | 66 | 51  | 0   | 48  | 47 | 69  | 48  | 58 | 55 | 58  | 29 | 91  | 64 | 70  | 33  | 8  |
| 9  | 60 | 13 | 84 | 71  | 62 | 78 | 57  | 86  | 0   | 18 | 60  | -2  | ∞  | 86 | 37  | 73 | 57  | 15 | 43  | -2  | 9  |
| 10 | -6 | 82 | 24 | -13 | 89 | 75 | -9  | -16 | -9  | 0  | -8  | -14 | 62 | 76 | -4  | 44 | -4  | 53 | -8  | -14 | 10 |
| 11 | 27 | 32 | ∞  | 35  | 66 | 12 | 72  | -8  | -1  | 71 | 0   | -5  | 13 | 11 | 47  | 18 | 4   | 66 | 62  | -3  | 11 |
| 12 | 78 | 88 | 76 | 67  | 64 | 24 | 5   | 42  | 63  | ∞  | 6   | 0   | 38 | 77 | 61  | 38 | 43  | 41 | 51  | 44  | 12 |
| 13 | -8 | 66 | 30 | -13 | 52 | -8 | -11 | -19 | -12 | 52 | -11 | -16 | 0  | -6 | -5  | 65 | 72  | -10| 78  | -16 | 13 |
| 14 | 12 | 33 | 4  | -15 | 76 | 40 | 37  | -9  | -2  | 35 | -1  | -6  | 9  | 0  | -6  | 77 | -6  | ∞  | 70  | -2  | 14 |
| 15 | 44 | 63 | 80 | 61  | 21 | 32 | 47  | 66  | 24  | 34 | 73  | 35  | 72 | 41 | 0   | 36 | 62  | 26 | 47  | ∞   | 15 |
| 16 | 83 | 12 | 57 | 64  | 56 | 15 | 76  | 8   | 55  | 7  | 52  | 0   | 53 | 54 | ∞   | 0  | 43  | 15 | 72  | 43  | 16 |
| 17 | 82 | 23 | 18 | ∞   | 70 | 30 | 7   | 45  | 58  | 74 | 5   | 47  | 27 | 16 | 74  | 82 | 0   | 27 | 24  | 81  | 17 |
| 18 | 58 | 42 | 73 | 59  | ∞  | 42 | -1  | 22  | 89  | 81 | 32  | -6  | 74 | 33 | 17  | 71 | 71  | 0  | 75  | -6  | 18 |
| 19 | 20 | 11 | 57 | 26  | 53 | ∞  | -1  | 47  | 38  | 57 | 62  | -6  | 58 | 81 | 51  | 42 | 82  | 55 | 0   | -6  | 19 |
| 20 | 8  | 50 | 18 | 85  | 79 | 78 | 43  | 37  | 45  | 74 | ∞   | 0   | 29 | 59 | 49  | 46 | 12  | 41 | 67  | 0   | 20 |
|    | 1  | 2  | 3  | 4   | 5  | 6  | 7   | 8   | 9   | 10 | 11  | 12  | 13 | 14 | 15  | 16 | 17  | 18 | 19  | 20  |    |

**P₆₀**

|    | 1  | 2 | 3 | 4 | 5 | 6 | 7    | 8  | 9  | 10 | 11 | 12    | 13 | 14 | 15 | 16 | 17   | 18 | 19 | 20 |    |
|----|----|---|---|---|---|---|------|----|----|----|----|-------|----|----|----|----|------|----|----|----|----|
| 1  |    |   |   |   |   |   |      |    |    |    |    |       |    |    |    |    |      |    |    |    | 1  |
| 2  | 20 |   |   | 1 |   |   | 12   | 11 |    |    | 12 | 11    |    |    |    |    | 11   |    |    | 11 | 2  |
| 3  |    |   |   |   |   |   |      |    |    |    |    |       |    |    |    |    |      |    |    |    | 3  |
| 4  |    |   |   |   |   |   |      |    |    |    |    |       |    |    |    |    |      |    |    |    | 4  |
| 5  | 20 |   |   | 1 |   |   | 11   | 1  | 11 | 11 |    | 12    | 20 | 11 | 11 | 4  |      | 4,11 |  | 10 | 18 | 5  |
| 6  |    |   |   |   |   |   | 7    |    |    |    |    |       |    |    |    |    |      |    |    |    | 6  |
| 7  |    |   |   |   |   |   |      |    |    |    |    |       |    |    |    |    |      |    |    |    | 7  |
| 8  |    |   |   |   |   |   |      |    |    |    |    |       |    |    |    |    |      |    |    |    | 8  |
| 9  |    |   |   |   |   |   |      |    |    |    |    | 20    |    |    |    |    |      |    |    |    | 9  |
| 10 | 20 |   |   | 1 |   |   | 12   | 11 | 11 |    | 12 | 20    |    |    | 4  |    | 4,11 |    |    | 19 | 10 |
| 11 |    |   |   |   |   |   |      |    |    |    |    | 20    |    |    |    |    |      |    |    | 9  | 11 |
| 12 |    |   |   |   |   |   |      |    |    |    |    |       |    |    |    |    |      |    |    |    | 12 |
| 13 |    |   |   | 1 |   |   | 6,12 | 11 | 11 |    |    | 11,20 |    |    | 7  |    |      | 6  |    | 18 | 13 |
| 14 |    |   |   |   |   |   |      | 11 | 11 |    | 17 | 11    |    |    | 4  |    | 4    |    |    | 11 | 14 |
| 15 |    |   |   |   |   |   |      |    |    |    |    |       |    |    |    |    |      |    |    |    | 15 |
| 16 |    |   |   |   |   |   |      |    |    |    |    |       |    |    |    |    |      |    |    |    | 16 |
| 17 |    |   |   |   |   |   |      |    |    |    |    |       |    |    |    |    |      |    |    |    | 17 |
| 18 |    |   |   |   |   |   | 12   |    |    |    |    | 20    |    |    |    |    |      |    |    |    | 18 |
| 19 |    |   |   |   |   |   | 12   |    |    |    |    | 20    |    |    |    |    |      |    |    |    | 19 |
| 20 |    |   |   |   |   |   |      |    |    |    |    |       |    |    |    |    |      |    |    |    | 20 |
|    | 1  | 2 | 3 | 4 | 5 | 6 | 7    | 8  | 9  | 10 | 11 | 12    | 13 | 14 | 15 | 16 | 17   | 18 | 19 | 20 |    |

The following are extensions pf paths from $D_8^{-1}M^-(60)$ to $D_8^{-1}M^-(80)$.

$j = 9$.
(5, 9)(9, 2): -3;  (14, 9)(9, 20): -4.
$j = 12$.
(2, 12)(12, 7): -11;  (14, 12)(12, 7): -1.
$j = 14$.
(2, 14)(14, 4): -15.
$j = 20$.

(2, 20)(20, 1): -4.

$$D_8^{-1}M^-(80)$$

|   | 7 | 8 | 11 | 17 | 18 | 19 | 5 | 1 | 4 | 12 | 20 | 2 | 9 | 13 | 16 | 6 | 10 | 14 | 3 | 15 |   |
|---|---|---|---|---|---|---|---|---|---|---|---|---|---|---|---|---|---|---|---|---|---|
|   | 1 | 2 | 3 | 4 | 5 | 6 | 7 | 8 | 9 | 10 | 11 | 12 | 13 | 14 | 15 | 16 | 17 | 18 | 19 | 20 |   |
| 1 | 0 | 50 | 43 | -7 | 63 | 77 | 5 | ∞ | 88 | 41 | 56 | 79 | 30 | 20 | 56 | 29 | 37 | 8 | 63 | 39 | 1 |
| 2 | -4 | 0 | 4 | -15 | 42 | 5 | -11 | -11 | -4 | 56 | -3 | -16 | 38 | 0 | 42 | -2 | -7 | 25 | 8 | -12 | 2 |
| 3 | 60 | 63 | 0 | 43 | 19 | 61 | 21 | 42 | 49 | 20 | 68 | 27 | 84 | 10 | 74 | 13 | 4 | 81 | ∞ | 63 | 3 |
| 4 | 85 | 79 | 29 | 0 | 58 | 84 | 83 | 20 | ∞ | 50 | 22 | 22 | 78 | 85 | 9 | 54 | 9 | 53 | 68 | 24 | 4 |
| 5 | -13 | -3 | -3 | -20 | 0 | -3 | -16 | -23 | -16 | -1 | -15 | -21 | -2 | -4 | -11 | 30 | -11 | -15 | -9 | -21 | 5 |
| 6 | 76 | 72 | 8 | 47 | 48 | 0 | -3 | -6 | 53 | 21 | 49 | -3 | 50 | 26 | 60 | ∞ | 18 | -2 | 9 | 72 | 6 |
| 7 | ∞ | 70 | 77 | 87 | 91 | 54 | 0 | -3 | 64 | 55 | 85 | 95 | 24 | 20 | 6 | 31 | 40 | 29 | 78 | 5 | 7 |
| 8 | 51 | ∞ | 24 | 23 | 77 | 66 | 51 | 0 | 48 | 47 | 69 | 48 | 58 | 55 | 58 | 29 | 91 | 64 | 70 | 33 | 8 |
| 9 | 60 | 13 | 84 | 71 | 62 | 78 | 57 | 86 | 0 | 18 | 60 | -2 | ∞ | 86 | 37 | 73 | 57 | 15 | 43 | -2 | 9 |
| 10 | -6 | 82 | 24 | -13 | 89 | 75 | -9 | -16 | -9 | 0 | -8 | -14 | 62 | 76 | -4 | 44 | -4 | 53 | -8 | -14 | 10 |
| 11 | 27 | 32 | ∞ | 35 | 66 | 12 | 72 | -8 | -1 | 71 | 0 | -5 | 13 | 11 | 47 | 18 | 4 | 66 | 62 | -3 | 11 |
| 12 | 78 | 88 | 76 | 67 | 64 | 24 | 5 | 42 | 63 | ∞ | 6 | 0 | 38 | 77 | 61 | 38 | 43 | 41 | 51 | 44 | 12 |
| 13 | -8 | 66 | 30 | -13 | 52 | -8 | -11 | -19 | -12 | 52 | -11 | -16 | 0 | -6 | -5 | 65 | 72 | -10 | 78 | -16 | 13 |
| 14 | 12 | 33 | 4 | -15 | 76 | 40 | -1 | -9 | -2 | 35 | -1 | -6 | 9 | 0 | -6 | 77 | -6 | ∞ | 70 | -4 | 14 |
| 15 | 44 | 63 | 80 | 61 | 21 | 32 | 47 | 66 | 24 | 34 | 73 | 35 | 72 | 41 | 0 | 36 | 62 | 26 | 47 | ∞ | 15 |
| 16 | 83 | 12 | 57 | 64 | 56 | 15 | 76 | 8 | 55 | 7 | 52 | 0 | 53 | 54 | ∞ | 0 | 43 | 15 | 72 | 43 | 16 |
| 17 | 82 | 23 | 18 | ∞ | 70 | 30 | 7 | 45 | 58 | 74 | 5 | 47 | 27 | 16 | 74 | 82 | 0 | 27 | 24 | 81 | 17 |
| 18 | 58 | 42 | 73 | 59 | ∞ | 42 | -1 | 22 | 89 | 81 | 32 | -6 | 74 | 33 | 17 | 71 | 71 | 0 | 75 | -6 | 18 |
| 19 | 20 | 11 | 57 | 26 | 53 | ∞ | -1 | 47 | 38 | 57 | 62 | -6 | 58 | 81 | 51 | 42 | 82 | 55 | 0 | -6 | 19 |
| 20 | 8 | 50 | 18 | 85 | 79 | 78 | 43 | 37 | 45 | 74 | ∞ | 0 | 29 | 59 | 49 | 46 | 12 | 41 | 67 | 0 | 20 |
|   | 1 | 2 | 3 | 4 | 5 | 6 | 7 | 8 | 9 | 10 | 11 | 12 | 13 | 14 | 15 | 16 | 17 | 18 | 19 | 20 |   |

**P₈₀**

|    | 1  | 2 | 3 | 4  | 5 | 6  | 7    | 8  | 9  | 10 | 11 | 12    | 13 | 14 | 15 | 16 | 17   | 18 | 19 | 20 |    |
|----|----|---|---|----|---|----|------|----|----|----|----|-------|----|----|----|----|------|----|----|----|----|
| 1  |    |   |   |    |   |    |      |    |    |    |    |       |    |    |    |    |      |    |    |    | 1  |
| 2  | 20 |   |   | 14 |   |    | 12   | 11 |    |    | 12 | 11    |    |    |    |    | 11   |    |    | 11 | 2  |
| 3  |    |   |   |    |   |    |      |    |    |    |    |       |    |    |    |    |      |    |    |    | 3  |
| 4  |    |   |   |    |   |    |      |    |    |    |    |       |    |    |    |    |      |    |    |    | 4  |
| 5  | 20 | 9 |   | 1  |   | 11 | 1    | 11 | 11 |    | 12 | 20    | 11 | 11 | 4  |    | 4,11 |    | 10 | 18 | 5  |
| 6  |    |   |   |    |   |    |      | 7  |    |    |    |       |    |    |    |    |      |    |    |    | 6  |
| 7  |    |   |   |    |   |    |      |    |    |    |    |       |    |    |    |    |      |    |    |    | 7  |
| 8  |    |   |   |    |   |    |      |    |    |    |    |       |    |    |    |    |      |    |    |    | 8  |
| 9  |    |   |   |    |   |    |      |    |    |    |    | 20    |    |    |    |    |      |    |    |    | 9  |
| 10 | 20 |   |   | 1  |   |    | 12   | 11 | 11 |    | 12 | 20    |    |    | 4  |    | 4,11 |    |    | 19 | 10 |
| 11 |    |   |   |    |   |    |      |    |    |    |    | 20    |    |    |    |    |      |    |    | 9  | 11 |
| 12 |    |   |   |    |   |    |      |    |    |    |    |       |    |    |    |    |      |    |    |    | 12 |
| 13 |    |   |   | 1  |   |    | 6,12 | 11 | 11 |    |    | 11,20 |    |    | 7  |    |      | 6  |    | 18 | 13 |
| 14 |    |   |   |    |   |    | 12   | 11 | 11 |    | 17 | 11    |    |    | 4  |    | 4    |    |    | 9  | 14 |
| 15 |    |   |   |    |   |    |      |    |    |    |    |       |    |    |    |    |      |    |    |    | 15 |
| 16 |    |   |   |    |   |    |      |    |    |    |    |       |    |    |    |    |      |    |    |    | 16 |
| 17 |    |   |   |    |   |    |      |    |    |    |    |       |    |    |    |    |      |    |    |    | 17 |
| 18 |    |   |   |    |   |    | 12   |    |    |    |    | 20    |    |    |    |    |      |    |    |    | 18 |
| 19 |    |   |   |    |   |    | 12   |    |    |    |    | 20    |    |    |    |    |      |    |    |    | 19 |
| 20 |    |   |   |    |   |    |      |    |    |    |    |       |    |    |    |    |      |    |    |    | 20 |
|    | 1  | 2 | 3 | 4  | 5 | 6  | 7    | 8  | 9  | 10 | 11 | 12    | 13 | 14 | 15 | 16 | 17   | 18 | 19 | 20 |    |

The following are the extensions of paths in $D_8^{-1}M^-(80)$ :

j = 2.
(5, 2)(2, 16): -5.
j = 3.
(2, 4)(4, 15): -6.

j = 7.
(2, 7)(7, 8): -14.
We now consider each of these cases.
(5, 16) has but one extension - to (5, 16)(16, 12): -5 - that is larger than a previous path.
(2, 15) can't be extended to a non-positive entry using the values in row 15.
(2, 8) can't be extended to values to a non-positive entry using values in row 8.
Thus, there are no paths that can be further extended. Checking each of the paths, [a, .... , b] we've obtained, none of them have the property that [a, ... , b, a] is a non-positive cycle. Therefore, $\sigma_{FWTSPOPT} = D_7$. Furthermore, since we were unable to obtain any cycles in Step 3(a), from the corollary to Theorem 6, $\sigma_{TSPOPT} = \sigma_{FWTSPOPT}$.